\documentclass[preprint,12pt]{elsarticle-mod3}

\usepackage{pgf,tikz,circuitikz}
\usetikzlibrary{arrows}
\usepackage{graphicx}
\usepackage{amssymb}
\usepackage{amsthm}
\usepackage{amsmath}

\biboptions{comma,square}

\journal{J. Combin. Theory A 118:3 (2011), 920-937}

\usepackage[breaklinks,pdfstartview=FitH,colorlinks,linktocpage,bookmarks=true]{hyperref}
\usepackage{geometry}
\usepackage{pgf,tikz}\usetikzlibrary{arrows}

\newtheorem{theorem}{Theorem}[section]
\newtheorem{lemma}[theorem]{Lemma}
\newtheorem{conjecture}[theorem]{Conjecture}
\newtheorem{corollary}[theorem]{Corollary}
\newtheorem{claim}[theorem]{Claim}

\theoremstyle{definition}
\newtheorem{definition}[theorem]{Definition}
\newtheorem{example}[theorem]{Example}
\newtheorem{problem}[theorem]{Problem}

\theoremstyle{remark}
\newtheorem{remark}[theorem]{Remark}

\DeclareMathOperator{\Real}{Re}
\DeclareMathOperator{\Imaginary}{Im}
\DeclareMathOperator{\Interior}{Int}
\DeclareMathOperator{\Area}{Area}

\begin{document}

\begin{frontmatter}



\title{Tiling by rectangles and alternating current}


\author[add1]{M. Prasolov}
\ead{0x00002a@gmail.com}
\address[add1]{Moscow State University, Faculty of Mechanics and mathematics, \\
Leninskie Gory, 1, GSP-1, Moscow, 119991, Russian Federation}

\author[add2,add3]{M. Skopenkov\corref{cor1}}
\ead{skopenkov@rambler.ru} %, http://skopenkov.ru} %, Tel.: +966 530 708 319, Fax: +966 2 802 0064}
\cortext[cor1]{Corresponding author}
\address[add2]{Institute for information transmission problems of the Russian Academy of Sciences,\\
Bolshoy Karetny per.~19, bld.~1, Moscow, 127994, Russian Federation}
\address[add3]{King Abdullah University of Science and Technology,\\
P.O.~Box 2187, 4700 KAUST, 23955-6900 Thuwal, Kingdom of Saudi Arabia}

\begin{abstract}
This paper is on tilings of polygons by rectangles.
A celebrated physical interpretation of such tilings
by R.L.~Brooks, C.A.B.~Smith, A.H.~Stone and W.T.~Tutte
uses direct-current circuits.
The new approach of this paper is an application of alternating-current circuits.
The following results are obtained:
\begin{itemize}
\item a necessary condition for a rectangle %of given side ratio can be tiled
to be tilable by rectangles of given shapes;
\item a criterion for a rectangle to be tilable by rectangles similar to it but not all homothetic to it;
\item a criterion for a ``generic'' polygon to be tilable by squares.
\end{itemize}
These results generalize those of C.~Freiling, R.~Kenyon, M.~Laczkovich, D.~Rinne, and G.~Szekeres. %The approach is also applied to discrete harmonic functions and random walks on graphs.
\end{abstract}


\end{frontmatter}

\section{Introduction}\label{intro}

%This paper is on combinatorial geometry and its relationship with graph theory, electrical circuit theory, random walks.

A rectangle $a\times b$, where $a$ and $b$ are integers, can be tiled by $a\cdot b$ squares. Thus a rectangle with rational side ratio can be tiled by squares. In 1903 M.~Dehn proved the reciprocal assertion:
%(see figure~\ref{bookshelf}):

\begin{theorem} \label{dehn} \textup{\cite{D}} A rectangle can be tiled by squares \textup{(}not necessarily equal\textup{)} if and only if
the ratio of two orthogonal sides of the rectangle is rational.
%it has rational side ratio.
\end{theorem}

%\begin{figure}[htbp]
%\begin{tabular}{cc}
%\includegraphics[width=3.5cm]{bookshelf.pdf} & %\includegraphics[width=4.6cm]{squaredsquare.pdf}
%\end{tabular}
%\caption{Examples of tilings by squares from \url{http://www.squaring.net}}
%\label{bookshelf}
%\end{figure}

Although this assertion is expected, the proof is complicated. After the original proof, many improvements have been made \cite{B,BSST,Ha57,P,Ya68}.
%of Hadwiger \cite{Ha57}, Pokrovskii \cite{P}, Boltianskii \cite{B}, Yaglom \cite{Ya68}.

The most interesting for us is the approach of R.~L.~Brooks, C.~A.~B.~Smith, A.~H.~Stone, and W.~T.~Tutte~\cite{BSST}. To a tiling of a rectangle they assign a direct-current circuit, and then deduce Theorem~\ref{dehn} from certain properties of the circuit.
%the uniqueness of current distribution in the circuit.
%We present this proof in \S\ref{prelim} as an introduction to our approach.
They also apply this technique to find a tiling of a square by squares of distinct sizes \cite{Dui92}; see %a story of Gardner
%~\cite{G} and
%\url{http://www.squaring.net} for a survey and artwork.
the figure on the front cover of the journal. %(J. Combin. Theory A). 
An elementary introduction to the approach is given in \cite{PrSk,SMD}.

%Later the minimal tiling of this kind was found by a computer search \cite{Dui92}, see \url{http://www.squaring.net} for a survey and artwork.
%see figure~\ref{bookshelf} to the right.
%In \S\ref{prelim} we describe this physical interpretation in detail.
%and present the physical proof of Theorem \ref{dehn} as an illustration.

%\begin{figure}[htbp]
%\includegraphics[width=4.5cm]{squaredsquare.pdf}
%\caption{A square can be tiled by $21$ distinct squares. Number $21$ is minimal with this property.}
%\label{squaredsquare}
%\end{figure}

%A classical related question is whether a square can be tiled by distinct squares. The minimal example of such a tiling has been found quite recently by a computer search, see figure~\ref{squaredsquare}.

%\begin{theorem}[Duijstvijn, 2007] A square can be tiled by $21$ distinct squares. Number $21$ is minimal with this property.
%\end{theorem}

In this paper we study finite tilings by arbitrary
rectangles. The sides of rectangles are assumed to be parallel to the coordinate axes, i.~e., either vertical or horizontal. By {\it the ratio} of a rectangle we %always
mean the length of the horizontal side divided by the length of the vertical side. We study the following %\emph{shape-tiling}
problem posed in~\cite[p.~218]{FLR} and~\cite[p.~3]{KeKi}:

\begin{problem} \label{mainproblem} Which rectangles can be tiled by rectangles of given ratios $c_1$, \dots,~$c_n$?
%When a rectangle of ratio $c$ can be tiled by rectangles of given ratios $c_1$, $c_2$, \dots,~$c_n$?
\end{problem}

A related problem of \emph{signed} tilings is solved in~\cite{KeKi}.

We do not put any restrictions on the number of rectangles in the tilings. %In particular, this number may be distinct from $n$.
Each of the ratios $c_1$, \dots,~$c_n$ can be used any number of times or not used at all.

For $n=1$ and $c_1=1$ the question of Problem~\ref{mainproblem} is answered by Theorem~\ref{dehn}.
%Notice that for each ratio $c_k$ we allow several rectangles in the tiling have this ratio.
A necessary condition for arbitrary $n$ was actually proved by M.~Dehn:
%if a required tiling exists
if a rectangle of ratio $c$ can be tiled by rectangles of ratios $c_1$, \dots,~$c_n$
then $c$ is (the value of) a rational function in $c_1$, $\dots$, $c_n$ with rational coefficients \cite[Lemma~4]{FR}.

This function depends only on the ``combinatorial structure'' of the tiling. For instance, if a rectangle of ratio $c$ is dissected into $2$ rectangles of ratios $c_1$ and $c_2$ by a vertical (respectively, horizontal) cut then $c(c_1,c_2)=c_1+c_2$ (respectively, $c(c_1,c_2)=\frac{c_1c_2}{c_1+c_2}$).
The problem reduces to description of possible functions $c(c_1,\dots,c_n)$. By the mentioned physical interpretation this is equivalent to a natural problem:
%\begin{problem}
\textit{describe possible formulas $c(c_1, \dots, c_n)$ expressing the conductance of a planar direct-current circuit through the conductances $c_1$, $\dots$, $c_n$ of individual resistors}.
%\end{problem}

The main idea of the paper is to apply {\it alternating}-current circuits (equivalently, circuits with complex-valued conductances) to the above problems. Our first result is the following theorem.

\begin{theorem}\label{th2}
Suppose that a rectangle of ratio $c$ can be tiled by rectangles of ratios $c_1$,  \dots, $c_n$.
%(in such a way that there is at least one rectangle with each ratio $R_k$).
Then $c=C(c_1,  \dots, c_n)$ for some rational function $C(z_1,\dots, z_n)$ having the following properties
\begin{enumerate}[(1)]
%\item $C(c_1, c_2, \dots, c_m)=c$;
\item \label{th2-1} \emph{rationality of coefficients}\textup{:}
    %$C(z_1, \dots, z_n)$ has rational coefficients, i.e.,
    $C(z_1, \dots, z_n)\in \mathbb{Q}(z_1,\dots,z_n)$\textup{;}
    %i.e., $C$ is a rational function with rational coefficients;
\item \label{th2-2} \emph{homogeneity}\textup{:}
    %$C(z_1, \dots, z_n)$ is degree $1$ homogeneous, i.e.,
    $C(tz_1, \dots, tz_n)= tC(z_1,\dots, z_n)$\textup{;}
\item \label{th2-3} \emph{positive reality}\textup{:} if $\Real z_1,\dots,\Real z_n>0$ then $\Real C(z_1, \dots, z_n)>0$.
\end{enumerate}
\end{theorem}

\begin{problem} Is the reciprocal theorem true for $n\ge 3$?
\end{problem}

%We prove this theorem in \S\ref{proofs}.
%Parts (\ref{th2-1}) and~(\ref{th2-2}) of Theorem~\ref{th2} were actually proved by M.~Dehn, see also~\cite[Lemma~4]{FR}.
Case $n=1$ (respectively, $n=2$) of both Theorem~\ref{th2} and its reciprocal %are known: these cases
is equivalent to Theorem~\ref{dehn} (respectively, to \cite[Theorem~5]{FLR} or else to Corollary~\ref{series-parallel} below).
%There is another condition, which the function $C$ must satisfy, but we do not know whether this condition follows from (1)--(3) for $n\ge 3$. \textbf{Dopisat'.}
%There is a heuristic argument showing that
For $n\ge 3$ the reciprocal theorem cannot be proved by our method; see Example~\ref{heuristic}.

Theorem~\ref{th2} has a clear physical meaning; see \S\ref{alternating}. But this theorem
(even together with its reciprocal) is not \emph{algorithmic}, i.~e., it does not give an algorithm to decide, if there exists a required tiling. Thus it is interesting to get less general but algorithmic results.

%The first
A result of this kind %after the Dehn theorem
was obtained independently by C.~Freiling, D.~Rinne in 1994 and M.~Laczkovich, G.~Szekeres in 1995. %It uses the following notion. An \emph{algebraic conjugate} of an algebraic number $c$ is a complex root of the minimal integral polynomial of $c$.

\begin{theorem} \label{lfrs} \textup{\cite{FR, LS}}
For $c>0$ the following $3$ conditions are equivalent\textup{:}
\begin{enumerate}[(1)]
\item \label{lfrs-1} a square can be tiled by rectangles of ratios $c$ and $1/c$\textup{;}
\item \label{lfrs-2} %the number $c$ is algebraic and all its algebraic conjugates have positive real parts
the number $c$ is a root of an integral polynomial such that all the other complex roots have positive real parts\textup{;}
\item \label{lfrs-3} for certain positive rational numbers $d_1$, \dots, $d_m$ we have
$$d_1 c +\cfrac1{d_2 c+\cfrac1{d_3 c+\dots+\cfrac1 {d_m c}}}=1.$$
\end{enumerate}
\end{theorem}

In particular, a square can be tiled by rectangles of ratios $2+\sqrt{2}$ and $\frac{1}{2+\sqrt{2}}$ but cannot be tiled by rectangles of ratios $1+\sqrt{2}$ and $\frac{1}{1+\sqrt{2}}$; see \cite{SMD} for an elementary proof.

We present a new short self-contained proof of Theorem~\ref{lfrs}.
This new proof %(announced in \cite{PrSk})
is a natural application of alternating-current circuits. We also get a new algorithmic result: % which does not seem to be tractable by methods available before:
%of \cite{FLR, FR, LS}:

\begin{theorem}\label{th1} %(cf.~\cite{FR,LS})
For a number $c>0$ the following $3$ conditions are equivalent:
\begin{enumerate}[(1)]
\item \label{th1-1} a rectangle of ratio $c$ can be tiled by rectangles of ratios $c$ and $1/c$
    so that there is at least one rectangle of ratio $1/c$ in the tiling\textup{;}
\item \label{th1-2} %the number $c^2$ is algebraic and all its algebraic conjugates distinct from $c^2$ are negative real numbers
the number $c^2$ is a root of an integral polynomial such that all the other complex roots are negative real numbers\textup{;}
\item \label{th1-3} for certain positive rational numbers $d_1$, \dots, $d_m$ we have
$$\cfrac{1}{d_1 c +\cfrac1{d_2 c+\dots+\cfrac1 {d_m c}}}=c.$$
\end{enumerate}
\end{theorem}

%This is proved in \S\ref{proofs}.
More algorithmic results can be found in~\cite[p.~224]{FLR} and \cite{Sharov}. For similar results on tiling by triangles see \cite{Szegedy}. For higher-dimensional generalizations %of tilings by rectangles
see \cite{P}. %\cite{BB,P}.
%The latter shows in particular that the case when not all ratios $c_n$ are ''used'' at least once in the tiling should be considered separately, thus giving a motivation to the statement of Theorem~\ref{th1}.

%\medskip
%\subsection{Tiling of polygons and electrical circuits with many terminals}
%In \S\ref{var}
We also consider tilings of arbitrary (not necessarily convex) polygons by rectangles. This generalization
%might seem a routine but in fact it is substantial and
reveals new connections between tilings and electrical circuits.

We apply direct-current circuits with several terminals to get a criterion for a ``generic'' polygon to be tilable by squares; see Theorem~\ref{squaring-polygons} below, again not algorithmic. This result generalizes Theorem~\ref{dehn} and~\cite[Theorems 9 and 12]{K}. An easier related problem of \emph{signed} tiling by squares is solved in \cite{FHTW,KK}.  %, see figure~\ref{chip}.

%\begin{figure}[htbp]
%\includegraphics[width=3.5cm]{chip.pdf}
%\caption{An electrical circuit with many terminals}
%\label{chip}
%\end{figure}

%Our approach to tilings of polygons uses direct-current circuits with several terminals.
%, studied by Colin de Verdi\`ere, Lawler, Sylvester and others
%\emph{Direct} and \emph{inverse problems} for such circuits is the subject of \emph{electrical impedance tomography} \cite{C, CGV, LSy}.

We apply alternating-current circuits with several terminals to get a short proof of a generalization of Theorem~\ref{lfrs} to polygons with rational vertices \cite{SuDi}; see Theorem~\ref{polygon-square} below.
%The proofs are based on certain known results of \emph{electrical impedance tomography}.
We also state a basic folklore result on \emph{electrical impedance tomography} for alternating-current circuits, cf.~\cite{CM, C, CGV, LSy}.

There is a close relationship among electrical circuits, discrete harmonic functions, and
%theory of discrete harmonic functions and
random walks on graphs \cite{DS, Lovasz, BeSch}. %, see a brilliant introduction in \cite{DS} and a survey in \cite{Lovasz}, see also \cite{BeSch}.
Our results have equivalent statements in the language of each of the theories; e.~g., see Corollary~\ref{random} below.
%either the language of electrical circuits or discrete harmonic functions or random walks.
%The reader can find some applications of our results to these subjects in the paper.
%We prefer the language of electrical circuits only because of physical clearness.

%\subsection{Organization of the paper} In \S\ref{prelim} we recall the basics of electrical circuits and present main ides of our proofs. In \S\ref{proofs} we prove Theorems~\ref{th2}, ~\ref{th1} and~\ref{lfrs}. In \S\ref{vars} we discuss some applications of our approach to electrical circuits with many terminals and to random walks.

%The paper is organized as follows. It splits naturally into two formally independent parts.

From here the paper splits into two formally independent parts: \S\S\ref{prelim}--\ref{proofs} and \S\S\ref{var}--\ref{finalproofs}.

The first part contains the proof of Theorems~\ref{th2}, \ref{lfrs}, \ref{th1}. In \S\ref{prelim} the basics of electrical circuits and their connection with tilings are recalled. In \S\ref{proofs} the results of~\S\ref{intro} are proved.

The second part concerns some variations. In \S\ref{var} the results on
tilings of polygons, electrical impedance tomography, and random walks are stated.
In \S\ref{final} the results of \S\ref{prelim} are generalized to electrical circuits with several terminals. In \S\ref{finalproofs}  the results of~\S\ref{var} are proved.

\section{Main ideas}\label{prelim}

\subsection{Electrical circuits}

Our approach is based on electrical circuits theory \cite{PrSk}. However, the reader is not assumed to be familiar with physics. In this section we recall all the required physical concepts  (although the presentation is formal and physical meaning is explained very briefly). This section does not contain new results.
For short proofs see~\S\ref{final}.
%All the results of \S\ref{prelim} are proved in \S\ref{final} in a more general situation.

%Informally, \textit{an electrical network} is a collection of conductors connected with each other. A fixed pair of points of the network can be joined with the poles of the battery.

%Formally,
An \textit{electrical network} is a connected graph
%(possibly having multiple edges)
with a nonnegative real number (\textit{conductance}) assigned to each edge, and two marked (\textit{boundary}) vertices.

% and two real numbers $U_1$ and $U_2$ assigned to each boundary
% vertex, called \textit{boundary voltage}.

For simplicity in this subsection we assume that the graph has neither multiple edges nor loops. Generalizations
for graphs with multiple edges are left to the reader. We say that an electrical network is \textit{planar}, if the graph is drawn in the unit disc in such a way that the boundary vertices belong to the boundary of the disc and the edges do not intersect each other.

%Denote by $V$ and $E$ the sets of vertices and
%edges, respectively, of the graph in question.
%Fix an orientation of each edge of the graph.

Fix an enumeration of the vertices $1$, $2$, $\dots$, $n$ such that $1$ and $2$ are the boundary ones. It is convenient to denote the number of boundary vertices by $b=2$. Denote by $m$ the number of edges. Denote by $c_{kl}$ the conductance of the edge between the vertices $k$ and $l$. Set $c_{kl}=0$, if there is no edge between $k$ and $l$.
%In particular, all $c_{kk}=0$.
%It is convenient to put also $c_{kk}=-\sum_{l=1}^n c_{kl}$.
%We often use shorthand \textit{conductor $c$} for the phrase \textit{an edge of conductance $c$}.

%Informally, \textit{an electrical circuit} is an electrical network whose boundary vertices are joint to a battery of certain specific voltage.

%Formally, by
An \textit{electrical circuit} is an electrical network along with two real numbers $U_1$ and $U_2$ (\textit{incoming voltages}) assigned to the boundary vertices.

Each electrical circuit gives rise to certain numbers $U_k$, where $1\le k\le n$ (\textit{voltages} at the vertices), and $I_{kl}$, where $1\le k,l\le n$ (\textit{currents} through the edges).
These numbers are defined by the following $2$ axioms:
%(0) \textit{Boundary condition.}
%$U(v_1)=U_1$, $U(v_n)=U_n$.
\begin{enumerate}[(1)]
\item[(C)]\textit{The Ohm law.} For each pair of vertices $k,l$ we have $I_{kl}=c_{kl}(U_k-U_l)$.
\item[(I)]\textit{The Kirchhoff current law.} For each vertex $k>b$ we have $\sum_{l=1}^{n} I_{kl}=0$.
\end{enumerate}
Informally law (I) means that electrical charge is not aggregated at the nonboundary vertices. In other words, these laws assert that $U_k$ is a \emph{discrete harmonic function}.
The numbers $U_k$ and $I_{kl}$ are well-defined by these axioms by the following classical result.
% proved in the end of \S\ref{prelim}.

\begin{theorem} \label{weyl} \textup{\cite{Weyl}} For any electrical circuit
the system of linear equations \textup{(C), (I)} in variables $U_k$, where $b< k\le n$, and $I_{kl}$, where $1\le k,l\le n$, has a unique solution.
\end{theorem}

%We prove a generalization of this theorem in \S\ref{proofs}.

\begin{figure}[htbp]
\begin{tabular}{cc}
\includegraphics[width=5cm]{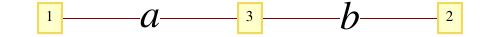} & \includegraphics[width=5cm]{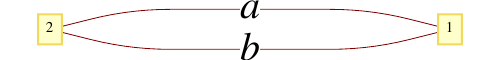}\\[9pt]
$C(a,b)=\frac{ab}{a+b}$ & $C(a,b)=a+b$
\end{tabular}
\caption{Series and parallel electrical networks}
\label{examples}
\end{figure}

Denote by $I_1=\sum_{k=1}^{n} I_{1k}$ the current flowing inside the circuit through vertex $1$.
\textit{The conductance} of an electrical circuit with $U_1\ne U_2$ is the number $C=I_1/(U_1-U_2)$. It is easy to see that the conductance does not depend on $U_1$ and $U_2$.
Thus \textit{the conductance} of an electrical network is well-defined. The reciprocal of conductance is called \emph{resistance}.
%Thus one can define \textit{the conductance} of an electrical network to be the conductance of any circuit made of the network.
Basic examples of networks and their conductances are shown in figure~\ref{examples}.

\subsection{Circuits and tilings}\label{ssec2-2}

%Let us explain informally the relationship between tilings by rectangles and electrical networks. Suppose that a rectangle is tiled by rectangles. \textbf{Plastinka. This has not been written yet.}

There is a close relationship between electrical circuits and tilings.
We say that an edge $kl$ of a circuit is \emph{essential}, if %for some (and thus for all) $U_1\ne U_2$ we have
$I_{kl}\ne 0$. Clearly, this property of an edge does not depend on $U_1$ and $U_2$, if $U_1\ne U_2$. Thus \emph{essential} edges in an electrical network are well-defined.

\begin{lemma} \label{correspondence-lemma} \textup{\cite[Theorem 1.4.1]{BSST,CFP}} %(see figure~\ref{tilingnetworks})
The following $2$ conditions are equivalent\textup{:}
\begin{enumerate}[(1)]
\item \label{c-l-1} a rectangle of ratio $c$ can be tiled by $m$ rectangles of ratios $c_1$, $\dots$, $c_m$\textup{;}
\item \label{c-l-2} there is a planar electrical network having conductance $c$ and consisting of $m$ essential edges of conductances
$c_1$, $\dots$, $c_m$.
\end{enumerate}
\end{lemma}

\begin{figure}[htbp]
\definecolor{qqwwff}{rgb}{0,0.4,1}
\definecolor{ffqqqq}{rgb}{1,0,0}
\definecolor{zzttqq}{rgb}{0.6,0.2,0}
\begin{tikzpicture}[line cap=round,line join=round,>=triangle 45,x=0.5cm,y=0.5cm]
\clip(-3.14,-4.24) rectangle (7.08,6.78);
\fill[color=zzttqq,fill=zzttqq,fill opacity=0.1] (-3,6) -- (6.9,6) -- (6.9,-3.6) -- (-3,-3.6) -- cycle;
\draw [color=zzttqq] (-3,6)-- (6.9,6);
\draw [color=zzttqq] (6.9,6)-- (6.9,-3.6);
\draw [color=zzttqq] (6.9,-3.6)-- (-3,-3.6);
\draw [color=zzttqq] (-3,-3.6)-- (-3,6);
\draw [color=zzttqq] (2.4,6)-- (2.4,0.6);
\draw [color=zzttqq] (2.4,0.6)-- (-3,0.6);
\draw [color=zzttqq] (6.9,1.5)-- (2.4,1.5);
\draw [color=zzttqq] (1.2,-3.6)-- (1.2,0.6);
\draw [color=zzttqq] (1.2,-0.6)-- (4.5,-0.6);
\draw (2.4,0.6)-- (2.4,-0.6);
\draw [color=zzttqq] (4.5,1.5)-- (4.5,-0.9);
\draw [color=zzttqq] (6.9,-0.9)-- (4.2,-0.9);
\draw [color=zzttqq] (4.2,-0.6)-- (4.2,-3.6);
\draw [line width=1pt,color=qqwwff] (1.98,6)-- (-0.3,0.6);
\draw [line width=1pt,color=qqwwff] (1.98,6)-- (4.66,1.5);
\draw [line width=1pt,color=qqwwff] (-0.3,0.6)-- (1.94,-3.6);
\draw [line width=1pt,color=qqwwff] (2.7,-0.6)-- (4.66,1.5);
\draw [line width=1pt,color=qqwwff] (-0.3,0.6)-- (2.7,-0.6);
\draw [line width=1pt,color=qqwwff] (5.56,-0.9)-- (4.66,1.5);
\draw [line width=1pt,color=qqwwff] (5.56,-0.9)-- (1.94,-3.6);
\draw [line width=1pt,color=qqwwff] (2.7,-0.6)-- (1.94,-3.6);
\draw [line width=1pt,color=qqwwff] (5.56,-0.9)-- (2.7,-0.6);
\fill [color=ffqqqq] (1.98,6) circle (2.0pt);
\draw[color=ffqqqq] (1.98,6.5) node {$+$};
\fill [color=ffqqqq] (-0.3,0.6) circle (2.0pt);
\fill [color=ffqqqq] (4.66,1.5) circle (2.0pt);
\fill [color=ffqqqq] (1.94,-3.6) circle (2.0pt);
\draw[color=ffqqqq] (1.96,-3.95) node {$-$};
\fill [color=ffqqqq] (2.7,-0.6) circle (2.0pt);
\fill [color=ffqqqq] (5.56,-0.9) circle (2.0pt);
\end{tikzpicture} %$\!\!\!\!$\input{metal-fig6c.tex}
\caption{Correspondence between tilings and electrical networks.}
\label{tilingnetworks}
\end{figure}
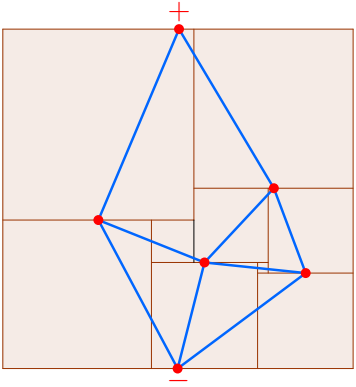

%Only assertion (1)$\implies$(2) is required for the proof of main results.
Let us sketch the proof of assertion (\ref{c-l-1})$\implies$(\ref{c-l-2}).
%, the only part of the lemma required for the proof of main results
Given a tiling as in~(\ref{c-l-1}) construct an electrical network as follows; see Figure~\ref{tilingnetworks}. Take a point in each horizontal cut of the tiling and in each horizontal side of the tiled rectangle. These points are vertices of the network. For each rectangle in the tiling draw an edge between the vertices in the cuts containing the horizontal sides of the rectangle. Set the conductance of the edge to be the ratio of the rectangle.
%Put the current through the edge equal to the horizontal side of the rectangle.
%Put the voltage at a vertex equal to the $y$-coordinate of the vertex.
%then laws (C),(I) are satisfied and
Then the conductance of the resulting network equals to the ratio of the tiled rectangle; see \S\ref{correspond} for the details.
%This network has conductance $c$.

%We prove Lemma~\ref{correspondence-lemma} in the end of this section.

%\textbf{Proof of lemma. Not written yet!!!}

%Formally, \textit{a circuit corresponding to a tiling} is defined as follows. \textbf{Correspondence. %This has not been written yet.}

\subsection{Formulas for conductance}

Let us summarize some useful properties of formulas for conductance.

\begin{lemma} \label{electrical} \label{main}
\textup{\cite{BSST,DS,Foster}}
Suppose that an electrical network consists of $m$ edges of conductances $c_{1}, \dots, c_{m}$.
Then the conductance $C(c_{1}, \dots, c_{m})$ of the network has the following properties\textup{:}
\begin{enumerate}[(1)]
\item \label{el-1} \emph{rationality}\textup{:} $C(c_{1}, \dots, c_{m})\in \mathbb{Q}(c_{1}, \dots, c_{m})$\textup{;}
\item \label{el-2} \emph{homogeneity}\textup{:} $C(tc_{1}, \dots, tc_{m})=tC(c_{1}, \dots, c_{m})$\textup{;}
\item \label{el-3} $\frac{\partial}{\partial c_j} C(c_{1}, \dots, c_{m})=(\frac{U_k-U_l}{U_1-U_2})^2$,
    where $k$ and $l$ are the endpoints of the edge $j$\textup{;}
\item \label{el-4} \emph{monotonicity}\textup{:} if $c_{1}, \dots, c_{m}>0$ then $\frac{\partial}{\partial c_j} C(c_{1}, \dots, c_{m} )\ge 0$; if the edge $j$ is essential then the latter inequality is strict\textup{;}
\item \label{el-5} \emph{positive reality}\textup{:} if $\Real c_{1},\dots, \Real c_{m}>0$ then $\Real C(c_{1}, \dots, c_{m})>0$.
\end{enumerate}
\end{lemma}

%\begin{problem} Is the converse theorem true?
%\end{problem}

%\begin{remark} (Akopyan, private communication) $C(c_{1}, \dots, c_{m})$ belongs to a convex cone with vertex at the origin, spanned by the points $c_{1}, \dots, c_{m}$ in $\mathbb{C}$.
%\end{remark}

\begin{remark} (A. Akopyan, private communication) Property (\ref{el-4}) follows from (\ref{el-1}), (\ref{el-2}), and~(\ref{el-5}).
Property~(\ref{el-5}) does not follow from~(\ref{el-1}), (\ref{el-2}), and~(\ref{el-4}); e.~g., the function $C(c_1,c_2)=(c_1+c_2)\frac{c_1^2+c_2^2}{c_1^2+2c_2^2}$ satisfies~(\ref{el-1}), (\ref{el-2}), (\ref{el-4}), but not~(\ref{el-5}).
\end{remark}

%Theorem \ref{dehn} follows directly from Lemma \ref{electrical}(1).
%Properties (\ref{el-1})--(\ref{el-4}) are well-known,

Property~(\ref{el-5}) concerns the extension of the function $C(c_1,\dots,c_m)$ to a complex domain.
%is more tricky because it involves an extension of $C(c_{1}, c_{2}, \dots, c_{m})$ to the complex plane. The latter result
%It does not seem that anyone has paid attention to this fundamental property for \textit{direct}-current circuits. Certainly it is well-known for
The meaning of this property becomes clear in the context of \textit{alternating}-current circuits. Short proof of the lemma is given in \S\ref{elcirc}.

\subsection{Alternating current}\label{alternating}

Let us explain the informal physical meaning of
the ``positive reality'' condition in
Lemma~\ref{main} and Theorem~\ref{th2}. This is not used elsewhere in the paper and the reader may easily skip this subsection. % and proceed to the formal proofs in the next one.

%It is quite natural to study a rational function $C(c_{ij})$ on the whole complex plane. A more surprising observation is that the function $C(c_{kl})$ for complex $c_{kl}$ has a clear physical meaning. Let us introduce the basics of \textit{alternating-current circuits}, which are logically equivalent to electrical circuits with complex-valued conductances.

%Electrical circuits with complex-valued conductances are ''equivalent'' toalternating-current circuits.

Informally, \textit{an alternating-current circuit} is a collection of conductors, capacitors, inductors and a single alternating-voltage source connected with each other.

%\textbf{Give the following algorithm applied to a specific example only?}

Formally, {\it an alternating-current circuit} is a %connected
graph with the following %additional
structure:
\begin{itemize}
\item two marked (\emph{boundary}) vertices;
\item two functions (\emph{voltages}) $\tilde U_1(t)=U\cos \omega t$ and $\tilde U_2(t)=0$ assigned to them; %;
\item division of the edges into $3$ types (\textit{conductors, capacitors,  and inductors});
\item a positive number $\tilde c_{kl}$ assigned to each edge (called \textit{conductance, capacitance}, or \textit{inductance}, depending on the type of edge).
\end{itemize}
Assume for simplicity that the graph has no multiple edges. The \emph{voltages} $\tilde U_k(t)$ and the \emph{currents} $\tilde I_{kl}(t)$ are  defined by the following axioms:
\begin{enumerate}[(1)]
\item[(\~ C)]\textit{The generalized Ohm law.} For each edge $kl$ we have
$$\tilde I_{kl}(t)=
\begin{cases}\tilde c_{kl}\left(\tilde U_k(t)-\tilde U_l(t)\right),& \text{if $kl$ is a conductor;}\\
\tilde c_{kl}\frac{d}{dt}\left(\tilde U_k(t)-\tilde U_l(t)\right),& \text{if $kl$ is a capacitor;}\\
\tilde c_{kl}\int_{\pi/2\omega}^t\left(\tilde U_k(t)-\tilde U_l(t)\right)dt,& \text{if $kl$ is an inductor.}
\end{cases}
$$
\item[(\~ I)]\textit{The Kirchhoff current law.} For each vertex $k\ne 1,2$ we have $\sum_{l=1}^{n} \tilde I_{kl}(t)=0$.
\end{enumerate}

The voltages and the currents can be found using the following well-known algorithm. Denote by $i=\sqrt{-1}$. Put $U_1=U$, $U_2=0$, and
$$c_{kl}=
\begin{cases}\quad\tilde c_{kl},  & \text{if $kl$ is a conductor;}\\
i\omega\tilde c_{kl},             & \text{if $kl$ is a capacitor;}\\
\frac{1}{i\omega}\tilde c_{kl},   & \text{if $kl$ is an inductor.}
\end{cases}
$$
Define the complex numbers $U_k$, where $3\le k\le n$, and $I_{kl}$, where $1\le k,l\le n$, by \emph{direct}-current laws (C), (I).
%Then the real parts of these numbers are the voltages and the currents in the circuit at the moment $t$. (Notice that the arguments of the complex numbers $U_k$ and $I_{kl}$ depend on the time parameter $t$.)
Then $\tilde U_k(t)=\Real (U_k e^{i\omega t})$, $\tilde I_{kl}(t)=\Real (I_{kl}e^{i\omega t})$.
In this sense alternating-current circuits are ``equivalent'' to direct-current circuits with complex-valued conductances (also called \emph{admittances}).

%\textbf{Example. Not written yet!!!}

%Let us explain physical meaning of Lemma~\ref{electrical}(5) and Theorem~\ref{th2}(3).
Notice that always $\Real c_{kl}\ge 0$. Physically this means nonnegative
\textit{energy dissipation} at the edge $kl$ (which equals to $\Real  c_{kl}|U_k-U_l|^2$).
%\begin{equation*}
%\int_{0}^{2\pi/\omega}\tilde I_{kl}(t)(\tilde U_k(t)-\tilde U_l(t))dt=\Real  (\bar I_{kl}(U_k-U_l))=\Real \bar c_{kl}|U_k-U_l|^2.
%\end{equation*}
%Thus a physical meaning of positive real part of \emph{admittance} $c_{kl}$ is positive energy dissipation.
Thus a physical meaning of positive reality is: ``a network consisting of elements dissipating energy also dissipates energy''.
%We make this intuitive argument precise in \S\ref{final}.

%Thus we come to the following formal definition.

\begin{figure}[h]
\hspace{4cm}
\makebox[1.5cm][r]{
  \begin{tabular}{l}
  \definecolor{zzttqq}{rgb}{0.6,0.2,0}
\begin{tikzpicture}[line cap=round,line join=round,>=triangle 45,x=0.3cm,y=0.3cm]
\clip(-3.42,-2.1) rectangle (4.32,5.92);
\fill[color=zzttqq,fill=zzttqq,fill opacity=0.1] (-2,5) -- (4,5) -- (4,-1) -- (-2,-1) -- cycle;
\draw [color=zzttqq] (-2,5)-- (4,5);
\draw [color=zzttqq] (4,5)-- (4,-1);
\draw [color=zzttqq] (4,-1)-- (-2,-1);
\draw [color=zzttqq] (-2,-1)-- (-2,5);
\draw [color=zzttqq] (-2,3.58)-- (4,3.56);
\draw [color=zzttqq] (-0.02,3.57)-- (0,-1);
\draw [color=zzttqq] (2,3.57)-- (2,-1);
\draw [color=zzttqq] (1,5)-- (1,3.57);
\end{tikzpicture}
  \end{tabular}
  }
  \makebox[3.5cm][c]{
  \begin{tabular}{c}
  %\documentclass{article}
%\usepackage{circuitikz}
%\begin{document}
\ctikzset{bipoles/length=.6cm}
\begin{circuitikz}[scale=0.6]
\draw[color=brown]
(0,0) to[short, *-*] (1,0)
      to[short, *-*] (2,0)
      to[short, *-] (4,0)
      to[sV](4,3)
      to[short, -*] (2,3)
      to[short, *-*] (0,3)
      to[C, *-*]     (0,2)
(0,0) to[L, *-*]     (0,2)
(1,0) to[L, *-*]     (1,2)
      to[short, *-*] (0,2)
(1,2) to[short, *-*] (2,2)        
(2,0) to[L, *-*]     (2,2)
      to[C, *-*]     (2,3)
;\end{circuitikz}
%\end{document}
  \end{tabular}
  }
\caption{Correspondence between tilings by rectangles and alternating-current networks.
}
\label{fig-perem}
\end{figure}
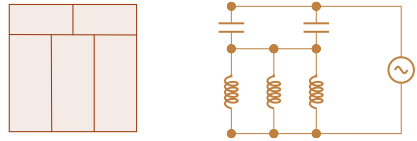

Furthermore, to each tiling of a rectangle %(or even a polygon) 
by \emph{similar} rectangles one can assign an alternating-current network; see Figure~\ref{fig-perem} (formally this is not used in the paper). 
The graph of the network is the same as in Section~\ref{ssec2-2}. But instead of resistors, the edges are now either unit capacitors or unit inductors, depending on the orientation of the corresponding rectangles. The constructed alternating-current network is ``equivalent'' to a direct-current network with complex edge conductances $z$ and $1/z$, where $z=i\omega$.

\subsection{Inverse problems}\label{inverse-problems}

%The proofs are based on the above physical interpretation.
%Let us prepare for the proof of other main results.
%This subsection is required only

\emph{The inverse problem} for electrical networks is to synthesize a network with given conductance from given elements.
%Let us state a simple inverse problem and outline its solution.
Let us state some classical results on this subject due to R.~M.~Foster and W.~Cauer. Short proofs are given in~\S\ref{basicposreal}.
This subsection is required only for the proof of assertions~(2)$\implies$(3) in Theorems~\ref{lfrs}--\ref{th1}.

\begin{theorem}[Foster's reactance theorem] \label{th-foster} \textup{\cite{Foster}} The following properties of a rational function $C(z)\in\mathbb{R}(z)$ are equivalent\textup{:}
\begin{enumerate}[(1)]
\item \label{th-foster-1}
    $C(i\omega)$ is the admittance of a network consisting of capacitors and inductors, considered as a function in the frequency $\omega$\textup{;}
\item \label{th-foster-2}
    $C(z)$ is the conductance of an electrical network such that the conductance of each edge $j$ is either $d_jz$ or $1/d_jz$ for some real numbers $d_1,\dots,d_m>0$\textup{;}
\item \label{th-foster-3} $C(z)$ is an odd \emph{positive real} rational function, i.~e., $\Real C(z)>0$, if $\Real z>0$.
\end{enumerate}
\end{theorem}

%\begin{problem} Construct a network from capacitors and inductors such that the conductance of the network is a given function of the frequency $\omega$.
%Equivalently, construct an electrical network such that the conductance of each edge $j$ is either $d_jz$ or $1/d_jz$ for some real numbers $d_j>0$ and the conductance of the network is a given function $C(z)$.
%\end{problem}

%Let us outline the proof of assertion~(\ref{th-foster-3})$\implies$(\ref{th-foster-2})  given in \cite{Cauer}.
%Let us outline the solution \cite{Foster, Cauer}. By Lemma~\ref{electrical} it follows that $C(z)$ must be an odd rational function, which is \emph{positive real}, i.~e., $\Real C(z)>0$, if $\Real z>0$.
%Using this property we construct a series-parallel electrical network with the conductance $C(z)$.

%The proof of assertions~(2)$\implies$(3) in Theorems~\ref{lfrs} and~\ref{th1} is based on the following analytical result.

The proof is based on the following analytic lemma.

\begin{lemma} \label{PosReal} \textup{\cite{Cauer}} For an odd rational function $C(z)\in \mathbb{R}(z)$ such that $\lim_{z\to \infty}C(z)\ne0$
the following $5$ conditions are equivalent\textup{:}
\begin{enumerate}[(1)]
\item\label{preserve} if $\Real z>0$ then $\Real C(z)>0$\textup{;}
%\item\label{preserve} $\Real C(z)>0$ if and only if $\Real z>0$;
\item\label{roots} if $C(z)=1$ then $\Real z>0$\textup{;}
%for some $x>0$ all roots of the equation $C(z)=x$ have positive real part;
\item\label{derivatives} if $C(z)=0$  then $\Real z=0$ and $C'(z)>0$\textup{;}
    %\textup{(}here $z\in\mathbb{C}$ or $z=\infty$\textup{)};
\item\label{factorization} %either $C(z)$ or $1/C(z)$ equals
$$
C(z)=d_1z\prod\limits_{k=1}^n\frac{z^2+a_k^2}{z^2+b_k^2},
$$
for some integer $n\ge 0$ and real numbers $d_1>0$, $a_1>b_1>a_2>\dots>b_n\ge 0$\textup{;}
%$$\frac{\prod\limits_{k=1}^n
%(z^2+b_k^2)}{d_1z\prod\limits_{k=1}^m(z^2+a_k^2)},
%$$
%for some integer numbers $m,n\ge 0$ and real numbers $a_k, b_k, d_1>0$ such that
%$m$ equals either $n{-}1$ or $n$ and
%$0<b_1<a_1<b_2<a_2<\dots$;
\item\label{continued} %either $C(z)$ or $1/C(z)$ equals
$$
C(z)=d_1 z +\cfrac1{d_2 z+\dots+\cfrac1 {d_m z}},
$$
%$$\cfrac1{d_1 z +\cfrac1{d_2 z+\dots+\cfrac1 {d_l z}}},$$
for some integer $m\ge1$ and real numbers $d_1,\dots,d_m>0$\textup{;}
%\item\label{series-parallel-network}
%$C(z)$ is the conductance of a series-parallel network with edges of conductances $d_1z$, $1/d_2z$, $d_3z$, $1/d_4z$, \dots, $(d_mz)^{(-1)^m}$ for some integer $m\ge1$ and real numbers $d_1,\dots,d_m>0$.
\end{enumerate}
\end{lemma}

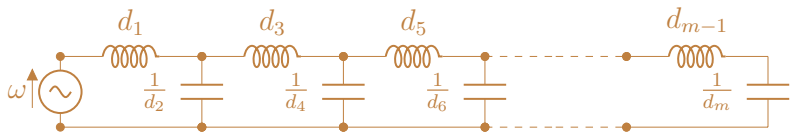
\begin{figure}[htb]
%\documentclass{article}
%\usepackage{circuitikz}
%\begin{document}
\ctikzset{bipoles/length=1.0cm}
\begin{circuitikz}[scale=1.0]
\draw[color=brown]
(0,1) to[L=$d_1$, *-*] (2,1)
      to[L=$d_3$, *-*] (4,1)
      to[L=$d_5$, *-*] (6,1)
(8,1) to[L=$d_{m-1}$, *-] (10,1)
(10,0) to[C=$\frac{1}{d_{m}}$] (10,1)
(10,0)to[short,-*] (8,0) 
(6,0)	to[short,*-*] (4,0)
(4,0)	to[short,*-*] (2,0)
(2,0)	to[short,*-*] (0,0)		        
(0,0)to[sV=$\omega$](0,1)
(2,0)to[C=$\frac{1}{d_{2}}$, *-] (2,1)
(4,0)to[C=$\frac{1}{d_{4}}$, *-] (4,1)
(6,0)to[C=$\frac{1}{d_{6}}$, *-] (6,1)
%(8,0)to[C=$\frac{1}{d_{8}}$, *-] (8,1)
;
\draw[dashed,color=brown]
(6,0) to[short, *-*] (8,0)      
(6,1) to[short, *-*] (8,1)  
;\end{circuitikz}
%\end{document}
\caption{A series-parallel network beyond Lemma~\ref{PosReal}}
\label{fig-ladder}
\end{figure}

Define inductively \emph{a series-parallel} electrical network. By definition, a network consisting of a single edge is series-parallel. If $a$ and $b$ are two series-parallel networks then both their series and parallel ``unions'' (see Figure~\ref{examples}) are series-parallel. One can see that condition~(\ref{continued}) of the lemma allows us to construct a series-parallel network of conductance $C(z)$ with edges of conductances $d_1z$, $1/d_2z$, $d_3z$, $1/d_4z$, \dots, $(d_mz)^{(-1)^m}$; see Figure~\ref{fig-ladder}. Assertion~(\ref{roots})$\implies$(\ref{continued}) was proved in~\cite[Lemma 4]{FLR} using the results of~\cite{Wall}; our proof following \cite{Cauer} is simpler.

%A rational function satisfying condition~(\ref{preserve}) of Lemma~\ref{PosReal} is called \emph{positive real}.
%Implication (\ref{preserve})$\implies$(\ref{factorization}) is known as \emph{the Foster reactance theorem}~\cite{Foster}.
%The rest of the lemma is proved in \cite{Cauer}.

%Short proofs of Lemma~\ref{PosReal} and the following corollary are given in~\S\ref{basicposreal}.

%A physical meaning of Lemma~\ref{PosReal} is synthesis of a circuit consisting of capacitors and inductors such that the conductance of the circuit is a given odd positive real function $C(z)$ in the frequency $z=i\omega$.

%Consider electrical circuits, in which all the edges have conductances $z$ and $1/z$, $\Real z>0$. (They have a natural physical meaning: circuits consisting of capacitors and inductors with incoming voltage of complex frequency $z/i$.) Let us describe possible conductances $C(z)$ of such electrical circuits. By Lemma~\ref{electrical}(\ref{el-1}), (\ref{el-2}) and~(\ref{el-5}) the functions $C(z)$ are \emph{positive real}, i.e., satisfy condition~(\ref{preserve}) of the following lemma.
%Denote by $\Real \infty=0$, $C(\infty)=\lim_{z\to 0}C(1/z)$ and $C'(\infty)=\lim_{z\to 0}(C(1/z))'$.
%for certain complex number $z$.
%The conductance of such a network is a function $C(z)$ in one variable $z$. By Lemma~\ref{electrical}(1), (2) and (5)
%this function satisfies property~(2) of the following lemma.

%\subsection{A non-series-parallel network}

%The following result shows that Theorem~\ref{th2} cannot be %proved by our method for $n\ge 3$.

\begin{corollary} \label{series-parallel} \textup{\cite{Cauer}}
If a function $C(c_1,c_2)$ satisfies conditions~\textup{(\ref{th2-1})--(\ref{th2-3})} of Theorem~\ref{th2} then
$C(c_1,c_2)$ is the conductance of a series-parallel electrical network with edge conductances $c_1$ and $c_2$.
\end{corollary}

\begin{example} \label{heuristic} A generalization of Corollary~\ref{series-parallel} to the case of $3$ variables $c_1$, $c_2$, $c_3$ is not true.
E.~g., take a network with $4$ vertices and edge conductances $c_{13}=c_1$, $c_{23}=c_2$, $c_{24}=c_1$, $c_{14}=c_2$, $c_{34}=c_3$. By Lemma~\ref{electrical}(\ref{el-3}) and a symmetry argument it follows that $\partial C(c_1,c_2,c_3)/\partial c_3=0$, if $c_1=c_2$. So $C(c_1,c_2,c_3)$ cannot be the conductance of a series-parallel network with edge conductances $c_1,c_2,c_3$ because all the edges of such networks are essential.
\end{example}

\section{Proof of main results} \label{proofs}

\subsection{Proof of Theorem \ref{th2}}

Hereafter in an \textit{electrical circuit} or \textit{network} we allow the conductances %$c_{kl}$
to be arbitrary complex numbers with positive real part. This generalization is motivated by \S\ref{alternating} (and describes both direct- and alternating-current circuits).
Theorem~\ref{th2} is a direct consequence of the results of~\S\ref{prelim}:

\begin{proof}[Proof of Theorem \ref{th2}]
Suppose that a rectangle of ratio $c$ can be tiled by rectangles of ratios $c_1$, $\dots$, $c_n$. By Lemma \ref{correspondence-lemma} there is an electrical network of conductance $c$ consisting of edges of conductances $c_1$, $\dots$, $c_n$. For each $k=1,\dots,n$ replace each edge of conductance $c_k$ in the network by an edge of complex conductance $z_k$, $\Real z_k>0$. Let $C(z_1,\dots,z_n)$ be the conductance of the resulting network.
The function $C(z_1,\dots,z_n)$
has properties~(\ref{th2-1})--(\ref{th2-3}) of Theorem~\ref{th2}
by Lemma \ref{electrical}(\ref{el-1}), (\ref{el-2}), and~(\ref{el-5}).
\end{proof}

\subsection{Proof of Theorem~\ref{lfrs}} \label{prooflfrs}

\begin{proof}[Proof of Theorem~\ref{lfrs}]
(\ref{lfrs-3})$\implies$(\ref{lfrs-1}) \cite{FLR} Suppose that condition~(\ref{lfrs-3}) of Theorem~\ref{lfrs} holds and, say, $m$ is odd.
%Let us construct a tiling of a square by rectangles of ratios $c$ and $1/c$.
Take a unit square. Cut off a rectangle of ratio $d_1c$ from the square by a vertical cut. The remaining part is a rectangle of ratio $$
1-d_1c=\cfrac1{d_2c +\dots+\cfrac1{d_mc}}.
$$
Now cut off a rectangle of ratio $1/d_2c$ from the remaining part by a horizontal cut. We get a rectangle of ratio
$$
d_3 c+\cfrac1{d_4 c+\dots+\cfrac1{d_mc}}.
$$
Continue this process alternating vertical and horizontal cuts. Condition~(\ref{lfrs-3}) guarantees that after step $(m-1)$ we get a rectangle of ratio $d_mc$. We obtain a tiling of the square by rectangles of ratios $d_1c$, $1/d_2c$, $d_3c$, $1/d_4c$, $\dots$, $d_mc$.
Since all $d_k\in\mathbb{Q}$, one can chop the tiling into rectangles of ratios $c$ and $1/c$.

(\ref{lfrs-1})$\implies$(\ref{lfrs-2}). Suppose that a square is tiled by rectangles of ratios $c$ and $1/c$.
%Compressing the figure we get a square dissected into squares and rectangles of ratio $R^2$.
By Lemma \ref{correspondence-lemma} there exists an electrical network of conductance $1$ with edge conductances $c$ and $1/c$.
% Consider the corresponding electrical circuit.
Replace each edge of conductance $c$ (respectively, $1/c$) in this network by an edge of conductance $z\in\mathbb{C}$ (respectively, $1/z$).
Let $C(z)$ be the conductance of the obtained network. Then $C(z)\in\mathbb{Q}(z)$ by Lemma \ref{electrical}(\ref{el-1})
and $C(z)$ is odd by Lemma \ref{electrical}(\ref{el-2}).
%Let $q(z)$ be the entire resistance of the obtained network (defined by the Kirchhoff laws).

Since $C(c)=1$ and $C(z)\in\mathbb{Q}(z)$ is nonconstant it follows that $c$ is algebraic.
Let $z$ be an algebraic conjugate of $c$. Then still $C(z)=1$.

Let us prove that $\Real z>0$.
First assume that $\Real z<0$. Then $\Real (-z)>0$ and $\Real (-1/z)>0$.
Thus by Lemma \ref{electrical}(\ref{el-5}) we have $0<\Real C(-z)=-\Real C(z)=-1$, a contradiction. Now assume that $\Real z=0$. Let $z_k\rightarrow z$, where each $\Real z_k<0$. Still $0<\Real C(-z_k)=-\Real C(z_k)\rightarrow-1$, a contradiction. Thus $\Real z>0$.

(\ref{lfrs-2})$\implies$(\ref{lfrs-3})
%\begin{lemma}\label{l1} Consider a rational function $$q(z)=\frac{\prod\limits_{k=1}^n (z^2+b_k^2)}{z\prod\limits_{k=1}^m(z^2+a_k^2)},$$ where $m$ equals $n{-}1$ or $n$, $a_k$ and $b_k$ are real numbers and $0<b_1<a_1<b_2<a_2<\dots$. Then $\Real q(z)>0$ if $\Real z>0$ and $\Real q(z)<0$ if $\Real z<0.$
%\end{lemma}
\cite{FLR}
Let $p(z)$ be a minimal polynomial of $c$. First assume that $\deg p(z)$ is odd. Put $C(z)=\frac{p(-z)-p(z)}{p(-z)+p(z)}$.
Then $C(c)=1$, $C(z)\in \mathbb{Q}(z)$, $\lim_{z\to\infty}C(z)\ne 0$, $C(z)$ is odd, and all the roots of the equation $C(z)=1$ have positive real parts.
By Lemma~\ref{PosReal}(\ref{roots})$\implies$(\ref{continued}) the function $C(z)$ satisfies condition~(\ref{continued})
of Lemma~\ref{PosReal}. Since $C(z)\in \mathbb{Q}(z)$ it follows by the Euclidean algorithm that all $d_k\in \mathbb{Q}$. Substituting $z=c$ we get the required condition. In the event that $\deg p(z)$ is even, replace $C(z)$ by $1/C(z)$ and argue as above.
\end{proof}

\subsection{Proof of Theorem~\ref{th1}}

We use the ideas of \S\ref{prooflfrs} and \S\ref{basicposreal}.

\begin{proof}[Proof of Theorem~\ref{th1}]
(\ref{th1-3})$\implies$(\ref{th1-1}) Analogously to the proof of Theorem~\ref{lfrs}, (\ref{lfrs-3})$\implies$(\ref{lfrs-1}).

(\ref{th1-1})$\implies$(\ref{th1-2}). Suppose that a rectangle of ratio $c$ is tiled by rectangles of ratios $c$ and $1/c$.
Rotating through the angle $\pi/2$ and stretching the figure we get a square
tiled by squares and rectangles of ratio $c^2$. By Lemma~\ref{correspondence-lemma} there exists an electrical network of conductance $1$ with edge conductances $1$ and $c^2$, in which all the edges are essential. Since there is at least one rectangle of ratio $1/c$ in the initial tiling, it follows that the network contains at least one edge of conductance $c^2$.
% Consider the corresponding electrical circuit.
Replace each edge of conductance $c^2$ (respectively, $1$) in the network by an edge of conductance $z\in\mathbb{C}$ (respectively, $w\in\mathbb{C}$).
Let $C(z,w)$ be the conductance of the obtained network. Denote by $C(z)=C(z,1)$.
%Let $q(z)$ be the entire resistance of the obtained network (defined by the Kirchhoff laws).

Let us prove that $c^2$ is algebraic. Indeed, by Lemma~\ref{electrical}(\ref{el-4}) we have $C'(c^2)>0$
because there is at least one essential edge of conductance $c^2$ in the network. Thus $C(z)$ is nonconstant.
By Lemma~\ref{electrical}(\ref{el-1}) it follows that $C(z)\in\mathbb{Q}(z)$. Since $C(c^2)=1$ it follows that $c^2$ is algebraic.

Let $z$ be an algebraic conjugate of $c^2$ distinct from $c^2$ itself. Then $C(z,1)=C(c^2)=1$.

Let us prove that $z$ is a negative real number. First assume that $\Imaginary  {z}<0$. Then $\Real  iz>0$. By Lemma~\ref{electrical}(\ref{el-2}) it follows that $\Real C(iz,i)=\Real  (iC(z,1))=\Real  i=0$. Since $C(iz,i)$ is a rational function it follows that any neighborhood of $iz$ contains a point $z'$ such that $\Real C(z',i)<0$. If the neighborhood is sufficiently small then $\Real  z'>0$ because $\Real  iz>0$. By continuity, a neighborhood of $i$ contains a point $w'$ such that $\Real  w'>0$ and still $\Real C(z',w')<0$. These inequalities contradict Lemma~\ref{electrical}(\ref{el-5}). Case $\Imaginary  {z}>0$ is treated similarly. Assume now that $z> c^2$. Then by Lemma~\ref{electrical}(\ref{el-4}) we have $1=C(z)>C(c^2)=1$, a contradiction. Case $0\le z<c^2$ is treated similarly. Thus $z<0$.

(\ref{th1-2})$\implies$(\ref{th1-3})
%\begin{lemma}\label{l1} Consider a rational function $$q(z)=\frac{\prod\limits_{k=1}^n (z^2+b_k^2)}{z\prod\limits_{k=1}^m(z^2+a_k^2)},$$ where $m$ equals $n{-}1$ or $n$, $a_k$ and $b_k$ are real numbers and $0<b_1<a_1<b_2<a_2<\dots$. Then $\Real q(z)>0$ if $\Real z>0$ and $\Real q(z)<0$ if $\Real z<0.$
%\end{lemma}
Let $p(z)$ be a minimal polynomial of $c^2$. Since the roots of a minimal polynomial are simple it follows that $p(z^2)=(z^2-c^2)\prod_{k=1}^n (z^2+b_k^2)$ for some $b_1>\dots>b_n>0.$  Take a polynomial $q(z)$ with rational coefficients such that $q(z)=z\prod_{k=1}^n(z^2+a_k^2)$, where $a_1>b_1>a_2\dots>b_n.$ Consider the rational function $C(z)=q(z)/(zq(z)-p(z^2))$. We have $C(c)=1/c$.

Let us check that the function $C(z)$ satisfies condition~(\ref{derivatives}) of Lemma~\ref{PosReal}.
Clearly, $C(z)$ is odd and $\lim_{z\to\infty}C(z)\ne0$.
The roots of $C(z)$ are the numbers $0, \pm ia_1,\dots,\pm ia_n$.
A direct evaluation shows that for each $l=1,\dots,n$
$$
C'(\pm ia_l)= -\frac{q'(\pm ia_l)}{p(-a_l^2)} =  \frac{2a_l^2}{(c^2+a_l^2)(a_l^2-b_l^2)}\prod_{k\ne l} \frac{a_k^2-a_l^2}{b_k^2-a_l^2}>0
$$
by the assumption $a_1>b_1>a_2>\dots>b_n>0$. Analogously $C'(0)=-q'(0)/p(0)>0$.
%Thus $C(z)$ satisfies condition~(\ref{derivatives}) of Lemma~\ref{PosReal}.

Then by Lemma~\ref{PosReal}(\ref{derivatives})$\implies$(\ref{continued}) the function $C(z)$ satisfies condition~(\ref{continued}) of Lemma~\ref{PosReal}. Since $C(z)\in \mathbb{Q}(z)$ it follows by the Euclidean algorithm that all $d_k\in \mathbb{Q}$. Substituting $z=c$ we get the required condition.
\end{proof}
%\begin{claim} The function $C(z)$ satisfies condition (3) of Lemma~\ref{PosReal}.
%$\Real z>0$ if and only if $\Real C(z)>0.$
%\end{claim}
%Analyzing the behavior of $C(z)$ near the poles we obtain that for each $y\in\mathbb{R}\setminus \{0\}$ the equation $C(ix)=iy$ has at least $2n{+}1$ real solutions and the equation $C(ix)=0$ has at least $2n$ real solutions. For each $w\in\mathbb{C}$ the equation $C(z)=w$ has at most $2n{+}1$ complex solutions and the equation $C(z)=0$ has at most $2n$ complex solutions. Therefore for each $x\in\mathbb{R}$ all solutions of the equation $C(iz)=ix$ are real. Note that $C(c)>0$ and $C(-c)<0$. Thus by continuity of the function $C(z)$ the claim follows.

\section{Variations}\label{var}

%In this section we study tilings of arbitrary (not necessarily convex) polygons by rectangles. This is equivalent to investigation of planar electrical networks with several terminals.
%This section consists mostly of open problems.

\subsection{Tilings of polygons by rectangles} \label{polygons}

In this subsection we study the following problem.
%generalization of Problem~\ref{main}.
%The results of \S\ref{var} are proved in \S\ref{finalproofs}.

\begin{problem}\label{general}
Which polygons
%When a given (not necessarily convex) polygon $P$
can be tiled by rectangles of given ratios $c_1, \dots, c_n$?
\end{problem}

A related problem of \emph{signed} tilings is solved in~\cite{KK}.

Case $n=1$, $c_1=1$ of the problem is a description of polygons which can be tiled by squares, a problem posed in \cite{FHTW}.
In the case of hexagons the description was obtained by R.~Kenyon; see Figure~\ref{ell}:
%; in the following theorem a misprint in the original statement is corrected). %, cf.~ \cite{FHTW,KK}.

\begin{figure}[htbp]
\definecolor{zzttqq}{rgb}{0.6,0.2,0}
\begin{tikzpicture}[line cap=round,line join=round,>=triangle 45,x=0.5cm,y=0.5cm]
\clip(-4.88,-1.69) rectangle (5.72,6.85);
\fill[color=zzttqq,fill=zzttqq,fill opacity=0.1] (1,6) -- (-4,6) -- (-4,-1) -- (5,-1) -- (5,4) -- (1,4) -- cycle;
\draw [color=zzttqq] (1,6)-- (-4,6);
\draw [color=zzttqq] (-4,6)-- (-4,-1);
\draw [color=zzttqq] (-4,-1)-- (5,-1);
\draw [color=zzttqq] (5,-1)-- (5,4);
\draw [color=zzttqq] (5,4)-- (1,4);
\draw [color=zzttqq] (1,4)-- (1,6);
\draw [color=zzttqq] (0,6)-- (0,3.04);
\draw [color=zzttqq] (0,-1)-- (0,3.48);
\draw [color=zzttqq] (0,4)-- (1,4);
\draw[color=zzttqq] (-2.35,2.94) node[anchor=north west] {$x$};
\draw[color=zzttqq] (0.05,5.4) node[anchor=north west] {$y$};
\draw[color=zzttqq] (2.2,2) node[anchor=north west] {$z$};
\fill [color=zzttqq] (1,6) circle (0.5pt);
\draw[color=zzttqq] (1.28,6.45) node {$A_2$};
\fill [color=zzttqq] (-4,6) circle (0.5pt);
\draw[color=zzttqq] (-3.96,6.45) node {$A_3$};
\fill [color=zzttqq] (-4,-1) circle (0.5pt);
\draw[color=zzttqq] (-4.34,-1.14) node {$A_4$};
\fill [color=zzttqq] (5,-1) circle (0.5pt);
\draw[color=zzttqq] (5.32,-1.06) node {$A_5$};
\fill [color=zzttqq] (5,4) circle (0.5pt);
\draw[color=zzttqq] (5.32,4.4) node {$A_6$};
\fill [color=zzttqq] (1,4) circle (0.5pt);
\draw[color=zzttqq] (1.46,4.4) node {$A_1$};
\end{tikzpicture} 
\caption{L-shaped hexagon.}
\label{ell}
\end{figure}
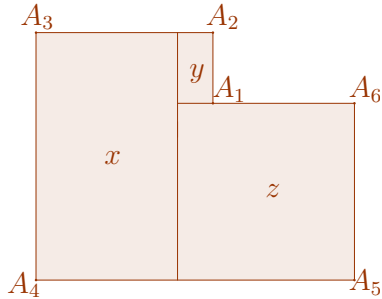

\begin{theorem} \label{th-kenyon}
\textup {(cf.~\cite[Theorem~9]{K})} Let $A_1A_2A_3A_4A_5A_6$ be a hexagon with right angles whose vertices are enumerated counterclockwise starting from the vertex of the nonconvex angle. The hexagon $A_1A_2A_3A_4A_5A_6$ can be tiled by squares if and only if the system
$$
\begin{cases}
A_3A_4\cdot x+ A_1A_2\cdot y&=A_2A_3,\\
A_5A_6\cdot z-A_1A_2\cdot y&=A_6A_1;
\end{cases}
$$
has a nonnegative rational solution $x,y,z$.
\end{theorem}

The proof of sufficiency is easy: given rational $x,y,z\ge 0$ satisfying the system one can dissect the hexagon into $3$ rectangles of ratios $x,y,z$ and then into squares; see Figure~\ref{ell}.

Our aim now is to give a similar criterion for a wide class of polygons.

%Let us give a necessary condition for an arbitrary ``generic'' polygon $P$ to be tilable by squares.

Hereafter $P$ is an \emph{orthogonal} polygon, i.~e., a polygon with the sides parallel to coordinate axes. Hereafter $P$ is \emph{simple}, i.~e., the boundary $\partial P$ is connected. Let $b$ be the number of the sides parallel to the $x$-axis. Enumerate the horizontal sides counterclockwise in $\partial P$.
Let $I_u$ be the \emph{signed length} of the side $u$, where the sign of $I_u$ is ``$+$'' (``$-$'') if the $P$ locally lies below (above) the side $u$. Let $U_u$ be the $y$-coordinate of the side $u$. Assume that $P$ is ``\emph{generic}'' in the following sense: the numbers $U_1,\dots,U_b$ are pairwise distinct.
%(this definition has nothing to do with a measure in the space of all polygons).

We need the following notion \cite{CM}. A sequence  $(p_1,\dots,p_{2k})$ of distinct integers in the interval $[1,b]$ is \emph{circularly ordered}, if a cyclic permutation of
the sequence is an increasing sequence.
Denote by $\Omega_b$ the set of all the real $b\times b$ matrices
$C_{uv}$ having the following properties:
\begin{enumerate}[(1)]
\item \emph{symmetry}: $C_{uv}=C_{vu}$;
    %the matrix $C_{uv}$ is symmetric;
\item \emph{zero sum of each row}: $\sum_{1\le u\le b} C_{uv}=0$;
    %the sum of the entries of the matrix $C_{uv}$ in each row is zero;
\item \emph{total non-negativity of certain ``minors''}: for each circularly ordered sequence $(p_1,\dots,p_{2k})$ we have:
$\det\{ - C_{p_ip_{2k-j+1}} \}_{i,j=1}^k\ge0$.
\end{enumerate}
%The latter property is a restriction on the signs of certain minors of the matrix $C_{uv}$.
E.~g., the set $\Omega_3$ consists of matrices of the form
${\tiny\left(\begin{matrix}
x+y & -x & -y\\
-x  & x+z& -z\\
-y  & -z &y+z
\end{matrix}\right) }$,
where $x,y,z\ge 0$.

%We need the following notation \cite{CM}.
%Let $\Omega_{b}$ be the set of symmetric $b\times b$ matrices $C_{uv}$ with zero sum of each column satisfying the following condition: for each $1\le k\le l\le b$ all the minors of the submatrix $\{-C_{uv}\}_{k\le u\le l, \ v<k \text{ or }v>l}$ are nonnegative. %Denote $-\bar\Omega_{b}=\{-C:C\in\bar\Omega_{b}\}$.

%Here $2b$ denotes the number of vertices in $P$.

\begin{theorem} \label{squaring-polygons}
Let $P$ be a ``generic'' simple orthogonal polygon with $b$ horizontal sides having signed lengths $I_1,\dots,I_b$ and $y$-coordinates $U_1,\dots,U_b$.
%Suppose that the polygon $P$ can be tiled by squares; then there is a matrix $C_{uv}\in \Omega_{b}$ with rational entries such that $I_v=\sum_{u=1}^{b} C_{uv}U_u$ for each $v=1,\dots, b$.
Then the following $2$ conditions are equivalent:
\begin{enumerate}[(1)]
\item \label{s-p-1} the polygon $P$ can be tiled by squares;
\item \label{s-p-2} there is a matrix $C_{uv}\in \Omega_{b}$ with rational entries such that $I_v=\sum_{u=1}^{b} C_{uv}U_u$ for each $v=1,\dots, b$.
\end{enumerate}
\end{theorem}

%\begin{problem} Is the converse theorem true for $b\ge 4$?
%\end{problem}

Cases $b=2$ and $b=3$ of this theorem
%and its converse is
are equivalent to Theorem~\ref{dehn}  and \cite[Theorem~9]{K}, respectively.
%\begin{example}\label{nongeneric}
Theorem~\ref{squaring-polygons} is algorithmic for the particular class of polygons $P$ such that $U_1,\dots,U_b$ are linearly independent over $\mathbb{Q}$.
Proof of the theorem is constructive, i.~e., gives an algorithm to construct the required tiling, if the latter exists. However, when $b\ge 4$ the algorithm becomes complicated compared to the cases of $b=2$ and $b=3$.
%and hard to implement because the problem is ill-posed.
Theorem~\ref{squaring-polygons} %(1)$\implies$(2)
does not necessarily hold for ``\emph{nongeneric}'' polygons, e.~g., for an orthogonal polygon with parameters
$$U_1=U_3=0, \quad U_2=2, \quad U_4=-4, \quad I_1=\sqrt{2},\quad I_2=2, \quad I_3=2-\sqrt{2},\quad I_4=-4.$$
%\end{example}
%For a generalization to tiling by rectangles see Remark~\ref{mostgeneral} below.

%\begin{example}\label{nongeneric} A \emph{nongeneric} polygon with $U_1=U_3=0$, $U_2=2$, $U_4=-4$, $I_1=\sqrt{2}$, $I_2=2$, $I_3=2-\sqrt{2}$, $I_4=-4$ satisfies condition~(1) of Theorem~\ref{squaring-polygons} but not~(2).
%\end{example}

%\begin{theorem} \label{rectangling-polygons} Suppose that a simple polygon $P$ can be tiled by rectangles of ratios $c_1, \dots, c_n$. Then $I_v=\sum_{u=1}^{b} C_{uv}(c_1, \dots, c_n)U_u$, $v=1,\dots,b$, for some rational functions $C_{uv}(z_1,\dots,z_n)$ such that
%\begin{enumerate}[(1)]
%\item $C_{uv}(z_1,\dots,z_n)\in\mathbb{Q}(z_1,\dots,z_n)^{b \times b}$;
%\item $C_{uv}(z_1,\dots,z_n)$ is degree $1$ homogeneous;
%\item if $\Real  z_k>0$ for all $k=1,\dots,n$ then $\Real C_{uv}(z_1,\dots,z_n)$ is almost positively definite;
%\item if $z_k >0$ for all $k=1,\dots,n$ then $C_{uv}(z_1,\dots,z_n)\in -\bar\Omega_{b}$.
%\end{enumerate}
%\end{theorem}

%\begin{problem} Is the converse theorem true?
%\end{problem}

%Neither of the above results is ''algorithmic''.
%However, it is possible to get ''algorithmic'' results in some particular cases of Problem~\ref{general}, cf.~ Theorem~\ref{lfrs} above:

We prove Theorem~\ref{squaring-polygons} in \S\ref{finalproofs}. We also give a short proof of the following result:

\begin{theorem} \textup{\cite{SuDi}} \label{polygon-square} A ``generic'' simple orthogonal polygon with rational vertices can be tiled by rectangles of ratios $c$ and $1/c$ if and only if a square can be tiled by rectangles of ratios $c$ and~$1/c$.
\end{theorem}

%We prove Theorems~\ref{squaring-polygons}  %\ref{rectangling-polygons}
%and~\ref{polygon-square} in the end of this section.

%\begin{theorem} \label{th3}
%Suppose that a polygon of area $S$ can be dissected into rectangles of ratios $R_1,\dots,R_n$. Suppose that the $x$-coordinates of all vertices of the polygon are rational. Then there exists a function $E(z_1,\dots,z_n)\in\mathbb{Q}(z_1,\dots,z_n)$ such that the following conditions hold:

%(1) $E(R_1, R_2, \dots, R_n)=S$;

%(2) $E(cz_1,cz_2, \dots, cz_n)\equiv cE(z_1,z_2, \dots, z_n)$;

%(3) if $\Real z_k>0$ for all $k=1,\dots,n$ then $\Real E(z_1,z_2, \dots, z_n)>0$.
%\end{theorem}

%\begin{conjecture} The converse theorem is also true.
%\end{conjecture}

%\subsection{Electrical circuits with several terminals} \label{responsesection}

\subsection{Electrical impedance tomography}\label{eit}

Our approach to Problem \ref{general} follows the idea of \cite{K, CGV} and uses electrical networks with several terminals.

Hereafter we allow electrical circuits to have several boundary vertices $1$, \dots, $b$ with prescribed voltages $U_1$, \dots, $U_b$.
If an electrical circuit is planar then we assume that the boundary vertices are enumerated counterclockwise along the boundary of the unit disc. We do \emph{not} assume that an electrical circuit is connected but require that each connected component contains a boundary vertex. The voltages and currents in such circuits are defined by the Ohm law~(C) and the Kirchhoff current law~(I) from~\S\ref{prelim}.
%Let $1,\dots,b$ be the boundary vertices of the electrical network., and $b+1,\dots,n$ are the remaining ({\it interior}) ones. Let $U_1,U_2,\dots,U_b$ be the given voltages at the boundary vertices, and $U_{b+1},\dots,U_n$ be the voltages at the interior vertices determined by the Ohm and Kirchhoff laws.
%Let $r_{ij}$ be the conductance of the edge between nodes $i$ and $j$.

Consider the linear map $\mathbb{C}^b\to \mathbb{C}^b$ which takes the vector of voltages $(U_1,\dots,U_b)$ to the vector of \emph{incoming currents} $(I_1,\dots, I_b)=(\sum_{k=1}^{n} I_{1k},\dots,\sum_{k=1}^{n} I_{bk})$ flowing inside the network through the vertices $1,\dots,b$, respectively. The matrix $C_{uv}$ of this linear map is called {\it the response} of the network.
This matrix is symmetric~\cite{CM}. For $b=2$ the response of a network is ${\tiny\left(\begin{matrix} C & -C\\ -C & C\end{matrix}\right)}$,
where $C$ is the conductance of the network.
%(This definition differs by a sign from the one in \cite{C}.)
%We have the following generalization of the results of \S\ref{prelim}, cf.~\cite{Duffin}.

%The value $E=\sum_{1\le u,v\le b} C_{uv} U_u\bar U_v$ is called {\it energy dissipation}.
%We consider $E(r_{ij})$ as the function of the conductances $r_{ij}$ for fixed voltages $v_i$.
%A symmetric $b\times b$ matrix $C$ is \textit{almost positively definite}, if it is non-negatively definite and $x^{T}Cx=0$ only if $Cx=0$.
%The following (probably new) results are proved completely analogously to the proof of the results of \S\ref{prelim} above.

%\subsection{Relation with tilings}

%\subsection{Electrical impedance tomography}

%Analogously to the continuous case the following problems for discrete harmonic functions are natural to study.
%The rest of \S\ref{var} consists of open problems and a few examples.

We reduce the results of \S\ref{polygons} to
\emph{the inverse problem} for electrical networks, which is
%following problems even more interesting in themselves:
%\begin{itemize}
%\item \emph{Direct problem.} Describe possible responses of electrical networks.
%\item \emph{Inverse problem.}
to synthesize a network having given response.
%to describe all the networks having a given response.
%\end{itemize}
This is a discrete analogue of \emph{electrical impedance tomography} \cite{Cal,SU}. This problem was posed in \cite{LSy} and solved for planar direct-current networks in \cite{CIM, CM, C, CGV}. Let us state certain deep results of Y.~Colin de Verdi\`ere, E.~B.~Curtis, and J.~A.~Morrow.
%For instance, the following deep results are known.
%Let us state the following deep results.

%We need the following definitions \cite{C,CM}. Let $G$ be a network and $u_1,\dots,u_k,v_1,\dots,v_k$ be a sequence of boundary vertices. The sequence is \emph{connected through $G$} if $G$ contains $k$ disjoint paths $p_1,\dots,p_k$ such that $p_l$ starts at $u_l$, ends at $v_l$ and passes through no other boundary vertices. Denote by $\pi(G)$ be the set of all the sequences connected through $G$. The network $G$ is \emph{minimal} if for each network $G^-$ obtained from $G$ by either deleting or contracting an edge we have $\pi(G^-)\ne \pi(G)$. %Denote by $\mathbb{R}_+$ the set of positive real numbers.

\begin{theorem} \label{CGV5} %\textup{\cite[Proposition 10]{C}}
\textup{\cite[Theorem~5]{CM, CIM, CGV}}
%The response of a planar electrical network with $b$ boundary vertices and positive edge conductances belongs to the set $\Omega_b$. Each element of the set $\Omega_b$ is the response of a minimal planar electrical network with $b$ boundary vertices and positive edge conductances.
The set of all possible responses of planar electrical networks with $b$ boundary vertices and positive edge conductances is the set $\Omega_b$ from~\S\ref{polygons}.
%The following $3$ conditions are equivalent:
%\begin{enumerate}[(1)]
%\item $C_{uv}$ is the response of a planar electrical network with $b$ boundary vertices and positive edge conductances;
%\item $C_{uv}$ is the response of a \textup{critical} planar electrical network with $b$ boundary vertices and positive edge conductances;
%\item $C_{uv}\in\Omega_b$.
%\end{enumerate}
\end{theorem}

An electrical network is \emph{minimal} (or \emph{critical}) if it has minimal number of edges among all planar electrical networks with positive edge conductances and with the same response. The minimality of a network depends only on its graph \cite{CGV}. In \cite[\S9]{CM, CIM} an algorithm for finding edge conductances in a minimal network with given response is presented. This algorithm implies easily the following result (D.~Ingerman, J.~A.~Morrow, private communication).

\begin{theorem} \label{inverse}
%\textup{\cite[\S6.4]{CM}}
Conductances of the edges in a minimal electrical network are uniquely determined by the response of the network. Each edge conductance is a rational function with rational coefficients in the entries of the response.
\end{theorem}

For \emph{alternating}-current networks the inverse problem is probably open \cite{Foster}.
Let us state a basic folklore result. The rest of~\S\ref{var} is not used in the proof of the above theorems. However, the authors believe that any
further progress in Problems~\ref{mainproblem} and~\ref{general} requires
%the solution of the inverse problem.
a generalization of the results stated below.

%Let us state some results for alternating-current circuits (the rest of \S\ref{var} is not used in the proofs of the above results). %of \S\ref{polygons}).

% The following notion is inappropriate because
% we need $C U=0$ instead of $Re C=0$ !!!
%A real $b\times b$ matrix $R$ is \emph{almost positively definite}, if for any $U\in\mathbb{R}^b$ we have either $U^T R U>0$ or else $R U=0$.

\begin{theorem}\label{detection} For $b=2$ or $b=3$ the following $2$ conditions are equivalent\textup{:}
\begin{enumerate}[(1)]
\item \label{det-1} $C_{uv}$ is the response of a connected electrical network with $b$ boundary vertices and with edge conductances having positive real parts\textup{;}
\item \label{det-2} $C_{uv}$ is a complex $b\times b$ matrix having the following properties\textup{:}
\begin{itemize}
\item \emph{symmetry}\textup{:} $C_{uv}=C_{vu}$\textup{;}
\item \emph{zero sum of each row}\textup{:} $\sum_{1\le u\le b}C_{uv}=0$\textup{;}
\item \emph{positive reality}\textup{:} $\sum_{1\le u,v\le b}(\Real C_{uv})U_u U_v\ge 0$ for any $U_1,\dots,U_b\in\mathbb{R}$\textup{;}
\item \emph{almost positively definiteness}\textup{:}    the latter inequality is strict unless $U_1=\dots=U_b$.
\end{itemize}
\end{enumerate}
\end{theorem}

\begin{problem} Does this result remain true for arbitrary $b\ge 4$?
\end{problem}

%Theorem~\ref{detection} is proved at the end of the section.

%The problem of \textit{determining the structure of a network by its response} might be of practical interest.
%One can see that for each direct-current network there is a (not necessarily planar) direct-current network without nonboundary vertices and with the same response. Whereas

%Unlike nonplanar direct-current networks  \emph{nonboundary vertices in alternating-current networks can be detected by the response}.
%For instance, by Theorem~\ref{detection} there are electrical networks with response
%{\tiny $\left(\begin{matrix}
%2 & 1 & -3\\
%1 & 2 & -3\\
%-3& -3& 6 \\
%\end{matrix}\right)$}; any such network
%necessarily has nonboundary vertices.

%(Indeed, the response of a network without nonboundary vertices cannot have positive entries outside the diagonal.)
%\end{example}

%A result concerning the inverse problem is presented in \S\ref{inverse}.

%\begin{problem}
%Which is the least possible number of vertices in an electrical network with a given response?
%\end{problem}

%\begin{problem} Suppose that $C_{uv}(c_{12}, c_{13}, \dots, c_{n-1,n} )$ is the response of an electrical network with $b$ boundary vertices and with edge conductances $c_{12}, c_{13}, \dots, c_{n-1,n}$. Let $C_{uv}$ be a complex $b\times b$ matrix. Is it true that the equation $C_{uv}(c_{12}, c_{13}, \dots, c_{n-1,n} )=C_{uv}$ has either no or unique or infinitely many (continuum of) solutions $(c_{12}, c_{13}, \dots, c_{n-1,n} )\in\mathbb{C}^{b(b-1)/2}$?
%\end{problem}

\subsection{Random walks}

A {\it random walk} on an electrical network (or on a \emph{weighted graph}) is the Markov chain with the transition matrix $P_{kl}=c_{kl}/\sum_{j=1}^{n}c_{jk}$. %Here we require the same assumptions on the graph $G$ and the numbers $c_{kl}$ as in \S\ref{prelim}. Under these assumptions the
Such Markov chain is ergodic and reversible.
Denote by $k_1l_1,\dots,k_ml_m$ all the edges of the Markov chain.
The following theorem allows to translate the results of \S1--\S2 to the language of random walks.

\begin{theorem} \textup{\cite[page 42]{DS}}
Let $P(c_{k_1l_1}, \dots,c_{k_ml_m})$ be the probability that a random walk starting at vertex $1$ reaches vertex $2$ before returning to $1$. Let $C(c_{k_1l_1}, \dots,c_{k_ml_m})$ be the conductance of the network \textup{(}with boundary vertices $1$ and $2$\textup{)}.
Then $P(c_{k_1l_1}, \dots,c_{k_ml_m})={C(c_{k_1l_1}, \dots,c_{k_ml_m})}/{(c_{12}+\dots+c_{1n})}$.
\end{theorem}

%\begin{figure}[htbp]
%\includegraphics{walksnetworks.pdf}
%\caption{Walks and networks}
%\label{walksnetworks}
%\end{figure}

%By Lemma~\ref{electrical}(1) $P(c_{kl})$ is a rational function in $c_{kl}$. A ``translation'' of Lemma~\ref{electrical}(5) is:

%\begin{theorem}
%If $\Real c_{kl}>0$ for all $1\le k<l\le n$ then $\Real \left((c_{12}+\dots+c_{1n})P(c_{kl})\right)>0$.
%\end{theorem}

%\begin{problem} Does this result remain true for nonreversible Markov chains?
%\end{problem}

For instance, a translation of Lemmas~\ref{electrical}(\ref{el-1}) and~(\ref{el-5}) is
the following result:

\begin{corollary} \label{random}
The probability $P(c_{k_1l_1}, \dots,c_{k_ml_m})$ is a rational function in %edge weights
$c_{k_1l_1}, \dots,c_{k_ml_m}$ such that\textup{:} if $\Real c_{k_1l_1}, \dots,\Real c_{k_ml_m}>0$ then $\Real \left((c_{12}+\dots+c_{1n})P(c_{k_1l_1}, \dots,c_{k_ml_m})\right)>0$.
\end{corollary}

The latter result does not necessarily hold for \emph{nonreversible} Markov chains; e.~g., for a Markov chain with vertices $1$, $2$, $3$, $4$ and oriented edges $14$, $42$, $43$.
%we have $P(c_{14},c_{42},c_{43})=c_{42}/(c_{42}+c_{43})$.
%A natural question is: \emph{does this result remain true for nonreversible Markov chains}?

In \cite{Fomin} a related result for non-reversible Markov chains was obtained by means of the results of electrical impedance tomography stated in~\S\ref{eit}.

Nonreversible planar Markov chains have a geometric interpretation as tilings of trapezoids by trapezoids \cite{K}.
Here a \emph{trapezoid} is a quadrilateral with two sides parallel to the $x$-axis. %The \emph{ratio} of the trapezoid is the length of the horizontal middle edge divided by the hight. 
Natural problems are: \emph{generalize the results of the paper to
%\begin{itemize}
tilings by trapezoids; infinite tilings; signed tilings}.
%\end{itemize}
%Probably, electrical networks with edges having different conductances in different directions (or conductances of arbitrary sign --- in case of signed tilings) might be useful for answering the question.

\section{Generalization of main ideas}\label{final}

\subsection{Electrical circuits} \label{elcirc}

%\subsection{Proof of the results of \S\ref{prelim}}

Our approach is based on a generalization of the results of \S\ref{prelim} to electrical circuits with $b$ terminals.
Short proofs of the results of \S\ref{prelim} are obtained in this section for the particular case of $b=2$. Our proof of Lemma~\ref{responseclaim}(\ref{rc-3}) generalizing Lemma~\ref{electrical}(\ref{el-3}) is probably new.
For simplicity in this subsection assume that the circuits do not have multiple edges.
%We prove the results of \S\ref{prelim} in this general situation.
%Our proofs of Lemma~\ref{correspondence-general} and Lemma~\ref{responseclaim}(3), generalizing
%Lemma~\ref{correspondence-lemma} and Lemma~\ref{electrical}(3), respectively, are probably new.
All the proofs are based on the following fundamental \textit{energy conservation law}.

\begin{lemma} \label{energy} \label{fullpower} Let $E(U,I)$ be a bilinear or sesquilinear function.
Consider an electrical network with the vertices $1, \dots, n$ such that $1, \dots, b$ are the boundary ones.
Suppose that the numbers $U_k$, where $1\le k\le n$, and $I_{kl}$, where $1\le k,l\le n$, obey the Ohm law \textup{(C)} and the Kirchhoff current law \textup{(I)} from \S\ref{prelim}. Set $I_u=\sum_{k=1}^n I_{uk}$. Then
$$\sum\limits_{1\le k<l\le n}E(U_k-U_l,I_{kl})=\sum\limits_{1\le u\le b}E(U_u,I_u).$$
\end{lemma}

%\begin{claim} \label{fullpower} Let $E(U,I)$ be a bilinear function.
%Suppose that the numbers $U_k$, $2< k\le n$, and $I_{kl}$, $1\le k,l\le n$ satisfy the laws \textup{(C),(I)}. Let $I_1=\sum_{k=2}^n I_{1k}$. Then
%$$\sum\limits_{1\le k<l\le n}E(U_k-U_l,I_{kl})=E(U_1-U_2,I_1).$$
%\end{claim}

We usually apply this lemma to the \emph{energy dissipation} function $E(U,I)=\Real (U\bar I)$.

\begin{proof}[Proof of Lemma \ref{energy}] By law (C) we have $I_{lk}=-I_{kl}$. Hence by law (I) we have
$$
\sum\limits_{1\le k<l\le n}E(U_k-U_l,I_{kl})=\sum\limits_{k=1}^n E(U_k,\sum\limits_{l=1}^n I_{kl})=
\sum\limits_{1\le u\le b}E(U_u,I_u).
$$
\end{proof}

%\begin{proof}[Proof of Claim \ref{fullpower}] Clearly, by law (I) we have $\sum_{k=1}^n I_{2k}=-I_1$. Thus
%$$\sum\limits_{1\le k<l\le n}E(U_k-U_l,I_{kl})=\sum\limits_{k=1}^n E(U_k,\sum\limits_{l=1}^n I_{kl})=E(U_1-U_2,I_1).$$
%\end{proof}

%\begin{claim} (Energy preserving law) \label{fullpower} Let the numbers $U_k$, $3\le k\le n$, and $I_{kl}$, $1\le k,l\le n$ satisfy the Ohm and the Kirchhoff laws (1)--(2). Let $I=-\sum_{k=2}^n I_{1k}$. Then
%$\sum\limits_{1\le k<m\le n}(U_m-U_k)\bar{I}_{km}=(U_2-U_1)\bar{I}$.
%\end{claim}

%\begin{proof} Clearly, by Kirchhoff law (2) we have $\sum_{k=1}^n I_{2k}=I$. Thus
%\sum\limits_{1\le k<m\le n}(U_m-U_k)\bar{I}_{km}=\sum\limits_{k=1}^n U_k\sum\limits_{m=1}^n (-\bar{I}_{km})=(U_2-U_1)\bar{I}.
%\end{proof}

Let us prove Theorem~\ref{weyl} for electrical circuits with $b$ boundary vertices and with complex edge conductances having positive real part.

\begin{proof}[Proof of Theorem \ref{weyl}] \textit{Uniqueness}. %Suppose that for each $k,l$ either $\Real c_{kl}>0$ or $c_{kl}=0$.
Suppose there are two collections of currents $I^{I,II}_{kl}$ and voltages $U^{I,II}_k$ obeying laws (C), (I). Then their differences $I_{kl}=I^{I}_{kl}-I^{II}_{kl}$ and $U_{k}=U^{I}_{k}-U^{II}_{k}$ obey laws (C), (I) for zero incoming voltages $U_1=\dots=U_b=0$. Then by %Claim~\ref{fullpower}
Lemma~\ref{energy} we have
$$
\sum\limits_{1\le k<l\le n}\Real \bar c_{kl}|U_k-U_l|^2=\sum\limits_{1\le k<l\le n}\Real ((U_k-U_l)\bar{I}_{kl})=
\sum\limits_{1\le u\le b}\Real (U_u\bar I_u)=0.
$$
For each $k,l$ we have either $\Real c_{kl}>0$ or $c_{kl}=0$. Thus each summand $\Real \bar c_{kl}|U_k-U_l|^2=0$. Since all the connected components of the circuit contain boundary vertices it follows that all the voltages $U_{k}$ are equal to each other. Hence $U_k=0$, $I_{kl}=0$, and thus $I^{I}_{kl}=I^{II}_{kl}$, $U^{I}_{k}=U^{II}_{k}$ for each $k,l$.

\textit{Existence}. The number of equations in the system (C), (I) equals the number of variables.
We have proved that
%By the above previous paragraph
the system has a unique solution for $U_1=\dots=U_b=0$. By the finite-dimensional Fredholm alternative it has a solution for any $U_1$, \dots, $U_b$.
\end{proof}

%\begin{remark} The assumption $\Real c_{kl}>0$ in Theorem~\ref{weyl} is essential by example in figure~\ref{examples} for $a=i,b=-i$.
%\end{remark}

%Let us summarize the properties of the formulas for response.
The following result generalizes Lemma~\ref{electrical}.

\begin{lemma} \label{responseclaim} Suppose that an electrical network has $b$ boundary vertices and $m$ edges of conductances $c_{1},\dots, c_{m}$. Then the response of the network $C_{uv}(c_{1},\dots, c_{m})$ has the following properties\textup{:}
\begin{enumerate}[(1)]
\item \label{rc-1}
    %\emph{rationality}\textup{:}
    $C_{uv}(c_{1},\dots, c_{m})\in \mathbb{Q}(c_{1},\dots, c_{m})^{b \times b}$\textup{;}
    %\item $C_{uv}(c_{1},\dots, c_{m})$ is linear-fractional with respect to each variable $c_{j}$;
\item \label{rc-2}
    %\emph{homogeneity}\textup{:}
    $C_{uv}(tc_{1},\dots, tc_{m})=tC_{uv}(c_{1},\dots, c_{m})$\textup{;}
\item \label{rc-3}
    $\frac{\partial}{\partial c_{j}}C_{uv}(c_{1},\dots, c_{m}) =(V_{ku}-V_{lu})(V_{kv}-V_{lv})$, where $k$ and $l$ are the endpoints of the edge $j$ and $V_{pq}$ is the matrix of the linear map $(U_1,\dots,U_b)\mapsto (U_1,\dots,U_n)$\textup{;}
\item \label{rc-4}
    %\emph{monotonicity}\textup{:}
    if $c_{1},\dots, c_{m}>0$ then $\frac{\partial}{\partial c_{j}}C_{uv}(c_{1},\dots, c_{m})$ is non-negatively definite\textup{;}
%\item $C_{uv}(c_{12}, c_{13}, \dots, c_{n-1,n} )$ is symmetric;
\item \label{rc-5}
    %\emph{}\textup{:}
    if $\Real c_{1},\dots,\Real c_{m}>0$ then $\Real C_{uv}(c_{1},\dots, c_{m})$ is non-negatively definite.
%\item \label{rc-6} $C_{uv}(c_{1},\dots, c_{m})=C_{vu}(c_{1},\dots, c_{m})$.
\end{enumerate}
\end{lemma}

\begin{proof}[Proof of Lemma \ref{responseclaim}] (\ref{rc-1})
By Theorem~\ref{weyl} and the Crammer rule the solution $\{I_{kl}(U_1,\dots,U_b)\}$ of the system of linear equations (C), (I)
consists of linear functions in $U_1, \dots, U_b$ with coefficients being rational functions in $c_{1}, \dots, c_m$.
So the entries of the matrix of the linear map $(U_1,\dots,U_b)\mapsto{\sum_{k=1}^n I_{uk}(U_1,\dots,U_b)}$ are rational functions in $c_{1}, \dots, c_m$.

%(1) Fix an edge $j$. Consider the matrix $M(c_j)$ of the system of linear equations (C), (I) as a function in the parameter $c_{j}$. The matrix $M(c_{j})$ contains a $2\times 2$ minor
%{\tiny $\left(
%\begin{matrix}
%-c_{j} & c_{j} \\
%c_j    & -c_{j}
%\end{matrix}
%\right)$}.
%The entries of $M(c_{j})$ outside this $2\times 2$ minor do not depend on $c_{j}$. Hence $\det M(c_{j})$ being the denominator of $C_{uv}(c_{1},\dots, c_{m})$ is linear in $c_{j}$. Analogously, the nominator is linear in $c_{j}$.

(\ref{rc-2}) Consider the system of linear equations obtained from laws (C), (I) by substituting $t c_{1},\dots, t c_{m}$ for $c_{1},\dots, c_{m}$. It defines the same voltages as the initial one and the currents are scaled by $t$. So $C(tc_{1},\dots, tc_{m})=tC(c_{1},\dots, c_{m})$. %\sum\limits_{k=2}^n tc_{1k}\frac {U_1-U_k}=tC(c_{kl})$.

(\ref{rc-3})
Fix the voltages $U_1,\dots,U_b$ and all the conductances $c_{pq}$ except $c_{kl}$.
Consider the voltages and the currents in the circuit as functions in $c_{kl}$. Set $E(U,I)=U\frac{\partial I}{\partial c_{kl}}-\frac{\partial{U}}{\partial c_{kl}}I$. %(\textit{Lagrange function}).
Then $E(U_k-U_l,I_{kl})=(U_k-U_l)^2$ and $E(U_p-U_q,I_{pq})=0$ for $pq\ne kl$.
Thus by Lemma~\ref{fullpower} we have
\begin{multline*}\sum_{1\le u,v\le b}\frac{\partial C_{uv}}{\partial c_{kl}}U_uU_v=\sum_{1\le u\le b}E(U_u,I_u)=\sum\limits_{1\le p<q\le n}E(U_p-U_q,I_{pq})=\\
=(U_k-U_l)^2=\sum_{1\le u,v\le b}(V_{ku}-V_{lu})(V_{kv}-V_{lv})U_uU_v.
\end{multline*}

%Let $E(U,I)=\frac{\partial{U}}{\partial c_{uv}}I-U\frac{\partial I}{\partial c_{uv}}$. %(\textit{Lagrange function}).
%Then $E(U_u-U_v,I_{uv})=(U_u-U_v)^2$ and $E(U_k-U_l,I_{kl})=0$ for $kl\ne uv$.
%Thus by Claim~\ref{fullpower} we have
%$$\frac{\partial C}{\partial c_{uv}}(U_1-U_2)^2=E(U_1-U_2,I)=\sum\limits_{1\le k<l\le n}E(U_k-U_l,I_{kl})=(U_u-U_v)^2.$$

(\ref{rc-4}) This follows directly from the latter formula.

(\ref{rc-5}) For each $k,l$ we have either $\Real c_{kl}>0$ or $c_{kl}=0$. Take $U_1,\dots,U_b\in\mathbb{R}$ (although the argument below works for complex voltages as well).
By Lemma~\ref{energy} we have
\begin{multline*}
\sum_{1\le u,v\le b} (\Real C_{uv})U_u \bar U_v=
\sum_{1\le u,v\le b} \Real(U_u \bar C_{uv}\bar U_v)=
\sum_{1\le u\le b} \Real (U_u \bar I_u)=\\
=\sum\limits_{1\le k<l\le n}\Real ((U_k-U_l)\bar I_{kl})=
\sum\limits_{1\le k<l\le n}\Real c_{kl}|U_k-U_l|^2\ge 0.
\end{multline*}
\end{proof}

\begin{remark} \label{strict} For a connected network the latter inequality is strict unless $U_1=\dots=U_b$. We have $C_{uv}=C_{vu}$ (which is also proved by an appropriate energy conservation argument).
\end{remark}

%\begin{remark}
%The classical proof of Lemma \ref{main}(4) is based on the following observation: \emph{if all $c_{kl}>0$ then the solution  $(U_{b+1},\dots,U_n)$ of the system of linear equations (C), (I) minimizes the energy dissipation function $\sum_{1\le k<l\le n}c_{kl}(U_k-U_l)^2$ for fixed $U_1,\dots,U_b$}. \cite{Rayleigh}
%\end{remark}

%\subsection{Inverse problem} \label{inverse}

%The following result is probably new.

%\begin{lemma} \label{inverse-lemma} NOT CORRECT!!! Suppose that $C_{uv}(c_1, \dots, c_{m})$ is the response of an electrical network with $b$ boundary vertices and with $m$ edges of conductances $c_{1}, \dots, c_{m}$. Then there is an integer number $0\le r\le m$ and $m-r$ functions
%$$
%c_{k}(c_1,\dots,c_r,C_{11},C_{12},\dots,C_{b-1,b})\in\mathbb{Q}(c_1,\dots,c_r,C_{11},C_{12},\dots,C_{b-1,b}),
%\qquad r+1\le k\le m
%$$
%such that the system of equations
%$$
%C_{uv}(c_1, \dots, c_{m})=C_{uv}, \qquad 1\le u\le v\le b,
%$$
%in variables $c_1, \dots, c_{m}$ is equivalent to
%$$
%c_k=c_{k}(c_1,\dots,c_r,C_{11},C_{12},\dots,C_{b-1,b}), \qquad r+1\le k\le m.
%$$
%\end{lemma}

\subsection{Circuits and tilings} \label{correspond}

We extend the ideas of \cite{BeSch,K} to get the following generalization of Lemma~\ref{correspondence-lemma}. The main difficulty is to guarantee the absence of ``saddle'' vertices in the network. 
%Part (\ref{c-g-2})$\implies$(\ref{c-g-1}) of the following result is probably new, cf.~\cite{BeSch,K}.

\begin{lemma} \label{correspondence-general} Let $P$ be a ``generic'' simple orthogonal polygon with horizontal sides of signed lengths $I_1,\dots,I_b$ and $y$-coordinates $U_1,\dots,U_b$. Then the following $2$ conditions are equivalent\textup{:}
\begin{enumerate}[(1)]
\item \label{c-g-1} the polygon $P$ can be tiled by $m$ rectangles of ratios $c_1,\dots,c_m$\textup{;}
%a polygon with sides parallel to the coordinate axes, vertical sides having lengths $I_1,\dots, I_b$ and $x$-coordinates $U_1, \dots, U_b$, can be tiled by rectangles of ratios $c_1,\dots,c_n$;
\item \label{c-g-2} there is a planar electrical circuit with $b$ boundary vertices, $m$ essential edges of conductances $c_1,\dots,c_m>0$, incoming voltages $U_1, \dots, U_b$ and incoming currents $I_1, \dots, I_b$.
\end{enumerate}
\end{lemma}

\begin{remark}
Condition~(\ref{c-g-2}) itself does not guarantee the existence of a orthogonal polygon with horizontal sides of signed lengths $I_1,\dots,I_b$ and $y$-coordinates $U_1,\dots,U_b$. %Example~\ref{nongeneric}
%\end{remark}
%\begin{remark}
%Example after Theorem~\ref{squaring-polygons} shows that
Lemma \ref{correspondence-general}(\ref{c-g-1})$\implies$(\ref{c-g-2}) is not necessarily true for ``\emph{nongeneric}'' polygons.
\end{remark}

\begin{figure}[htbp]
\definecolor{qqwwff}{rgb}{0,0.4,1}
\definecolor{ffqqqq}{rgb}{1,0,0}
\definecolor{zzttqq}{rgb}{0.6,0.2,0}
\begin{tikzpicture}[line cap=round,line join=round,>=triangle 45,x=0.7cm,y=0.7cm]
\clip(-1.66,0.8) rectangle (11.1,5.19);
\fill[color=zzttqq,fill=zzttqq,fill opacity=0.1] (-1.5,3) -- (1.5,3) -- (1.5,5) -- (4.5,5) -- (4.5,1) -- (-1.5,1) -- cycle;
\draw [color=zzttqq] (-1.5,3)-- (1.5,3);
\draw [color=zzttqq] (1.5,3)-- (1.5,5);
\draw [color=zzttqq] (1.5,5)-- (4.5,5);
\draw [color=zzttqq] (4.5,5)-- (4.5,1);
\draw [color=zzttqq] (4.5,1)-- (-1.5,1);
\draw [color=zzttqq] (-1.5,1)-- (-1.5,3);
\draw [color=zzttqq] (1.5,3)-- (3.5,3);
\draw [color=zzttqq] (3.5,5)-- (3.5,1);
\draw [line width=1pt,color=ffqqqq] (2,5)-- (4,5);
\draw [line width=1pt,color=ffqqqq] (-1,3)-- (3,3);
\draw [line width=1pt,color=ffqqqq] (-1,1)-- (4,1);
\draw [line width=1pt,dash pattern=on 4pt off 4pt,color=qqwwff] (4,5)-- (4,1);
\draw [line width=1pt,dash pattern=on 4pt off 4pt,color=qqwwff] (2.5,5)-- (2.5,3);
\draw [line width=1pt,dash pattern=on 4pt off 4pt,color=qqwwff] (1,3)-- (1,1);
\draw (5.48,3.7) node[anchor=north west] {$\mathbf{ \rightarrow}$};
\draw [color=zzttqq,fill=zzttqq,fill opacity=0.1] (9,3) circle (1.98);
\draw [line width=1pt,dash pattern=on 4pt off 4pt,color=qqwwff] (9,4.98)-- (9,1.02);
\draw [line width=1pt,dash pattern=on 4pt off 4pt,color=qqwwff] (9,4.98)-- (7.02,3);
\draw [line width=1pt,dash pattern=on 4pt off 4pt,color=qqwwff] (7.02,3)-- (9,1.02);
\fill [color=ffqqqq] (9,4.98) circle (2.0pt);
\fill [color=ffqqqq] (7.02,3) circle (2.0pt);
\fill [color=ffqqqq] (9,1.02) circle (2.0pt);
\end{tikzpicture} 
\caption{Construction of an electrical network.}
\label{construction}
\end{figure}
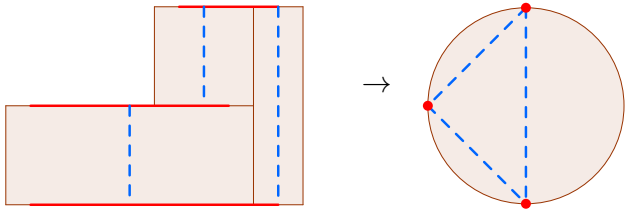

\begin{proof}[Proof of Lemma~\ref{correspondence-general}] (\ref{c-g-1})$\implies$(\ref{c-g-2}). %Denote by $1,2,\dots,b$ all the vertical sides of the polygon $P$. Denote by $b+1,\dots,n$ all the maximal vertical cuts in the tiling.
%Let us construct a graph with vertices $1,\dots,n$ in the plane as follows.
Take a ``generic'' polygon $P$ tiled by rectangles.
%It is convenient to consider it as a subset of the following auxiliary polygon $\hat P$. For each vertical side $I$ of $P$ take an isosceles triangle $\hat I$ with base $I$ and a small height such that $P\cap \hat I=I$. Set $\hat P$ to be the union of $P$ and all the triangles $\hat I$.

Let us construct the graph of the required network; see Figure~\ref{construction}.
%Take the union of the drawn segments and the horizontal sides of all rectangles $P_e$.
%Remove the ''bottom'' horizontal sides of the polygon $P$.
Consider the union of the horizontal sides of all the rectangles in the tiling. This union splits into several disjoint segments called \emph{horizontal cuts}. Paint red (bold) each horizontal cut except small neighborhoods of its endpoints. Paint blue (dashed) the vertical centerline of each rectangle in the tiling.

Contract all red segments. Then the blue set ``becomes'' a graph $G$ and the polygon $P$ ``becomes'' a topological disc $D$
(since the $y$-coordinates of the horizontal sides of $P$ are distinct it follows that the intersection of each red segment  with the boundary $\partial P$ is connected).
Denote by $1,\dots,b$ the vertices of the graph $G$
belonging to the boundary $\partial D$ and denote by $b+1,\dots,n$ all the other vertices.
%Clearly, $G\subset D$, $G\cap \partial D=\{1,\dots,b\}$ and
Clearly, each connected component of $G$ contains a boundary vertex. Thus $G$ is a graph of a planar network.

Let us define voltages, currents, and conductances in the network. For each vertex $k=1,\dots,n$ of the graph $G$ set $U_k$ to be the $y$-coordinate of the horizontal red segment contracted to the vertex. For each edge $kl$ of the graph $G$, obtained from the vertical centerline of a rectangle in the tiling, set $I_{kl}$ and $c_{kl}$ to be the horizontal side (with an appropriate sign) and the ratio of the rectangle, respectively. The laws (C), (I) are now checked directly. The constructed network is the required network.

\begin{figure}[htbp]
\definecolor{qqwwff}{rgb}{0,0.4,1}
\definecolor{zzttqq}{rgb}{0.6,0.2,0}
\begin{tikzpicture}[line cap=round,line join=round,>=triangle 45,x=0.7cm,y=0.7cm]
\clip(0,-1.5) rectangle (15.85,6.21);
\fill[color=zzttqq,fill=zzttqq,fill opacity=0.1] (8.5,3) -- (11.5,3) -- (11.5,5) -- (14.5,5) -- (14.5,1) -- (8.5,1) -- cycle;
\draw [color=zzttqq] (8.5,3)-- (11.5,3);
\draw [color=zzttqq] (11.5,3)-- (11.5,5);
\draw [color=zzttqq] (11.5,5)-- (14.5,5);
\draw [color=zzttqq] (14.5,5)-- (14.5,1);
\draw [color=zzttqq] (14.5,1)-- (8.5,1);
\draw [color=zzttqq] (8.5,1)-- (8.5,3);
\draw [color=zzttqq] (11.5,3)-- (13.5,3);
\draw [color=zzttqq] (13.5,5)-- (13.5,1);
\draw (5.82,3.51) node[anchor=north]
{$\mathbf{ \rightarrow}$};
\draw [color=zzttqq,fill=zzttqq,fill opacity=0.1] (3,3) circle (1.4cm);
\draw [line width=1pt,color=qqwwff] (2.97,5)-- (1,3);
\draw [line width=1pt,color=qqwwff] (2.97,5)-- (3,1);
\draw [line width=1pt,color=qqwwff] (3,1)-- (1,3);
\draw [->] (7.5,0) -- (7.5,6);
\draw [->] (7.5,0) -- (15.5,0);
\draw (6.89,5.44) node[anchor=north] {$U_1$};
\draw (6.89,3.44) node[anchor=north] {$U_2$};
\draw (6.89,1.39) node[anchor=north] {$U_3$};
\draw (8.5,-0.13) node[anchor=north] {$I_{f(3)}$};
\draw (11.51,-0.13) node[anchor=north] {$I_{f(2)}$};
\draw (13.5,-0.13) node[anchor=north] {$I_f$};
\draw (14.51,-0.13) node[anchor=north] {$I_{f(1)}$};
\draw (3.90,3.62) node[anchor=north] {$f(1)$};
\draw (0.55,5.49) node[anchor=north] {$f(2)$};
\draw (0.55,1.77) node[anchor=north] {$f(3)$};
\draw [color=qqwwff](12.60,4.51) node[anchor=north] {$P_{12}$};
\draw [color=qqwwff](11.05,2.58) node[anchor=north] {$P_{23}$};
\draw [color=qqwwff](13.98,3.62) node[anchor=north] {$P_{13}$};
\draw [->] (1.06,4.67) -- (1.65,4.08);
\draw [->] (1.06,1.42) -- (1.67,1.97);
\draw (2.19,3.62) node[anchor=north] {$f$};
\draw[color=zzttqq] (10.60,4.97) node {$P$};
\fill [color=qqwwff] (2.97,5) circle (2.0pt);
\draw[color=qqwwff] (3.17,5.55) node {1};
\fill [color=qqwwff] (3,1) circle (2.0pt);
\draw[color=qqwwff] (3.14,0.49) node {3};
\fill [color=qqwwff] (1,3) circle (2.0pt);
\draw[color=qqwwff] (0.52,3.31) node {2};
\draw [color=black] (7.5,5)-- ++(-1.0pt,0 pt) -- ++(2.0pt,0 pt) ++(-1.0pt,-1.0pt) -- ++(0 pt,2.0pt);
\draw [color=black] (7.5,3)-- ++(-1.0pt,0 pt) -- ++(2.0pt,0 pt) ++(-1.0pt,-1.0pt) -- ++(0 pt,2.0pt);
\draw [color=black] (7.5,1)-- ++(-1.0pt,0 pt) -- ++(2.0pt,0 pt) ++(-1.0pt,-1.0pt) -- ++(0 pt,2.0pt);
\draw [color=black] (8.5,0)-- ++(-1.0pt,0 pt) -- ++(2.0pt,0 pt) ++(-1.0pt,-1.0pt) -- ++(0 pt,2.0pt);
\draw [color=black] (11.5,0)-- ++(-1.0pt,0 pt) -- ++(2.0pt,0 pt) ++(-1.0pt,-1.0pt) -- ++(0 pt,2.0pt);
\draw [color=black] (13.5,0)-- ++(-1.0pt,0 pt) -- ++(2.0pt,0 pt) ++(-1.0pt,-1.0pt) -- ++(0 pt,2.0pt);
\draw [color=black] (14.5,0)-- ++(-1.0pt,0 pt) -- ++(2.0pt,0 pt) ++(-1.0pt,-1.0pt) -- ++(0 pt,2.0pt);
\end{tikzpicture} 
\caption{Construction of a tiling.}
\label{construction2}
\end{figure}
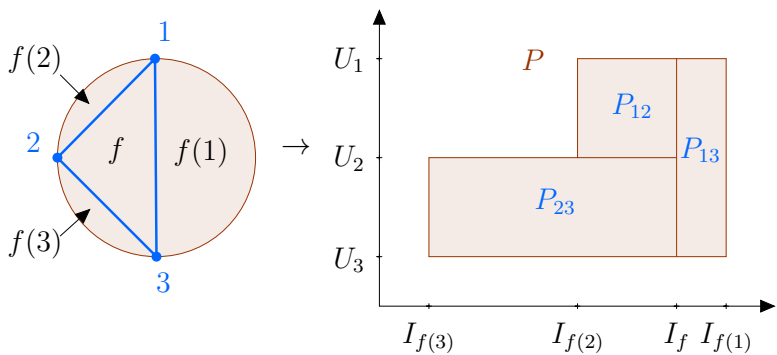

(\ref{c-g-2})$\implies$(\ref{c-g-1}). Take an electrical network as in~(\ref{c-g-2}).
%Suppose that $U_1$, $\dots$, $U_b$ and $I_1, \dots, I_b$ are the incoming voltages and currents, respectively.
%Assume without loss of generality that the minimal $x$-coordinate of the points of the polygon $P$ is $0$.
%Let $P$ be the rectangle with the vertices
%$(0,U_1), (0,U_2), (I_1,U_2), (I_1,U_1)$.
Construct a tiling of $P$ as follows.

Let $e$ be an edge of the network. Denote by $e\uparrow$ ($e\downarrow$) the endpoint of $e$ with higher (lower) voltage (it is well-defined because all the edges are essential). By  \emph{a face} we mean a connected component of the complement to the network in the unit disc $D$.
Denote by $e\leftarrow$ ($e\rightarrow$) the
face %connected component of the complement to the network in the disc,
that borders the edge $e$ from the left-hand (right-hand) side while one moves along the edge $e$ from $e\uparrow$ to $e\downarrow$.

By law (I) it follows that to each face $f$ one can assign a number $I_f$ in such a way that $I_{kl\rightarrow}-I_{kl\leftarrow}=|I_{kl}|$. Without loss of generality assume $\min_{f}I_f=\min_{(x,y)\in P}x$, where the minimum in the left-hand side is over all the faces $f$ meeting $\partial D$.

Let $P_e$ be the rectangle with the vertices $(I_{e\leftarrow},U_{e\uparrow}), (I_{e\leftarrow},U_{e\downarrow}),
(I_{e\rightarrow},U_{e\uparrow}), (I_{e\rightarrow},U_{e\downarrow})$. The rectangles $P_e$, where $e$ runs through all the edges of the network, tile the polygon $P$ by the following two lemmas (the rectangles $P_e$ cover $P$ by Lemma~\ref{cla} and do not overlap by Lemma~\ref{clb}).
\end{proof}
%It suffices to prove the following:
%\begin{enumerate}[(a)]
%\item $\partial \bigcup P_e\subset\partial P$;
%\item $\sum_e\mathop{Area}(P_e)=\mathop{Area}(P)$.
%\end{enumerate}
%Indeed,

\begin{lemma}\label{cla} $\bigcup_e P_e=P$.
\end{lemma}

\begin{proof} It suffices to prove that $\partial \bigcup_e P_e\subset\partial P$. Since $\partial P$ is connected and $\bigcup_e P_e$ is bounded, the lemma will follow.

%Add a vertex $0$ lying outside the disc $D$ to the network. Join the vertex $0$ with the boundary vertices by $b$ \emph{new} edges. %Set $U_0=0$ and define $I_{0k}=I_k$ to be the current flowing inside the network through a boundary vertex $k$.

We need the following description of the boundary $\partial P$; see Figure~\ref{construction2}.
Boundary vertices split $\partial D$ into $b$ arcs. Start from vertex $b$ and move along the circle $\partial D$
counterclockwise. Enumerate the arcs in the order they appear in the motion.
Denote by $f(v)$ the face containing the arc $v$. Denote by $H_v$ the segment joining the points $(I_{f(v)},U_v)$ and $(I_{f({v+1})},U_v)$, where we set $f(b+1)=f(1)$.
Denote by $V_v$ the segment joining the points $(I_{f(v)},U_{v-1})$ and $(I_{f({v})},U_{v})$, where we set $U_0=U_b$.
Clearly, $\partial P=\bigcup_{v=1}^{b}(H_v\cup V_v)$.

Take a ``general position'' point $p\in\partial\bigcup_e P_e$, say, in %the interior of
a horizontal side of the ``polygon'' $\bigcup_e P_e$.
Without loss of generality assume that
the point $p$ belongs to the top side of a rectangle $P_e$. Denote by $v=e\uparrow$ the vertex of the edge $e$ of higher voltage.
%Denote by $\mathop{St}v$ the set of edges containing the vertex $v$.
%Denote by $S$ the set of all the rectangles $P_d$ such that $v$ is an endpoint of the edge $d$.

Draw a horizontal line $H$ through the top side of the rectangle $P_e$. %Consider all the rectangles $P_d$ such that $d\in\mathop{St}v$.
We say that a rectangle $P_d$ is \emph{adjacent} if the vertex $v$ is an endpoint of the edge $d$.
Adjacent rectangles border upon the line $H$ either from above or from below.

First assume that the vertex $v$ is nonboundary. Let us show that each point of $H$ (besides a finite set) is bordered by the same number of adjacent rectangles from above and from below. Indeed, the intersection of an adjacent rectangle $P_d$ and the line $H$ is the \emph{touching} segment with the endpoints points $(I_{d\leftarrow},U_v)$ and $(I_{d\rightarrow},U_v)$. Orient the segment from the left to the right (the right to the left), if the rectangle $P_d$ borders upon the line $H$ from above (below). Then the ``head'' of each touching segment is a ``tail'' of another touching segment. Hence
almost each point of the line $H$ belongs to the same number of touching segments of each orientation.

Now, since the rectangle $P_e$ borders upon the point $p$ from below and $p$ is in ``general position'' it follows that some adjacent rectangle $P_d$ borders upon it from above. Thus $p$ belongs to $\Interior (P_e\cup P_d)\subset \Interior \bigcup_e P_e$, a contradiction.

So $v$ is a boundary vertex. Analogously to the above, each point of $H-H_v$ (besides a finite set) is bordered by the same number of adjacent rectangles from above and from below. Hence $p\in H_v$ and thus $p\in\partial P$.
\end{proof}
%Assume that say $I_{e\leftarrow}<x<I_{e\rightarrow}$, $y=U_{e\uparrow}$. Let $u=e\uparrow$. Let us show that $u$ is a boundary vertex and $I_{u \leftarrow}<x<I_{u\rightarrow}$.

%Let $e=e_1,\dots,e_p$ be all the edges for which $u=e_k\uparrow$ and $e_{p+1},\dots,e_q$ be all the edges for which $u=e_k\downarrow$. Consider \emph{positive} segments $[I_{e_k\leftarrow},I_{e_k\rightarrow}]$, $1\le k\le p$, and \emph{negative} segments $[I_{e_k\rightarrow},I_{e_k\leftarrow}]$, $p+1\le k\le q$.

%Assume that $u$ is nonboundary. Then clearly any point of the line $\mathbb{R}$ is covered by the same number of positive and negative segments. In particular, there is a negative segment $[I_{e_k\rightarrow},I_{e_k\leftarrow}]$ containing $x$. Then the point $(x,y)$ belongs to the horizontal side of the rectangle $P_{e_k}$ and lies in the interior of $P_e\cup P_{e_k}$, a contradiction.
%So $u$ is a boundary vertex.

%Assume $x\notin [I_{u \leftarrow},I_{u\rightarrow}]$. Still any point of the set $\mathbb{R}-[I_{u \leftarrow},I_{u\rightarrow}]$ is covered by the same number of positive and negative segments. Thus we again get a contradiction. Hence $I_{u \leftarrow}<x<I_{u\rightarrow}$. So $(x,y)$ belongs to a horizontal side of $P$.

%Analogously one can prove that a point $(x,y)$ in the interior of a vertical side of the polygon $\bigcup_e P_e$ necessarilylies in a vertical side of $P$. This proves the claim.

\begin{lemma}\label{clb}
$\sum_e\Area(P_e)=\Area(P)$.
\end{lemma}

\begin{proof} This follows immediately from Lemma~\ref{energy} because $\Area(P_{kl})=(U_k-U_l)I_{kl}$ and $\Area(P)=\sum_{1\le u\le b}U_uI_u$.
\end{proof}

\subsection{Inverse problems}\label{basicposreal}

First let us present a proof of Lemma~\ref{PosReal} following~\cite{Cauer}. For generalizations to the case $b>2$ see \cite{Foster} and references therein.

\begin{proof}[Proof of Lemma~\ref{PosReal}]
(\ref{preserve})$\implies$ (\ref{roots}). Indeed, if $\Real z\le 0$ then $\Real C(z)=-\Real C(-z)\le 0$ and thus $C(z)\ne 1$.
%So (\ref{roots}) holds.
%\cite[Lemma 4]{FLR} Let $r ( z ) = C(z) + C ( - z )$. We wish to show
%that $r ( z )$ is identically zero. To begin, note that $r(z)$ takes the imaginary axis into itself.
%Thus $r ( i y ) = C(iy) + C ( - i y ) = C(iy) +  C(\overline{iy}) = C ( i y ) + \overline{C(iy)} = 0$. Since the
%numerator of $r ( z )$ has an infinite number of roots, the numerator is identically zero,
%and hence so is $r(z)$. Finally, suppose $C(z) = x$ has a root with real part less than or
%equal to zero. If this root has real part less than zero, then~(\ref{preserve}) is treated. If this root is
%on the imaginary axis, then by continuity we can pick $z$ sufficiently close to the root but
%with negative real part so that $C(z)$ is close enough to $x$ so that $\Real (C(z)) > 0$, also
%violating~(\ref{preserve}).

(\ref{roots})$\implies$(\ref{preserve}). %Since $C(z)$ is odd it follows that all roots of $C(z)=-1$ have negative real parts.
Consider the equation $C(z)=w$. Move $w$ continuously in the half-plane $\Real w>0$. The roots cannot cross the line $\Real z=0$ (because $\Real z=0$ implies $\Real C(z)=0$ for an odd function $C(z)\in\mathbb{R}(z)$). Thus for each $w$ in the half-plane $\Real w>0$ all roots of $C(z)=w$ are in the half-plane $\Real z>0$. Since $C(z)$ is odd it follows that the same is true for the half-planes $\Real w<0$, $\Real z<0$. So (\ref{preserve}) holds.

(\ref{preserve})$\implies$(\ref{derivatives}). %Since $C(z)$ is odd it follows that $\Real z<0$ implies $\Real C(z)<0$.
Suppose that $C(z)=0$. %, where $z\in \mathbb{C}$.
Then $\Real z=0$
because $\Real z>0\implies\Real C(z)>0$ and $\Real z<0\implies\Real C(z)=-\Real C(-z)<0$. Since implication~(\ref{preserve}) and its converse hold in a neighborhood of the point $z$ it follows that $C'(z)>0$.
%A simple limiting argument proves the same for $z=\infty$.

(\ref{derivatives})$\implies$(\ref{factorization})
%Assume for simplicity that $C(\infty)\ne0$.
Let $z_1,\dots,z_m$ be the roots of $C(z)$. Since $C'(z_k)>0$ it follows that the roots are simple. Thus $C(z)$ has not more than $m$ poles. The roots split the projective line $\Real z=0$ into $m$ ``segments''. Since $C'(z_k)>0$ it follows that for sufficiently small $\epsilon>0$ we have $i\cdot C(z_k-i\epsilon)>0$ and $i\cdot C(z_k+i\epsilon)<0$. By the intermediate value theorem it follows that each ``segment'' contains a pole of $C(z)$. Thus all the  $m$ poles of $C(z)$ belong to the projective line $\Real z=0$ and alternate with the roots. So~(\ref{factorization}) holds.
%Thus $C(z)$ satisfies condition~(\ref{factorization}).

(\ref{factorization})$\implies$(\ref{continued}). %The \emph{height} of a rational function is the sum of the degrees of the nominator and the denominator.
Denote by $\mathop{ht}C(z)$ the sum of the degrees of the nominator and the denominator of $C(z)$.
The proof is by induction over $\mathop{ht}C(z)$. If $\mathop{ht}C(z)=1$ then there is nothing to prove. Now assume that $n\ge 1$ and, say, $b_n\ne 0$ in condition~(\ref{factorization}).
%Let
%$$
%C(z)=\frac{\prod\limits_{k=1}^n
%(z^2+b_k^2)}{d_1z\prod\limits_{k=1}^n(z^2+a_k^2)}.
%$$
%(The other alternatives from condition (\ref{factorization}) are considered similarly.)

Denote by $r(z)=1/(C(z)-d_1z)$ and $q(z)=1/C(z)$. Let us prove that $r(z)$ satisfies condition~(\ref{derivatives}).
Indeed, clearly $\lim_{z\to\infty}r(z)\ne 0$ and $r(z)$ is odd.
The roots of $r(z)$ are the numbers $\pm ib_1,\dots,\pm ib_n$.
For each $l=1,\dots,n$
$$r'(\pm ib_l)=q'(\pm ib_l)=\frac{2}{d_1(a_l^2-b_l^2)}\prod\limits_{k\ne l}\frac{b_k^2-b_l^2}{a_k^2-b_l^2}>0
$$
by the condition $a_1>b_1>a_2>\dots>b_n\ge 0$.

Hence by Lemma~\ref{PosReal}(\ref{derivatives})$\implies$(\ref{factorization}) it follows that
$r(z)$ satisfies condition~(\ref{factorization}) as well. On the other hand, $\mathop{ht}r(z)<\mathop{ht}C(z)$. By the inductive hypothesis, $r(z)$ satisfies condition~(\ref{continued}).
Thus
$
C(z)=d_1z+1/r(z)
$
also satisfies condition~(\ref{continued}).

%(\ref{continued})$\implies$(\ref{series-parallel-network}). The construction of the required network is straightforward, cf.~\S\ref{prooflfrs}.

(\ref{continued})$\implies$(\ref{preserve}).
This is proved by induction over $m$.
%This follows from Lemma~\ref{electrical}(\ref{el-5}).
\end{proof}

\begin{proof}[Proof of Theorem~\ref{th-foster}]
(\ref{th-foster-1})$\implies$(\ref{th-foster-2}). This immediately follows from the definitions; see~\S\ref{alternating}.

(\ref{th-foster-2})$\implies$(\ref{th-foster-3}). This follows from Lemma~\ref{electrical}(\ref{el-2}) and~(\ref{el-5}).

(\ref{th-foster-3})$\implies$(\ref{th-foster-1}).
First assume $\lim_{z\to\infty}C(z)\ne0$. Then by
Lemma~\ref{PosReal}(\ref{preserve})$\implies$(\ref{continued})
condition~(\ref{continued}) holds.
Then the required series-parallel network with edges of conductances $d_1z$, $1/d_2z$, $d_3z$, $1/d_4z$, \dots, $(d_mz)^{(-1)^m}$ is constructed directly; cf.~\S\ref{prooflfrs}.

Now assume $\lim_{z\to\infty}C(z)=0$. Apply
Lemma~\ref{PosReal}(\ref{preserve})$\implies$(\ref{continued})
to the function $1/C(z)$. Then the function $C(z)$ satisfies condition~(\ref{continued}) of Lemma~\ref{PosReal} with the only difference that $d_1=0$.
The construction of the previous paragraph still leads to the required network.
\end{proof}

\begin{remark}\label{rationality} If $C(z)\in\mathbb{Q}(z)$ then we may assume the network in condition~(\ref{th-foster-2}) of Theorem~\ref{th-foster}
is series-parallel and each number $d_j$ is rational. %Indeed, by The Euclidean algorithm it follows that all $d_j\in\mathbb{Q}$ in expression~(\ref{continued}) of Lemma~\ref{PosReal}.
\end{remark}

\begin{proof}[Proof of Corollary~\ref{series-parallel}] %Let $C(c_1,c_2)$ be the conductance of the given network. Then by Lemma~\ref{electrical}(2) we get $C(c_1,c_2)=\sqrt{c_1c_2}\, C(z,1/z)$, where $z=\sqrt{c_1/c_2}$.
By condition~(\ref{th2-2}) of Theorem~\ref{th2} we have $C(c_1,c_2)=\sqrt{c_1c_2}\, C(z)$, where $C(z)=C(z,1/z)$ and $z=\sqrt{c_1/c_2}$.
By Theorem~\ref{th-foster}(\ref{th-foster-3})$\implies$(\ref{th-foster-1})
and Remark~\ref{rationality} it follows that
there is a series-parallel network of conductance $C(z)$ such that the conductance of each edge $j$
is either $d_jz$ or $1/d_jz$ for some $d_1,\dots,d_m\in\mathbb{Q}$. Replace each edge by an appropriate series-parallel subnetwork to get a new network with edge conductances $z$ and $1/z$. Multiplying the conductance of each edge by $\sqrt{c_1c_2}$ we get the required network.
%By Lemma~\ref{PosReal}(\ref{preserve})$\implies$(\ref{continued}) it satisfies condition (\ref{continued}).
%Therefore, say, for $m$ even and $C(0)=0$,
%$$
%C(c_1,c_2)=d_1c_1+\cfrac1{\cfrac{d_2}{c_2} +\cfrac1{d_3 c_1+\dots+\cfrac{d_m}{ c_2}}}.
%$$
%All the numbers $d_k\in \mathbb{Q}$ by the Euclidean algorithm. Now the required series-parallel network is constructed analogously to the proof of Theorem~\ref{lfrs}(\ref{lfrs-3})$\implies$(\ref{lfrs-1}).
\end{proof}

\section{Proof of variations}\label{finalproofs}

\subsection{Proof of Theorem \ref{squaring-polygons}} \label{proof-general}

%The proof of Theorem~\ref{squaring-polygons} follows the idea of \cite[\S8.1]{K} and uses the following results of \cite{CGV}.

%\begin{theorem} \label{CGV5} \textup{\cite[Theorem~5]{CGV}} The set $-\bar\Omega_b$ is the set of all possible responses of planar networks with $b$ boundary vertices and positive edge conductances.
%\end{theorem}

%\begin{theorem} \label{CGV3} \textup{\cite[Theorem~3]{CGV}} Let $G$ be a network with a given response $C_{uv}$ having minimal number of edges. Consider the response of $G$ as a matrix-valued function in edge conductances $c_{kl}$. Then this function is injective.
%\end{theorem}

\begin{proof}[Proof of Theorem~\ref{squaring-polygons}] (\ref{s-p-1})$\implies$(\ref{s-p-2}).
Let the polygon $P$ be tiled by squares.
%rectangles of ratios $c_1, \dots,c_n$.
By Lemma~\ref{correspondence-general} there is a planar electrical circuit with edge conductances $1$,
%$c_1, \dots,c_n$,
incoming voltages $U_1$, \dots, $U_b$, and incoming currents $I_1$, \dots, $I_b$. Let
$C_{uv}$ be the response of the circuit. Then $I_v=\sum C_{uv} U_u$. By Lemma~\ref{responseclaim}(\ref{rc-1}) all the entries of $C_{uv}$ are rational. By Theorem~\ref{CGV5} we have $C_{uv}\in \Omega_b$.
%\end{proof}

(\ref{s-p-2})$\implies$(\ref{s-p-1}). Let $C_{uv}\in \Omega_b$ be a matrix with rational entries such that $I_v=\sum C_{uv} U_u$. By Theorem~\ref{CGV5} there are planar electrical networks with the response $C_{uv}$. Take a minimal network with this property. By Theorem~\ref{inverse}
the conductance of each edge of the network is rational.
%Let us prove that all the conductances $c_{kl}\in\mathbb{Q}$. Indeed, otherwise perform a nontrivial automorphism of the field $\mathbb Q(c_{12},\dots,c_{n,n-1})$. The obtained network has the same graph and the same response by Lemma~\ref{responseclaim}(1). This contradicts to minimality of the number of edges in the circuit by Theorem~\ref{CGV3}. Thus all $c_{kl}\in\mathbb Q$.
Set the incoming voltages to be $U_1$, \dots, $U_b$. Then the incoming currents are $I_1$, \dots, $I_b$. Delete all the unessential edges from the circuit. By Lemma~\ref{correspondence-general} it follows that the polygon $P$ can be tiled by rectangles of rational ratio, and hence by squares.
\end{proof}

%\begin{remark} The converse theorem would follow from Lemma~\ref{correspondence-general}, \cite[Theorem~5]{CGV} and the following conjecture: \emph{if the response of a network with real edge conductances is rational then there is a network with rational edge conductances and the same response}.
%\end{remark}

\begin{corollary} \label{mostgeneral} \textup{(of Lemmas~\ref{responseclaim}, \ref{correspondence-general}, and Theorem~\ref{CGV5})} Let $P$ be a ``generic'' orthogonal polygon $P$ with $b$ horizontal sides having
signed length $I_1,\dots,I_b$ and $y$-coordinates $U_1,\dots,U_b$. Suppose that the polygon $P$ can be tiled by rectangles of ratios $c_1$, $\dots$, $c_n$. Then there is a function $C_{uv}(z_1,\dots,z_n)$ satisfying conditions~\textup{(\ref{rc-1}), (\ref{rc-2}),} and \textup{(\ref{rc-5})} of Lemma~\ref{responseclaim} such that $C_{uv}(c_1,\dots,c_n)\in \Omega_b$ and $I_v=\sum_{1\le u\le b} C_{uv}(c_1,\dots,c_n)U_u$ for each $v=1,\dots,b$.
\end{corollary}

\subsection{Proof of Theorem~\ref{polygon-square}}

\begin{proof}[Proof of Theorem~\ref{polygon-square}] $\impliedby$. This holds because a polygon with rational vertices can be tiled by squares.

$\implies$. Let $P$ be tiled by rectangles of ratios $c$ and $1/c$. Let us prove analogously to the proof of Theorem~\ref{lfrs}(\ref{lfrs-1})$\implies$(\ref{lfrs-2}) that all algebraic conjugates of $c$ have positive real parts. Then Theorem~\ref{polygon-square} will follow from Theorem~\ref{lfrs}(\ref{lfrs-2})$\implies$(\ref{lfrs-1}).

Consider the circuit given by Lemma \ref{correspondence-general}.
Replace each edge of conductance $c$ (respectively, $1/c$) in the circuit by an edge of conductance $z\in\mathbb{C}$ (respectively, $1/z$). Let $C_{uv}(z)$ be the response of the obtained circuit. Consider the \emph{energy dissipation} function $E(z)=\sum_{1\le u,v\le b} C_{uv}(z)U_u U_v$. Since each $U_u\in\mathbb{Q}$ it follows by Lemma~ \ref{responseclaim}(\ref{rc-1}) that $E(z)\in\mathbb{Q}(z)$. By Lemma~ \ref{responseclaim}(\ref{rc-2}) it follows that
$E(z)$ is odd.
Clearly, $E(c)=\sum_{1\le u\le b} I_u U_u=\Area(P)$. Thus $E(c)\in\mathbb{Q}$ and $E(c)>0$.
%We have $E(c)>0$ by Lemma \ref{responseclaim}(4) (in fact $E(c)$ is the area of $P$).

Since $E(z)\in\mathbb{Q}(z)$, $E(z)$ is nonconstant, and $E(c)\in\mathbb{Q}$ it follows that $c$ is algebraic. Let $z$ be an algebraic conjugate of $c$. Then $E(z)=E(c)>0$.

Let us prove that $\Real z>0$.
Indeed, first assume that $\Real z<0$. Then by Lemma \ref{responseclaim}(\ref{rc-5}) we have $0\le\Real E(-z)=-\Real E(z)<0$, a contradiction. A simple limiting argument shows that assumption $\Real z=0$ also leads to a contradiction.
Thus $\Real z>0$.
\end{proof}

%\begin{proof}[Proof of Theorem~\ref{rectangling-polygons}]
%Suppose that a polygon $P$ is tiled by rectangles of ratios $c_1,\dots,c_n$.
%Consider the electrical circuit given by Lemma~\ref{correspondence-general}.
%For each $k=1,\dots,n$ replace each edge of conductance $c_k$ in the network by an edge of conductance $z_k$. Let $C_{uv}(z_1,\dots,z_n)$ be the response of the obtained network. Then Theorem~\ref{rectangling-polygons}(1)--(4)
%follows from Lemma~\ref{responseclaim} and \cite[Theorem~5]{CGV}.
%\end{proof}

\subsection{Proof of Theorem~\ref{detection}}

\begin{proof}[Proof of Theorem~\ref{detection}] (\ref{det-1})$\implies$(\ref{det-2}). This follows from Lemma~\ref{responseclaim}(\ref{rc-5}) and Remark~\ref{strict}.

(\ref{det-2})$\implies$(\ref{det-1}). For $b=2$ there is nothing to prove. Assume that $b=3$. Let $\delta>0$ be a small number, $r_{uv}=-\Real C_{uv}-\delta$, $m_{uv}=-\Imaginary C_{uv}$. By the almost positively definiteness it follows that
$\Real ${\tiny $\left(\begin{matrix}
C_{11} & C_{12} \\
C_{21} & C_{22}
\end{matrix}\right)$}
is positively definite. Thus
{\tiny $\left(\begin{matrix}
r_{31}+r_{12} & -r_{12} \\
-r_{12}       & r_{12}+r_{23}
\end{matrix}\right)$}
is positively definite for sufficiently small $\delta$. Hence $r_{12}+r_{23}, r_{31}+r_{12}, r_{12}r_{23}+r_{23}r_{31}+r_{31}r_{12}> 0$. Analogously $r_{23}+r_{31}>0$. Thus at least two of the numbers $r_{12}, r_{23}, r_{31}$ are positive.

If $r_{12}, r_{23}, r_{31}>0$ then the required network is a triangle with vertices $1$, $2$, $3$ and with edge conductances $c_{kl}=r_{kl}+im_{kl}+\delta$.

Now assume that exactly one of the numbers $r_{12}, r_{23}, r_{31}$, say, $r_{31}$ is nonpositive. Take a large number $M$ and denote by $\Delta_M=r_{12}r_{23}+r_{23}r_{31}+r_{31}r_{12}+iM(r_{23}+r_{12})$.
The required network is a complete graph on the vertices $1$, $2$, $3$, $4$ with edge conductances
$$
\begin{tabular}{llll}
$c_{12}=$& $im_{12} +\delta$,
\quad & $c_{14}=$& $\Delta_M/r_{23}$,\\
$c_{23}=$& $im_{23} +\delta$,
\quad & $c_{34}=$& $\Delta_M/r_{12}$,\\
$c_{31}=$& $im_{31} +\delta-iM$,
\quad & $c_{24}=$& $\Delta_M/(r_{31}+iM)$.
\end{tabular}
$$
Clearly, for $M^2>(r_{12}r_{23}+r_{23}r_{31}+r_{31}r_{12})|r_{31}|/(r_{23}+r_{12})$ we have each $\Real  c_{kl}> 0$.

Let us show by an electrical transformation that the network has response $C_{uv}$. Replace the ``letter $Y$'' formed by the edges $14$, $24$, and $34$ by a ``triangle $\Delta$'' formed by $3$ new edges of conductances $c'_{12}=r_{12}$, $c'_{23}=r_{23}$, and $c'_{31}=r_{31}+iM$. This \emph{$Y\Delta$-transformation} does not change the response \cite[page 12]{K}. The obtained network has $3$ pairs of multiple edges. Thus it has the same response as a triangle with edge conductances $r_{12}+im_{12}+\delta$, $r_{23}+im_{23}+\delta$, $r_{31}+im_{31}+\delta$. So the network has the response $C_{uv}$.
\end{proof}

\section*{Acknowledgements} The authors are grateful to A.~Akopyan, M.~Houston, and B.~George for useful discussions.
M.~Skopenkov was supported in part by INTAS grant 06-1000014-6277,
Russian Foundation of Basic Research grant
06-01-72551-NCNIL-a,
%07-01-00648-a,
%President of the Russian Federation grant NSh-4578.2006.1,
%Agency for Education and Science grant RNP-2.1.1.7988,
Moebius Contest Foundation for Young Scientists and Euler Foundation.

%\section*{References}

\bibliographystyle{elsarticle-num}

\begin{thebibliography}{99}

\bibitem{BeSch} I.~Benjamini and O.~Schramm, \textrm{Random walks and harmonic functions on infinite planar graphs using square tilings}, Ann.~Prob.~\textbf{24:3} (1996), 1219--1238.

%\bibitem{Bol } B. Bollobas, \textrm{Graph theory: an introductory course}, Springer, 1979.

\bibitem{B} V.~G. Boltianskii, \textrm{Hilbert's Third Problem (trans.~by R.~Silverman)}, V.~H.~Winston \& Sons, Washington D.~C., 1978.

\bibitem{BSST} R.~L.~Brooks, C.~A.~B.~Smith, A.~H.~Stone, and W.~T.~Tutte, \textrm{The dissection of rectangles into squares}, Duke Math.~J.~\textbf{7} (1940), 312--340.

\bibitem{Cal} A.~P.~Calder\'on, \textrm{On an inverse boundary value problem},
    Comput.~Appl.~Math.~\textbf{25:2-3} (2006), 133--138. %\url{http://www.scielo.br/scielo.php?pid=S0101-82052006000200002&script=sci_arttext&tlng=en}

\bibitem{CFP} J.~Cannon, W.~Floyd, W.~Parry, \textrm{Squaring rectangles: the finite Riemann mapping theorem},
    Contemp.~Math.~\textbf{169} (1994), 133--211.

\bibitem{Cauer} W.~Cauer, \textrm{Die Verwirklichung der Wechselstromwiderst \"ande vorgeschriebener Frequenzabh \"angigkeit}, Archiv f\"ur Elektrotechnik \textbf{17} (1926), 355--388 (in German).

\bibitem{C} Y.~Colin de Verdi\`ere, \textrm{R\'eseaux \'electriques planaires I},
    %\textrm{Planar resistor networks I},
    Comm.~Math.~Helv.~\textbf{69:1} (1994), 351--374 (in French).


\bibitem{CGV} Y.~Colin de Verdi\`ere, I.~Gitler, and
    D.~ Vertigan, \textrm{R\'eseaux \'electriques planaires II},
    %\textrm{Planar resistor networks II},
    Comm.~Math.~Helv.~\textbf{71:1} (1996), 144--167 (in French).

\bibitem{CIM} E.~B.~Curtis, D.~Ingerman, and~J.~A.~Morrow, \textrm{Circular planar graphs and resistor networks}, Lin.~Alg.~Appl.~\textbf{283} (1998), 115--150.

\bibitem{CM} E.~B.~Curtis and J.~A.~Morrow, \textrm{Inverse problems for electrical networks}, Series on Appl. Math. \textbf{13}, World Scientific, Singapore, 2000.

\bibitem{D} M.~Dehn, \textrm{\"Uber die Zerlegung von Rechtecken in Rechtecke}, Math.~Ann.~\textbf{57} (1903), 314--332 (in German).

\bibitem{DS} P.~G.~Doyle and J.~L.~Snell,
  \textrm{Random walks and electric networks},
  Mathematical Association of America, 1984,
  \url{http://arxiv.org/abs/math.PR/0001057}.

%\bibitem{Duffin} R.~J.~Duffin, \textrm{Elementary operations which generate network matrices}, Proc.~AMS \textbf{6:3} (1955), 335--339.

\bibitem{Dui92} A.~J.~W.~Duijvestijn, \textrm{Simple perfect squared square of lowest order},
    J.~Comb.~Theory B \textbf{25} (1978), 240--243.

%C. J. Bouwkamp and A. J. W. Duijvestijn, \textrm{Catalogue of Simple Perfect Squared Squares of Orders $21$ through $25$}, Eindhoven Univ. Technology, Dept. of Math., Report 92-WSK-03, 1992.

%\bibitem{BB} I. Feshchenko, D. Radchenko, L. Radzivilovsky and M. Tantsiura, \textrm{Dissecting a brick into bars}, Geometriae Dedicata, DOI	 10.1007/s10711-009-9413-y, \url{http://arXiv.org/abs/math/0809.1883v1}.

\bibitem{Fomin} S.~Fomin, \textrm{Loop-erased walks and total positivity}, Trans.~AMS \textbf{353:9} (2001), 3563--3583.

%\bibitem{Foster} R.~M.~Foster, \textrm{A reactance theorem}, Bell System Techn.~J.~\textbf{3} (1924), 259--267.

\bibitem{Foster} R.~M.~Foster, \textrm{Academic and Theoretical Aspects of Circuit Theory}, Proc. IRE \textbf{50:5} (1962), 866--871.

\bibitem{FHTW} C.~Freiling, R.~Hunter, C.~Turner, and
    R.~ Wheeler, \textrm{Tiling with Squares and Anti-Squares}, Amer.~Math.~Monthly \textbf{107:3} (2000), 195--204.

\bibitem{FLR} C.~Freiling, M.~Laczkovich, and D.~Rinne,
    \textrm{Rectangling a rectangle},
    Discr.~Comp.~Geom.~\textbf{17} (1997), 217-225.

\bibitem{FR} C.~Freiling and D.~Rinne,
    \textrm{Tiling a square with similar rectangles},
    Math.~Res.~Lett.~\textbf{1} (1994), 547--558.

%\bibitem{GKLU} A. Greenleaf, Y. Kurylev, M. Lassas, G. Uhlmann, \textrm{Invisibility and Inverse Problems}, \url{http://arXiv.org.0810.0263v1}

%\bibitem{G} M. Gardner, \textrm{Squaring the square}. In: \textrm{The 2nd Scientific American Book of Mathematical Puzzles and Diversions}, University of Chicago Press, 1987, 256 p.

\bibitem{Ha57} H.~Hadwiger, \textrm{Vorlesungen Uber Inhalt, Oberfl\"ahe und Isoperimetrie}, Springer--Verlag, 1957 (in German).

%\bibitem{HH} F.V. Henle, J.M. Henle, \textrm{Squaring the plane}, Amer. Math. Monthly \textbf{115:1} (2008), 3--12.

\bibitem{KeKi} K.~Keating and J.~L.~King, \textrm{Shape tiling}, Elect.~J.~Comb.~\textbf{4:2} (1997), R12.

\bibitem{KK} K.~Keating and J.~L.~King, \textrm{Signed tilings with squares}, J.~Comb.~Theory A~\textbf{85:1} (1999), 83--91.

\bibitem{K} R.~Kenyon, \textrm{Tilings and discrete Dirichlet problems}, Israel J.~Math.~\textbf{105:1} (1998), 61--84.

\bibitem{LS} M.~Laczkovich and G.~Szekeres,
    \textrm{Tiling of the square with similar rectangles},
    Discr.~Comp.~Geom.~\textbf{13} (1995), 569--572.

\bibitem{LSy} G.~F.~Lawler and J.~Sylvester, \textrm{Determining resistances from boundary measurements in finite networks},
    SIAM J.~Discr.~Math.~\textbf{2} (1989), 211--239.

\bibitem{Lovasz} L.~Lovasz, \textrm{Random walks on graphs: a survey}, In: Combinatorics, Paul Erdos is Eighty, D.~Milos, V.~T.~Sos, and T.~Szony, Eds., Budapest, Hungary: Janos Bolyai Math.~Soc., 1996, 353--398.

%\bibitem{Na} A.I. Nachman, \textrm{Reconstructions From Boundary Measurements}, Ann. Math. Second Series \textbf{128:3} (1988), 531--576.

\bibitem{P} V.~G.~Pokrovskii, \textrm{Slicings of $n$-dimensional parallelepipeds}, Math.~Notes \textbf{33:2} (1983), 137--140.
%Matematicheskie Zametki \textbf{33:2} (1983), 273--280 (in Russian).

\bibitem{Sharov} F. Sharov, \textrm{Dissection of a rectangle into rectangles with given side ratios}, Mat. Prosveschenie 3rd ser. \textbf{20} (2016), 200--214 (in Russian) \url{http://arxiv.org/abs/1604.00316} (in English and in Russian).

\bibitem{PrSk} M.~Skopenkov, M.~Prasolov and S.~Dorichenko, \textrm{Dissections of a metal rectangle}, Kvant \textbf{3} (2011), 10--16 (in Russian) \url{http://arxiv.org/abs/1011.3180}.

\bibitem{SMD} M.~Skopenkov, O.~Malinovskaya and S.~Dorichenko, \textrm{Compose a square}, Kvant \textbf{2} (2015), 6--11 (in Russian) \url{http://arxiv.org/abs/1305.2598}. 

%\bibitem{Rayleigh} J.W.S. Rayleigh, \textrm{On the theory of resonance}, In: Collected scientific papers \textbf{1} (1899), 33--75.

%Richards, A special class of functions...

\bibitem{SU} J.~Sylvester and G.~Uhlmann, \textrm{A global uniqueness theorem for an inverse boundary value problem}, Ann.~Math.~\textbf{125} (1987), 153--169.

\bibitem{SuDi} Z.~Su and R.~Ding, \textrm{Tilings of orthogonal polygons with similar rectangles or triangles}, J.~Appl.~Math.~Comp.~\textbf{17:1} (2005), 343--350.

\bibitem{Szegedy} B.~Szegedy, \textrm{Tilings of the square with similar right triangles}, Combinatorica \textbf{21:1} (2001), 139--144.

\bibitem{Wall} H.~S.~Wall, \textrm{Analytic theory of continued fractions}, Chelsea Pub.~Co., Bronx,
    N.~Y., 1967, 433 p.

\bibitem{Weyl} H.~Weyl, \textrm{Repartici\'on de corriente en uno red conductora},
    Rev.~Mat.~Hisp.~Amer.~\textbf{5} (1923),
    153--164 (in Spanish).

\bibitem{Ya68} I.~M.~Yaglom, \textrm{How to dissect a square?} Mathematicheskaya bibliotechka, Nauka, Moscow, 1968, 112 p. (in Russian), \url{http://ilib.mirror1.mccme.ru/djvu/yaglom/square.htm}.
\end{thebibliography}

\end{document}

%\renewcommand{\refname}{Список литературы}
%\abstractname{Аннотация}

\begin{frontmatter}

%% Title, authorsи addresses

%%и the form \ead[url] for the home page:

%\title{Tiling by rectanglesи alternating current}
\title{Разрезания на прямоугольники и переменный ток}


%\author[add1]{M. Prasolov}
\author[add1]{М. Прасолов}
\ead{0x00002a@gmail.com}
%\address[add1]{Moscow State University, Faculty of Mechanicsи mathematics, \\
\address[add1]{Механико-математический факультет МГУ им. М.В. Ломоносова, \\
%Leninskie Gory, 1, GSP-1, Moscow, 119991, Russian Federation}
119991, Москва, Ленинские Горы, 1, ГСП-1}

%\author[add2,add3]{M. Skopenkov\corref{cor1}}
\author[add2,add3]{М. Скопенков\corref{cor1}}
\ead{skopenkov@rambler.ru} %, http://skopenkov.ru} %, Tel.: +966 530 708 319, Fax: +966 2 802 0064}
\cortext[cor1]{Corresponding author}
%\address[add2]{Institute for information transmission problems of the Russian Academy of Sciences,\\
\address[add2]{Институт проблем передачи информации им. А.А. Харкевича РАН,\\
%Bolshoy Karetny per.~19, bld.~1, Moscow, 127994, Russian Federation}
127994, Москва, Большой каретный пер., 19, с.~1}
\address[add3]{King Abdullah University of Science and Technology,\\
P.O.~Box 2187, 4700 KAUST, 23955-6900 Thuwal, Kingdom of Saudi Arabia}

\begin{abstract}
%This paper is on tilings of polygons by rectangles.
%A celebrated physical interpretation of such tilings
%by R.L.~Brooks, C.A.B.~Smith, A.H.~Stoneи W.T.~Tutte
%uses direct-current circuits.
%The new approach of this paper is an application of alternating-current circuits.
%The following results are obtained:
Эта статья посвящена разрезаниям многоугольников на прямоугольники. Знаменитая физическая интерпретация таких разрезаний, принадлежащая Р.~Л.~Бруксу, К.~А.~Б.~Смиту, А.~Г.~Стоуну и У.~Т.~Татту, использует электрические цепи постоянного тока. В данной статье описан новый подход к физической интерпретации разрезаний при помощи цепей переменного тока. Получены следующие результаты:
\begin{itemize}
%\item a necessary condition for a rectangle %of given side ratio can be tiled
%to be tilable by rectangles of given shapes;
\item необходимое условие существования разрезания прямоугольника на прямоугольники заданной формы;
%\item a criterion for a rectangle to be tilable by rectangles similar to it but не all homothetic to it;
\item критерий существования разрезания данного прямоугольника на подобные ему же прямоугольники, некоторые из которых не гомотетичны данному;
%\item a criterion for a ``generic'' polygon to be tilable by squares.
\item критерий существования разрезания многоугольника \2общего положения\3 на квадраты.
\end{itemize}
%These results generalize those of C.~Freiling, R.~Kenyon, M.~Laczkovich, D.~Rinne,и G.~Szekeres. %The approach is also applied to discrete harmonic functionsи random walks on graphs.
Эти результаты обобщают результаты Фрайлинга, Кеньёна, Лацковича, Ринна и Секереша.
\end{abstract}

\begin{keyword}
%% At most 6 for J Comb Theory A

%Tiling \sep rectangle \sep orthogonal polygon \sep alternating current \sep electrical network
Разрезание \sep прямоугольник \sep ортогональный многоугольник \sep переменный ток \sep электрическая цепь

%Tiling \sep partition \sep dissection \sep rectangle \sep orthogonal polygon \sep positive real function \sep continued fraction \sep algebraic number \sep electrical network \sep electrical circuit \sep alternating current \sep conductance \sep admittance \sep discrete harmonic function \sep random walk



\MSC[2010] 52C20 \sep 94C05 \sep 31C20 \sep 30C15 \sep 60J10

\end{keyword}

\end{frontmatter}

%\section{Introduction}\label{intro}
\section{Введение}\label{intro}

%This paper is on combinatorial geometryи its relationship with graph theory, electrical circuit theory, random walks.

%A rectangle ${a\times b}$, where ${a$и ${b}$ are integers, can be tiled by ${a\cdot b}$ squares. Thus a rectangle with rational side ratio can be tiled by squares. In 1903 M.~Dehn proved the reciprocal assertion:
Прямоугольник \fo{a\times b}, где \fo{a,b\in\mathbb{N}}, можно разрезать на \fo{a\cdot b} квадратов. Поэтому прямоугольник с рациональным отношением сторон также можно разрезать на квадраты. В 1903~г. М.~Ден доказал и обратное утверждение:
%(see figure~\ref{bookshelf}):

%\begin{theorem} \label{dehn} \textup{\cite{D}} A rectangle can be tiled by squares \textup{(}not necessarily equal\textup{)} еслии only если
%the ratio of two orthogonal sides of the rectangle is rational.
\begin{theorem} \label{dehn} \textup{\cite{D}} Прямоугольник может быть разрезан на квадраты \textup{(}не обязательно равные\textup{)} тогда и только тогда, когда отношение его сторон рационально.
%it has rational side ratio.
\end{theorem}

%\begin{figure}[htbp]
%\begin{tabular}{cc}
%\includegraphics[width=3.5cm]{bookshelf.pdf} & %\includegraphics[width=4.6cm]{squaredsquare.pdf}
%\end{tabular}
%\caption{Examples of tilings by squares from \url{http://www.squaring.net}}
%\label{bookshelf}
%\end{figure}

%Although this assertion is expected, the proof is complicated. After the original proof, many improvements have been made \cite{B,BSST,Ha57,P,Ya68}.
Хоть это утверждение и кажется вполне естественным, его доказательство оказалось довольно сложным. Первоначальное доказательство со временем несколько раз упрощалось \cite{B,BSST,Ha57,P,Ya68}.
%of Hadwiger \cite{Ha57}, Pokrovskii \cite{P}, Boltianskii \cite{B}, Yaglom \cite{Ya68}.

%The most interesting for us is the approach of R.~L.~Brooks, C.~A.~B.~Smith, A.~H.~Stone,и W.~T.~Tutte~\cite{BSST}. To a tiling of a rectangle they assign a direct-current circuit,и тогда deduce Theorem~\ref{dehn} from certain properties of the circuit.
Для нас наиболее интересен подход Р.~Л.~Брукса, К.~А.~Б.~Смита, А.~Г.~Стоуна и У.~Т.~Татта~\cite{BSST}. Разрезанию прямоугольника они сопоставили цепь постоянного тока и доказали теорему~\ref{dehn}, пользуясь свойствами электрических цепей.
%the uniqueness of current distribution in the circuit.
%We present this proof in \S\ref{prelim} as an introduction to our approach.
%They also apply this technique to find a tiling of a square by squares of distinct sizes \cite{Dui92}; see %a story of Gardner
%~\cite{G}и
%\url{http://www.squaring.net} for a surveyи artwork.
Они также применили эту технику к нахождению разрезания квадрата на квадраты с различными сторонами \cite{Dui92} (см. увлекательную историю об этом в~\cite{G}).
%the figure on the front cover of the journal. %(J. Combin. Theory A). 
%An elementary introduction to the approach is given in \cite{PrSk,SMD}.
Элементарное описание их метода см. в~\cite{PrSk,SMD}.

%Later the minimal tiling of this kind was found by a computer search \cite{Dui92}, see \url{http://www.squaring.net} for a surveyи artwork.
%see figure~\ref{bookshelf} to the right.
%In \S\ref{prelim} we describe this physical interpretation in detail.
%and present the physical proof of Theorem \ref{dehn} as an illustration.

%\begin{figure}[htbp]
%\includegraphics[width=4.5cm]{squaredsquare.pdf}
%\caption{A square can be tiled by ${21}$ distinct squares. Number ${21}$ is minimal with this property.}
%\label{squaredsquare}
%\end{figure}

%A classical related question is whether a square can be tiled by distinct squares. The minimal example of such a tiling has been found quite recently by a computer search, see figure~\ref{squaredsquare}.

%\begin{theorem}[Duijstvijn, 2007] A square can be tiled by ${21}$ distinct squares. Number ${21}$ is minimal with this property.
%\end{theorem}

%In this paper we study finite tilings by arbitrary
%rectangles. The sides of rectangles are assumed to be parallel to the coordinate axes, i.~e., either vertical or horizontal. By {\it the ratio} of a rectangle we %always
%mean the length of the horizontal side divided by the length of the vertical side. We study the following %\emph{shape-tiling}
%problem posed in~\cite[p.~218]{FLR}и~\cite[p.~3]{KeKi}:
В этой статье мы рассмотрим разрезания на конечное число произвольных прямоугольников. Предполагается, что стороны прямоугольников параллельны осям координат --- то есть либо вертикальны, либо горизонтальны. \textit{Отношением сторон} прямоугольника мы будем считать отношение горизонтальной стороны к вертикальной. Рассмотрим следующую задачу, поставленную в~\cite[с.~218]{FLR} и~\cite[с.~3]{KeKi}:

%\begin{problem} \label{mainproblem} Which rectangles can be tiled by rectangles of given ratios ${c_1}$, \dots,~${c_n}$?
\begin{problem} \label{mainproblem} Какие прямоугольники можно разрезать на прямоугольники с отношениями сторон ${c_1}$,~\dots,~${c_n}$?
%When a rectangle of ratio $c$ can be tiled by rectangles of given ratios ${c_1}$, ${c_2}$, \dots,~${c_n}$?
\end{problem}

%A related problem of \emph{signed} tilings is solved in~\cite{KeKi}.
Близкая задача о замощениях \emph{\2с учётом знаков\3} решена в~\cite{KeKi}.

%We do не put any restrictions on the number of rectangles in the tilings. %In particular, this number may be distinct from ${n}$. ch of the ratios ${c_1}$, \dots,~${c_n}$ can be used any number of times or не used at all.
Мы не делаем каких-либо ограничений на количество прямоугольников в разрезании, кроме конечности этого числа: каждое из отношений ${c_1}$,~\dots,~${c_n}$ может быть использовано сколько угодно раз или не использовано совсем.

%For ${n=1$и ${c_1=1}$ the question of Problem~\ref{mainproblem} is answered by Theorem~\ref{dehn}.
%%Notice that for each ratio ${c_k}$ we allow several rectangles in the tiling have this ratio.
%A necessary condition for arbitrary ${n}$ was actually proved by M.~Dehn:
%%if a required tiling exists
%if a rectangle of ratio $c$ can be tiled by rectangles of ratios ${c_1}$, \dots,~${c_n$
%тогда $c$ is (the value of) a rational function in ${c_1}$, ${\dots}$, ${c_n}$ with rational coefficients \cite[Lemma~4]{FR}.
Для случая \fo{n=1} и \fo{c_1=1} ответ на вопрос задачи~\ref{mainproblem} даёт теорема~\ref{dehn}. Необходимое условие для произвольного ${n}$ было доказано М.~Деном: если прямоугольник с отношением сторон $c$ можно разрезать на прямоугольники с отношениями сторон ${c_1}$,~\dots,~${c_n}$, то число $c$ представимо как рациональная функция (точнее, как значение этой функции) с рациональными коэффициентами от ${c_1}$, ${\dots}$, ${c_n}$.

%This function depends only on the ``combinatorial structure'' of the tiling. For instance, если a rectangle of ratio $c$ is dissected into ${2}$ rectangles of ratios ${c_1$и ${c_2}$ by a vertical (respectively, horizontal) cut тогда ${c(c_1,c_2)=c_1+c_2}$ (respectively, ${c(c_1,c_2)=\frac{c_1c_2}{c_1+c_2}}$).
%The problem reduces to description of possible functions ${c(c_1,\dots,c_n)}$. By the mentioned physical interpretation this is equivalent to a natural problem:
%%\begin{problem}
%\textit{describe possible formulas ${c(c_1, \dots, c_n)}$ expressing the conductance of a planar direct-current circuit through the conductances ${c_1}$, ${\dots}$, ${c_n}$ of individual resistors}.
%%\end{problem}
Эта функция зависит только от \2комбинаторной структуры\3 разрезания: например, если прямоугольник с отношением сторон $c$ можно разрезать на 2 прямоугольника с отношениями сторон ${c_1}$ и ${c_2}$ вертикальным (или, соответственно, горизонтальным) разрезом, то \fo{c(c_1,c_2)=c_1+c_2} (или, соответственно, ${c(c_1,c_2)=\frac{c_1c_2}{c_1+c_2}}$).

%The main idea of the paper is to apply {\it alternating}-current circuits (equivalently, circuits with complex-valued conductances) to the above problems. Our first result is the following theorem.
Главная идея нашей статьи --- применить цепи \textit{переменного} тока (что эквивалентно цепям постоянного тока с комплексными сопротивлениями) к задаче~\ref{mainproblem}. Первым нашим результатом является следующая теорема.

\begin{theorem}\label{th2}
%Suppose that a rectangle of ratio $c$ can be tiled by rectangles of ratios ${c_1}$,  \dots, ${c_n}$.
%%(in such a way that there is at least one rectangle with each ratio ${R_k}$).
%тогда ${c=C(c_1,  \dots, c_n)}$ for some rational function ${C(z_1,\dots, z_n)}$ having the following properties
Пусть прямоугольник с отношением сторон $c$ можно разрезать на прямоугольники с отношениями сторон ${c_1}$,~\dots,~${c_n}$. Тогда \fo{c=C(c_1, \dots, c_n)} для некоторой рациональной функции \fo{C(z_1,\dots, z_n)}, обладающей следующими свойствами:
\begin{enumerate}[(1)]
%%\item ${C(c_1, c_2, \dots, c_m)=c}$;
%\item \label{th2-1} \emph{rationality of coefficients}\textup{:}
\item \label{th2-1} \emph{рациональность коэффициентов}\textup{:}
%    %$C(z_1, \dots, z_n)}$ has rational coefficients, i.e.,
    ${C(z_1, \dots, z_n)\in\mathbb{Q}(z_1, \dots, z_n)}$\textup{;}
%    %i.e., $c$ is a rational function with rational coefficients;
%\item \label{th2-2} \emph{homogeneity}\textup{:}
\item \label{th2-2} \emph{однородность}\textup{:}
%    %$C(z_1, \dots, z_n)}$ is degree ${1}$ homogeneous, i.e.,
    ${C(tz_1, \dots, tz_n)= tC(z_1,\dots, z_n)}$\textup{;}
%\item \label{th2-3} \emph{positive reality}\textup{:} если ${\Real z_1,\dots,\Real z_n>0}$ тогда ${\Real C(z_1, \dots, z_n)>0}$.
\item \label{th2-3} \emph{положительность вещественной части}\textup{:} если ${{\Real z_1,\dots,\Real z_n>0}}$, то ${{\Real C(z_1, \dots, z_n)>0}}$.
\end{enumerate}
\end{theorem}

%\begin{problem} Is the reciprocal theorem true for ${n\ge 3}$?
\begin{problem} Верно ли обратное утверждение для ${n\ge 3}$?
\end{problem}

%%We prove this theorem in \S\ref{proofs}.
%%Parts (\ref{th2-1})и~(\ref{th2-2}) of Theorem~\ref{th2} were actually proved by M.~Dehn, see also~\cite[Lemma~4]{FR}.
%Кse ${n=1}$ (respectively, ${n=2}$) of both Theorem~\ref{th2}и its reciprocal %are known: these cases
%is equivalent to Theorem~\ref{dehn} (respectively, to \cite[Theorem~5]{FLR} or else to Corollary~\ref{series-parallel} below).
%%There is another condition, which the function $c$ must satisfy, but we do не know whether this condition follows from (1)--(3) for ${n\ge 3}$. \textbf{Dopisat'.}
%%There is a heuristic argument showing that
Случай ${n=1}$ (соответственно ${n=2}$) теоремы~\ref{th2} и обратной к ней эквивалентен теореме~\ref{dehn} (соответственно теореме~5 в \cite{FLR} и следствию~\ref{series-parallel} ниже). Для ${n\ge 3}$ обратная теорема не может быть доказана нашим методом --- см. пример~\ref{heuristic}.

Теорема~\ref{th2} имеет простой физический смысл --- см.~\S\ref{alternating}. Но эта теорема
(а также обратная к ней) не является \emph{алгоритмической}, то есть, если требуемое разрезание существует, теорема не даёт алгоритма построения этого разрезания. Поэтому интересно получить менее общий, но зато алгоритмический результат.

%The first
Результат такого вида был независимо получен Фрайлингом и Ринном в 1994\7г. и Лацковичем и Секерешем в 1995\7г.

\begin{theorem} \label{lfrs} \textup{\cite{FR, LS}}
Для ${c>0}$ следующие ${3}$ условия эквивалентны\textup{:}
\begin{enumerate}[(1)]
\item \label{lfrs-1} квадрат можно разрезать на прямоугольники с отношениями сторон $c$ и ${1/c}$\textup{;}
\item \label{lfrs-2} число $c$ является корнем %ненулевого 
многочлена с целыми коэффициентами, у которого все комплексные корни имеют положительную действительную часть\textup{;}
\item \label{lfrs-3} для некоторых положительных рациональных чисел ${d_1}$, \dots, ${d_m}$ выполнено %равенство
$$d_1 c +\cfrac1{d_2 c+\cfrac1{d_3 c+\dots+\cfrac1 {d_m c}}}=1.$$
\end{enumerate}
\end{theorem}

В частности, квадрат можно разрезать на прямоугольники с отношениями сторон $2+\sqrt{2}$ и $\frac{1}{2+\sqrt{2}}$, но нельзя разрезать на прямоугольники с отношениями сторон $1+\sqrt{2}$ and $\frac{1}{1+\sqrt{2}}$; см. элементарное доказательство в~\cite{SMD}.

Мы представим новое короткое самодостаточное доказательство  теоремы~\ref{lfrs}. Это новое доказательство %(announced in \cite{PrSk})
--- естественное приложение цепей переменного тока. Мы также получим новый алгоритмический результат: % which does не seem to be tractable by methods available before:
%of \cite{FLR, FR, LS}:

\begin{theorem}\label{th1} %(cf.~\cite{FR,LS})
Для числа ${c>0}$ следующие ${3}$ условия эквивалентны:
\begin{enumerate}[(1)]
\item \label{th1-1} прямоугольник с отношением сторон $c$ можно разрезать на прямоугольники с отношениями сторон $c$ и ${1/c}$
    так, что в разрезании будет хотя бы один прямоугольник с отношением сторон ${1/c}$\textup{;}
\item \label{th1-2} число ${c^2}$ является корнем %ненулевого 
многочлена с целыми коэффициентами, у которого все остальные корни являются отрицательными действительными числами\textup{;}
\item \label{th1-3} для некоторых положительных рациональных чисел ${d_1}$, \dots, ${d_m}$ выполнено %равенство
$$\cfrac{1}{d_1 c +\cfrac1{d_2 c+\dots+\cfrac1 {d_m c}}}=c.$$
\end{enumerate}
\end{theorem}

%This is proved in \S\ref{proofs}.
Другие алгоритмические результаты можно найти в~\cite[с.~224]{FLR} и \cite{Sharov}. Аналогичные результаты о разрезаниях на треугольники см. в\7\cite{Szegedy}. Многомерные обобщения обсуждаются в\7\cite{P}. %\cite{BB,P}.
%The latter shows in particular that the case when не all ratios ${c_n}$ are ''used'' at least once in the tiling should be considered separately, thus giving a motivation to the statement of Theorem~\ref{th1}.

%\medskip
%\subsection{Tiling of polygonsи electrical circuits with many terminals}
%In \S\ref{var}
Мы также рассмотрим разрезания произвольных (необязательно выпуклых) многоугольников на прямоугольники. Это обобщение
%might seem a routine but in fact it is substantialи
раскрывает новые связи между разрезаниями и электрическими цепями.

Мы применим цепи постоянного тока с несколькими входами, чтобы получить критерий возможности разрезания \2общего\3 многоугольника на квадраты --- см. теорему~\ref{squaring-polygons} ниже. Она не является алгоритмической. Этот результат обобщает теорему~\ref{dehn} и~\cite[теоремы 9 и 12]{K}. Более простая родственная задача о разрезаниях на квадраты \emph{\2с учётом знаков\3} решена в \cite{FHTW,KK}.  %, see figure~\ref{chip}.

%\begin{figure}[htbp]
%\includegraphics[width=3.5cm]{chip.pdf}
%\caption{An electrical circuit with many terminals}
%\label{chip}
%\end{figure}

%Our approach to tilings of polygons uses direct-current circuits with several terminals.
%, studied by Colin de Verdi\`ere, Lawler, Sylvesterи others
%\emph{Direct}и \emph{inverse problems} for such circuits is the subject of \emph{electrical impedance tomography} \cite{C, CGV, LSy}.

Цепи переменного тока с несколькими входами мы применим для короткого доказательства обобщения теоремы~\ref{lfrs} на многоугольники с рациональными координатами вершин \cite{SuDi} --- см. теорему~\ref{polygon-square} ниже.
%The proofs are based on certain known results of \emph{electrical impedance tomography}.
Также мы докажем один фольклорный результат, касающийся \emph{электроимпедансной томографии}, для цепей переменного тока, ср.~\cite{CM, C, CGV, LSy}.

Существует тесная связь между электрическими цепями, дискретными гармоническими функциями и
%theory of discrete harmonic functionsи
случайными блужданиями на графах \cite{DS, Lovasz, BeSch}. %, see a brilliant introduction in \cite{DS}и a survey in \cite{Lovasz}, see also \cite{BeSch}.
Наши результаты имеют эквивалентные утверждения на языке каждой из этих теорий --- см., например, следствие~\ref{random} ниже.
%either the language of electrical circuits or discrete harmonic functions or random walks.
%The reader can find some applications of our results to these subjects in the paper.
%We prefer the language of electrical circuits only because of physical clearness.

%\subsection{Organization of the paper} In \S\ref{prelim} we recall the basics of electrical circuitsи present main ides of our proofs. In \S\ref{proofs} we prove Theorems~\ref{th2}, ~\ref{th1}и~\ref{lfrs}. In \S\ref{vars} we discuss some applications of our approach to electrical circuits with many terminalsи to random walks.

%The paper is organized as follows. It splits naturally into two formally independent parts.

Далее наша статья разделена на две формально независимых части: \S\S\ref{prelim}--\ref{proofs} и \S\S\ref{var}--\ref{finalproofs}.

Первая часть содержит доказательства теорем~\ref{th2}, \ref{lfrs} и \ref{th1}. В \S\ref{prelim} напоминаются базовые определения и свойства электрических цепей, а также объясняется их связь с разрезаниями. В \S\ref{proofs} эти результаты доказываются.

%The second part concerns some variations. In \S\ref{var} the results on
Во второй части статьи приводятся вариации на эти темы. В \S\ref{var} приводятся формулировки результатов, касающихся разрезаний многоугольников, электроимпедансной томографии и случайных блужданий.
В \S\ref{final} результаты \S\ref{prelim} обобщаются на электрические цепи с несколькими входами. В \S\ref{finalproofs} доказываются результаты, сформулированные в~\S\ref{var}.

\section{Главные идеи}\label{prelim}

\subsection{Электрические цепи}

Наш подход основан на теории электрических цепей \cite{PrSk}, но мы не предполагаем, что читатель знаком с физикой --- в этом параграфе мы напомним все нужные нам физические понятия. Несмотря на применение физических соображений, наше изложение производится на математическом уровне строгости и физический смысл происходящего объяснён очень кратко. В этом параграфе новых результатов не содержится.
Короткие доказательства даны в~\S\ref{final}.
%All the results of \S\ref{prelim} are proved in \S\ref{final} in a more general situation.

%Informally, \textit{an electrical network} is a collection of conductors connected with each other. A fixed pair of points of the network can be joined with the poles of the battery.

%Formally,
Итак, \textit{электрической цепью} называется связный граф,
%(possibly having multiple edges)
каждому ребру которого сопоставлено некоторое положительное действительное число (\textit{проводимость}) и две вершины которого отмечены (\textit{граничные вершины}). Для простоты в этом пункте мы будем считать, что граф не содержит ни кратных рёбер, ни петель. Обобщения для графов с кратными рёбрами оставляем читателю.

% и two real numbers ${U_1}$ и ${U_2}$ assigned to each boundary
% vertex, called \textit{boundary voltage}.

Электрическую цепь будем называть \textit{плоской}, если её можно изобразить на единичном круге так, что граничные вершины цепи принадлежат границе круга и никакие два ребра не имеют общих точек, кроме, возможно, их общей вершины.

%Юnote by ${V}$ и ${E}$ the sets of vertices и нges, respectively, of the graph in question.
%Fix an orientation of each edge of the graph.

Зафиксируем нумерацию вершин цепи числами ${1}$, ${2}$, ${\dots}$, ${n}$ так, чтобы вершины ${1}$ и ${2}$ были бы граничными. Обозначим количество граничных вершин в цепи через ${b:=2}$, количество рёбер в цепи --- через ${m}$, проводимость ребра между вершинами ${k}$ и ${l}$ --- через ${c_{kl}}$. Положим ${c_{kl}:=0}$, если между вершинами ${k}$ и ${l}$ нет ребра.
%In particular, all ${c_{kk}=0}$.
%It is convenient to put also ${c_{kk}=-\sum_{l=1}^n c_{kl}}$.
%We often use shorthand \textit{conductor ${c$} for the phrase \textit{an edge of conductance ${c$}.

%Informally, \textit{an electrical circuit} is an electrical network whose boundary vertices are joint to a battery of certain specific voltage.

%Formally, by
\textit{Замкнутой электрической цепью} называется электрическая цепь вместе с двумя действительными числами ${U_1}$ и ${U_2}$ --- \textit{входными потенциалами}, сопоставленными граничным вершинам цепи.

Каждой вершине ${k}$ (${1\le k\le n}$) замкнутой электрической цепи сопоставим число ${U_k}$ --- \textit{потенциал} на вершине, а каждому ребру ${kl}$ (${1\le k,l\le n}$) --- число ${I_{kl}}$, называемое \textit{током} через ребро.
Эти числа определяются следующими двумя аксиомами:
%(0) \textit{Boundary condition.}
%$U(v_1)=U_1}$, ${U(v_n)=U_n}$.
\begin{enumerate}[(1)]
\item[(C)]\textit{Закон Ома.} Для каждой пары вершин ${k,l}$ выполнено равенство ${I_{kl}=c_{kl}(U_k-U_l)}$.
\item[(I)]\textit{Правило Кирхгофа.} Для каждой вершины ${k>b}$ выполнено равенство ${{\sum_{l=1}^{n} I_{kl}=0}}$.
\end{enumerate}
Неформально, правило (I) означает, что электрический заряд не накапливается в неграничных вершинах. Другими словами, эти аксиомы утверждают, что ${U_k}$ является \emph{дискретной гармонической функцией}.
Числа ${U_k}$ и ${I_{kl}}$ однозначно определяются этими аксиомами, что доказывается следующей классической теоремой.
% proved in the end of \S\ref{prelim}.

\begin{theorem} \label{weyl} \textup{\cite{Weyl}} Для любой замкнутой электрической цепи
система линейных уравнений \textup{(C),\7(I)} от переменных ${U_k}$, где ${b< k\le n}$, и ${I_{kl}}$, где ${1\le k,l\le n}$, имеет единственное решение.
\end{theorem}

%We prove a generalization of this theorem in \S\ref{proofs}.

\begin{figure}[htbp]
\begin{tabular}{cc}
\includegraphics[width=5cm]{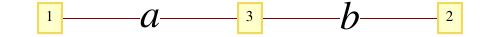} & 
\includegraphics[width=5cm]{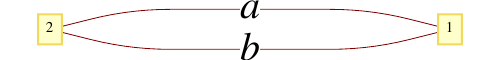}\\[9pt]
${C(a,b)=\frac{ab}{a+b}}$ & ${C(a,b)=a+b}$
\end{tabular}
\caption{последовательная и параллельная электрические цепи.}
\label{examples}
\end{figure}

Обозначим через ${I_1:=\sum_{k=1}^{n} I_{1k}}$ ток, текущий внутрь замкнутой электрической цепи через вершину ${1}$.
\textit{Проводимостью} замкнутой электрической цепи с входными потенциалами ${U_1\ne U_2}$ называется число ${C=I_1/(U_1-U_2)}$. Легко видеть, что проводимость не зависит от ${U_1}$ и ${U_2}$.
Поэтому проводимость электрической цепи определена однозначно. Величина, обратная к проводимости, называется \emph{сопротивлением}.
%Thus one can define \textit{the conductance} of an electrical network to be the conductance of any circuit made of the network.
Базовые примеры электрических цепей и вычисления их проводимостей изображены на рис.~\ref{examples}.

\subsection{Электрические цепи и разрезания}\label{ssec2-2}

%Let us explain informally the relationship between tilings by rectangles и electrical networks. Suppose that a rectangle is tiled by rectangles. \textbf{Plastinka. This has не been written yet.}

Как отмечалось, существует тесная связь между электрическими цепями и разрезаниями.
Мы будем говорить, что ребро ${kl}$ замкнутой цепи \emph{существенно}, если %for some (and thus for all) ${U_1\ne U_2}$ we have
${I_{kl}\ne 0}$. Ясно, что эта характеристика ребра не зависит от ${U_1}$ и ${U_2}$, если ${U_1\ne U_2}$. Поэтому существенные рёбра в (незамкнутой) электрической цепи определены корректно.

\begin{lemma} \label{correspondence-lemma} \textup{\cite[Теорема 1.4.1]{BSST,CFP}} %(see figure~\ref{tilingnetworks})
Следующие два условия эквивалентны\textup{:}
\begin{enumerate}[(1)]
\item \label{c-l-1} прямоугольник с отношением сторон $c$ можно разрезать на ${m}$ прямоугольников с отношениями сторон ${c_1}$, ${\dots}$, ${c_m}$\textup{;}
\item \label{c-l-2} существует плоская электрическая цепь с проводимостью $c$, состоящая из ${m}$ существенных рёбер с проводимостями
$c_1$, ${\dots}$, ${c_m}$.
\end{enumerate}
\end{lemma}

%Only assertion (1)$\implies$(2) is required for the proof of main results.
Приведём здесь набросок доказательства утверждения (\ref{c-l-1})$\implies$(\ref{c-l-2}).
%, the only part of the lemma required for the proof of main results
По разрезанию, описанному в~(\ref{c-l-1}), можно построить электрическую цепь, описанную в (\ref{c-l-2}), --- см. рис.~\ref{tilingnetworks}. На каждой горизонтальной линии разреза, а также на обоих горизонтальных сторонах разрезаемого прямоугольника отметим по точке --- эти точки будут вершинами электрической цепи. Для каждого прямоугольника в разрезании проведём ребро между вершинами на разрезах, которые содержат горизонтальные стороны рассматриваемого прямоугольника. Проводимость этого ребра положим равной отношению сторон прямоугольника.
%Put the current through the edge equal to the horizontal side of the rectangle.
%Put the voltage at a vertex equal to the $y$-coordinate of the vertex.
%тогда laws (C),(I) are satisfied и
Тогда проводимость всей цепи равна отношению сторон разрезаемого прямоугольника --- более подробно см. в\7\S\ref{correspond}.
%This network has conductance $c$.

%We prove Lemma~\ref{correspondence-lemma} in the end of this section.

%\textbf{Proof of lemma. не written yet!!!}

%Formally, \textit{a circuit corresponding to a tiling} is defined as follows. \textbf{Correspondence. %This has не been written yet.}

\begin{figure}[htbp]
\definecolor{qqwwff}{rgb}{0,0.4,1}
\definecolor{ffqqqq}{rgb}{1,0,0}
\definecolor{zzttqq}{rgb}{0.6,0.2,0}
\begin{tikzpicture}[line cap=round,line join=round,>=triangle 45,x=0.5cm,y=0.5cm]
\clip(-3.14,-4.24) rectangle (7.08,6.78);
\fill[color=zzttqq,fill=zzttqq,fill opacity=0.1] (-3,6) -- (6.9,6) -- (6.9,-3.6) -- (-3,-3.6) -- cycle;
\draw [color=zzttqq] (-3,6)-- (6.9,6);
\draw [color=zzttqq] (6.9,6)-- (6.9,-3.6);
\draw [color=zzttqq] (6.9,-3.6)-- (-3,-3.6);
\draw [color=zzttqq] (-3,-3.6)-- (-3,6);
\draw [color=zzttqq] (2.4,6)-- (2.4,0.6);
\draw [color=zzttqq] (2.4,0.6)-- (-3,0.6);
\draw [color=zzttqq] (6.9,1.5)-- (2.4,1.5);
\draw [color=zzttqq] (1.2,-3.6)-- (1.2,0.6);
\draw [color=zzttqq] (1.2,-0.6)-- (4.5,-0.6);
\draw (2.4,0.6)-- (2.4,-0.6);
\draw [color=zzttqq] (4.5,1.5)-- (4.5,-0.9);
\draw [color=zzttqq] (6.9,-0.9)-- (4.2,-0.9);
\draw [color=zzttqq] (4.2,-0.6)-- (4.2,-3.6);
\draw [line width=1pt,color=qqwwff] (1.98,6)-- (-0.3,0.6);
\draw [line width=1pt,color=qqwwff] (1.98,6)-- (4.66,1.5);
\draw [line width=1pt,color=qqwwff] (-0.3,0.6)-- (1.94,-3.6);
\draw [line width=1pt,color=qqwwff] (2.7,-0.6)-- (4.66,1.5);
\draw [line width=1pt,color=qqwwff] (-0.3,0.6)-- (2.7,-0.6);
\draw [line width=1pt,color=qqwwff] (5.56,-0.9)-- (4.66,1.5);
\draw [line width=1pt,color=qqwwff] (5.56,-0.9)-- (1.94,-3.6);
\draw [line width=1pt,color=qqwwff] (2.7,-0.6)-- (1.94,-3.6);
\draw [line width=1pt,color=qqwwff] (5.56,-0.9)-- (2.7,-0.6);
\fill [color=ffqqqq] (1.98,6) circle (2.0pt);
\draw[color=ffqqqq] (1.98,6.5) node {$+$};
\fill [color=ffqqqq] (-0.3,0.6) circle (2.0pt);
\fill [color=ffqqqq] (4.66,1.5) circle (2.0pt);
\fill [color=ffqqqq] (1.94,-3.6) circle (2.0pt);
\draw[color=ffqqqq] (1.96,-3.95) node {$-$};
\fill [color=ffqqqq] (2.7,-0.6) circle (2.0pt);
\fill [color=ffqqqq] (5.56,-0.9) circle (2.0pt);
\end{tikzpicture} %}$\!\!\!\!}$\input{metal-fig6c.tex}
\caption{соответствие между разрезаниями и электрическими цепями.}
\label{tilingnetworks}
\end{figure}

\subsection{Свойства проводимости}

Коротко изложим основные свойства проводимости.

\begin{lemma} \label{electrical} \label{main}
\textup{\cite{BSST,DS,Foster}}
Пусть электрическая цепь содержит ${m}$ рёбер с проводимостями ${c_{1}, \dots, c_{m}}$.
Тогда проводимость ${C(c_{1}, \dots, c_{m})}$ всей цепи имеет следующие свойства\textup{:}
\begin{enumerate}[(1)]
\item \label{el-1} \emph{рациональность}\textup{:} ${C(c_{1}, \dots, c_{m})\in \mathbb{Q}(c_{1}, \dots, c_{m})}$\textup{;}
\item \label{el-2} \emph{однородность}\textup{:} ${C(tc_{1}, \dots, tc_{m})=tC(c_{1}, \dots, c_{m})}$\textup{;}
\item \label{el-3} ${\frac{\partial}{\partial c_j} C(c_{1}, \dots, c_{m})=(\frac{U_k-U_l}{U_1-U_2})^2}$,
    где ${k}$ и ${l}$ --- концы ребра ${j}$\textup{;}
\item \label{el-4} \emph{монотонность}\textup{:} если ${c_{1}, \dots, c_{m}>0}$, то ${\frac{\partial}{\partial c_j} C(c_{1}, \dots, c_{m} )\ge 0}$; если ребро ${j}$ существенно, то последнее неравенство строгое\textup{;}
\item \label{el-5} \emph{положительность вещественной части}\textup{:} если ${\Real c_{1},\dots, \Real c_{m}>0}$, то ${\Real C(c_{1}, \dots, c_{m})>0}$.
\end{enumerate}
\end{lemma}

%\begin{problem} Is the converse theorem true?
%\end{problem}

%\begin{remark} (Akopyan, private communication) ${C(c_{1}, \dots, c_{m})}$ belongs to a convex cone with vertex at the origin, spanned by the points ${c_{1}, \dots, c_{m}}$ in ${\mathbb{C}}$.
%\end{remark}

\begin{remark} (А. Акопян, частное сообщение) Свойство (\ref{el-4}) следует из (\ref{el-1}), (\ref{el-2}) и~(\ref{el-5}).
Свойство~(\ref{el-5}) не следует из~(\ref{el-1}), (\ref{el-2}) и~(\ref{el-4}), --- например, функция ${{C(c_1,c_2)=(c_1+c_2)\frac{c_1^2+c_2^2}{c_1^2+2c_2^2}}}$ удовлетворяет~(\ref{el-1}), (\ref{el-2}), (\ref{el-4}), но свойство~(\ref{el-5}) не выполняется.
\end{remark}

%Theorem \ref{dehn} follows directly from Lemma \ref{electrical}(1).
%Properties (\ref{el-1})--(\ref{el-4}) are well-known,

Свойство~(\ref{el-5}) касается продолжения функции ${C(c_1,\dots,c_m)}$ в комплексную область.
%is more tricky because it involves an extension of ${C(c_{1}, c_{2}, \dots, c_{m})}$ to the complex plane. The latter result
%It does не seem that anyone has paid attention to this fundamental property for \textit{direct}-current circuits. Certainly it is well-known for
Смысл этого свойства станет ясным в контексте цепей \textit{переменного} тока. Короткое доказательство леммы см. в \S\ref{elcirc}.

\subsection{Переменный ток}\label{alternating}

Объясним физический смысл
условия положительности вещественной части в
лемме~\ref{main} и теореме~\ref{th2}. Это объяснение далее в статье не используется и читатель может пропустить этот пункт. % и proceed to the formal proofs in the next one.

%It is quite natural to study a rational function ${C(c_{ij})}$ on the whole complex plane. A more surprising observation is that the function ${C(c_{kl})}$ for complex ${c_{kl}}$ has a clear physical meaning. Let us introduce the basics of \textit{alternating-current circuits}, which are logically equivalent to electrical circuits with complex-valued conductances.

%Electrical circuits with complex-valued conductances are ''equivalent'' toalternating-current circuits.

Неформально, \textit{замкнутая цепь переменного тока} --- это набор соединённых между собой резисторов, конденсаторов, катушек индуктивности и единственного источника переменного тока.

%\textbf{Give the following algorithm applied to a specific example only?}

Формально, {\it замкнутая цепь переменного тока} --- это %connected
граф со следующей ­ditional
структурой:
\begin{itemize}
\item две отмеченные (\emph{граничные}) вершины;
\item две функции (\emph{функции потенциала}) ${\tilde U_1(t):=U\cos \omega t}$ и ${\tilde U_2(t):=0}$, соответствующие им; %;
\item разделение рёбер на три типа (\textit{резистор, конденсатор и катушка индуктивности});
\item положительное число ${\tilde c_{kl}}$, соответствующее каждому ребру (называется \textit{проводимостью, ёмкостью} или \textit{индуктивностью}, в зависимости от типа ребра).
\end{itemize}
Предположим для простоты, что такой граф не имеет кратных рёбер. \emph{Потенциалы} ${\tilde U_k(t)}$ и \emph{токи} ${\tilde I_{kl}(t)}$ определяются следующими аксиомами:
\begin{enumerate}[(1)]
\item[(\~ C)]\textit{Закон Ома для переменного тока.} Для каждого ребра ${kl}$ выполнено равенство
$$\tilde I_{kl}(t)\begin{cases}\tilde c_{kl}\left(\tilde U_k(t)-\tilde U_l(t)\right),& \text{если ${kl}$ --- резистор;}\\
\tilde c_{kl}\frac{d}{dt}\left(\tilde U_k(t)-\tilde U_l(t)\right),& \text{если ${kl}$ --- конденсатор;}\\
\tilde c_{kl}\int_{\pi/2\omega}^t\left(\tilde U_k(t)-\tilde U_l(t)\right)dt,& \text{если ${kl}$ --- катушка индуктивности.}
\end{cases}
$$
\item[(\~ I)]\textit{Правило Кирхгофа.} Для каждой вершины ${k\ne 1,2}$ выполнено равенство ${\sum_{l=1}^{n} \tilde I_{kl}(t)=0}$.
\end{enumerate}

Потенциалы и токи можно найти, используя следующий широко известный алгоритм. Положим ${U_1:=U}$, ${U_2:=0}$, ${i:=\sqrt{-1}}$ и
$$c_{kl}:\begin{cases}\quad\tilde c_{kl},  & \text{если ${kl}$ --- резистор;}\\
i\omega\tilde c_{kl},             & \text{если ${kl}$ --- конденсатор;}\\
\frac{1}{i\omega}\tilde c_{kl},   & \text{если ${kl}$ --- катушка индуктивности.}
\end{cases}
$$
Буквально по таким же законам (C), (I), как для цепей постоянного тока, найдём числа ${U_k}$, где ${3\le k\le n}$, и ${I_{kl}}$, где ${1\le k,l\le n}$.
%тогда the real parts of these numbers are the voltages и the currents in the circuit at the moment ${t}$. (Notice that the arguments of the complex numbers ${U_k}$ и ${I_{kl}}$ depend on the time parameter ${t}$.)
Тогда ${\tilde U_k(t)=\Real (U_k e^{i\omega t})}$, ${\tilde I_{kl}(t)=\Real (I_{kl}e^{i\omega t})}$.
В этом смысле цепи переменного тока \2эквивалентны\3 цепям постоянного тока с комплексными проводимостями (также называемыми \emph{адмитансами}).

%\textbf{Example. не written yet!!!}

%Let us explain physical meaning of Lemma~\ref{electrical}(5) и Theorem~\ref{th2}(3).
Заметим, что всегда ${\Real c_{kl}\ge 0}$. Физически это означает неотрицательность
\textit{диссипации энергии} на ребре ${kl}$ (которая равна ${\Real  c_{kl}|U_k-U_l|^2}$ за период).
%\begin{equation*}
%\int_{0}^{2\pi/\omega}\tilde I_{kl}(t)(\tilde U_k(t)-\tilde U_l(t))dt=\Real  (\bar I_{kl}(U_k-U_l))=\Real \bar c_{kl}|U_k-U_l|^2.
%\end{equation*}
%Thus a physical meaning of positive real part of \emph{admittance} ${c_{kl}}$ is positive energy dissipation.
Поэтому физический смысл положительности вещественной части такой: \2цепь, состоящая из элементов, диссипирующих энергию, тоже диссипирует энергию\3.
%We make this intuitive argument precise in \S\ref{final}.

%Thus we come to the following formal definition.

Каждому разрезанию прямоугольника  %(or even a polygon) 
на \emph{подобные} прямоугольники можно сопоставить цепь переменного тока, --- см. рис.~\ref{fig-perem} (формально это в статье не используется). 
Граф цепи такой же, как в пункте~\ref{ssec2-2}. Но вместо резисторов, рёбра теперь являются либо конденсаторами, либо катушками индуктивности, в зависимости от ориентации рассматриваемых прямоугольников. Построенная цепь переменного тока \2эквивалентна\3 цепи постоянного тока с комплексными проводимостями рёбер ${z}$ и ${1/z}$, где ${z:=i\omega}$.

\begin{figure}[h]
\hspace{4cm}
\makebox[1.5cm][r]{
  \begin{tabular}{l}
  \definecolor{zzttqq}{rgb}{0.6,0.2,0}
\begin{tikzpicture}[line cap=round,line join=round,>=triangle 45,x=0.3cm,y=0.3cm]
\clip(-3.42,-2.1) rectangle (4.32,5.92);
\fill[color=zzttqq,fill=zzttqq,fill opacity=0.1] (-2,5) -- (4,5) -- (4,-1) -- (-2,-1) -- cycle;
\draw [color=zzttqq] (-2,5)-- (4,5);
\draw [color=zzttqq] (4,5)-- (4,-1);
\draw [color=zzttqq] (4,-1)-- (-2,-1);
\draw [color=zzttqq] (-2,-1)-- (-2,5);
\draw [color=zzttqq] (-2,3.58)-- (4,3.56);
\draw [color=zzttqq] (-0.02,3.57)-- (0,-1);
\draw [color=zzttqq] (2,3.57)-- (2,-1);
\draw [color=zzttqq] (1,5)-- (1,3.57);
\end{tikzpicture}
  \end{tabular}
  }
  \makebox[3.5cm][c]{
  \begin{tabular}{c}
  %\documentclass{article}
%\usepackage{circuitikz}
%\begin{document}
\ctikzset{bipoles/length=.6cm}
\begin{circuitikz}[scale=0.6]
\draw[color=brown]
(0,0) to[short, *-*] (1,0)
      to[short, *-*] (2,0)
      to[short, *-] (4,0)
      to[sV](4,3)
      to[short, -*] (2,3)
      to[short, *-*] (0,3)
      to[C, *-*]     (0,2)
(0,0) to[L, *-*]     (0,2)
(1,0) to[L, *-*]     (1,2)
      to[short, *-*] (0,2)
(1,2) to[short, *-*] (2,2)        
(2,0) to[L, *-*]     (2,2)
      to[C, *-*]     (2,3)
;\end{circuitikz}
%\end{document}
  \end{tabular}
  }
\caption{соответствие между разрезаниями на прямоугольники и цепями переменного тока.}
\label{fig-perem}
\end{figure}

\subsection{Обратные задачи}\label{inverse-problems}

%The proofs are based on the above physical interpretation.
%Let us prepare for the proof of other main results.
%This subsection is required only

\emph{Обратной задачей} для электрических цепей называется задача построения цепи с данной проводимостью из данных элементов.
%Let us state a simple inverse problem и outline its solution.
Приведём классический результат Р.\7М.\7Фостера и У.\7Кауэра в этой области. Его короткое доказательство можно найти в~\S\ref{basicposreal}.
Всё сказанное в данном пункте в дальнейшем используется только в доказательстве утверждений~(2)$\implies$(3) в теоремах~\ref{lfrs}--\ref{th1}.

\begin{theorem}[теорема Фостера о реактивном сопротивлении] \label{th-foster} \textup{\cite{Foster}} Следующие свойства рациональной функции ${C(z)\in\mathbb{R}(z)}$ эквивалентны\textup{:}
\begin{enumerate}[(1)]
\item \label{th-foster-1}
    ${C(i\omega)}$ является проводимостью электрической цепи, состоящей из конденсаторов и катушек индуктивности, где проводимость рассматривается как функция от частоты ${\omega}$\textup{;}
\item \label{th-foster-2}
    ${C(z)}$ является проводимостью электрической цепи, в которой проводимость каждого ребра ${j}$ равна либо ${d_jz}$, либо ${1/d_jz}$ для некоторых действительных чисел ${d_1,\dots,d_m>0}$\textup{;}
\item \label{th-foster-3} ${C(z)}$ является нечётной рациональной функцией, для которой ${\Real C(z)>0}$, если ${\Real z>0}$.
\end{enumerate}
\end{theorem}

%\begin{problem} Construct a network from capacitors и inductors such that the conductance of the network is a given function of the frequency ${\omega}$.
%Equivalently, construct an electrical network such that the conductance of each edge ${j}$ is either ${d_jz}$ or ${1/d_jz}$ for some real numbers ${d_j>0}$ и the conductance of the network is a given function ${C(z)}$.
%\end{problem}

%Let us outline the proof of assertion~(\ref{th-foster-3})$\implies$(\ref{th-foster-2})  given in \cite{Cauer}.
%Let us outline the solution \cite{Foster, Cauer}. By Lemma~\ref{electrical} it follows that ${C(z)}$ must be an odd rational function, which is \emph{positive real}, i.~e., ${\Real C(z)>0}$, если ${\Real z>0}$.
%Using this property we construct a series-parallel electrical network with the conductance ${C(z)}$.

%The proof of assertions~(2)$\implies$(3) in Theorems~\ref{lfrs} и~\ref{th1} is based on the following analytical result.

Доказательство теоремы основано на следующей лемме из анализа.

\begin{lemma} \label{PosReal} \textup{\cite{Cauer}} Для нечётной рациональной функции ${C(z)\in \mathbb{R}(z)}$, удовлетворяющей условию ${\lim_{z\to \infty}C(z)\ne0}$,
следующие ${5}$ условий эквивалентны\textup{:}
\begin{enumerate}[(1)]
\item\label{preserve} если ${\Real z>0}$, то ${\Real C(z)>0}$\textup{;}
%\item\label{preserve} ${\Real C(z)>0}$ если и only если ${\Real z>0}$;
\item\label{roots} если ${C(z)=1}$, то ${\Real z>0}$\textup{;}
%for some ${x>0}$ all roots of the equation ${C(z)=x}$ have positive real part;
\item\label{derivatives} если ${C(z)=0}$, то ${\Real z=0}$ и ${C'(z)>0}$\textup{;}
    %\textup{(}here ${z\in\mathbb{C}}$ or ${z=\infty}$\textup{)};
\item\label{factorization} %either ${C(z)}$ or ${1/C(z)}$ equals
$$
C(z)=d_1z\prod\limits_{k=1}^n\frac{z^2+a_k^2}{z^2+b_k^2},
$$
для некоторого целого числа ${{n\ge 0}}$ и действительных чисел ${{d_1>0}}$, ${{a_1>b_1>a_2>\dots>b_n\ge 0}}$\textup{;}
%$$\frac{\prod\limits_{k=1}^n
%(z^2+b_k^2)}{d_1z\prod\limits_{k=1}^m(z^2+a_k^2)},
%$$
%for some integer numbers ${m,n\ge 0}$ и real numbers ${a_k, b_k, d_1>0}$ such that
%$m}$ equals either ${n{-}1}$ or ${n}$ и
%$0<b_1<a_1<b_2<a_2<\dots}$;
\item\label{continued} %either ${C(z)}$ or ${1/C(z)}$ equals
$$
C(z)=d_1 z +\cfrac1{d_2 z+\dots+\cfrac1 {d_m z}},
$$
%$$\cfrac1{d_1 z +\cfrac1{d_2 z+\dots+\cfrac1 {d_l z}}},$$
для некоторого целого числа ${m\ge1}$ и действительных чисел ${d_1,\dots,d_m>0}$\textup{;}
%\item\label{series-parallel-network}
%$C(z)}$ is the conductance of a series-parallel network with edges of conductances ${d_1z}$, ${1/d_2z}$, ${d_3z}$, ${1/d_4z}$, \dots, ${(d_mz)^{(-1)^m}}$ for some integer ${m\ge1}$ и real numbers ${d_1,\dots,d_m>0}$.
\end{enumerate}
\end{lemma}

Определим по индукции \emph{последовательно-параллельные} электрические цепи. Цепь, содержащая только одно ребро является последовательно-параллельной. Если ${a}$ и ${b}$ --- две последовательно-параллельных цепи, то и их последовательное и параллельное \2соединения\3 (см. рис.~\ref{examples}) последовательно-параллельны. Видно, что условие~(\ref{continued}) леммы\7\ref{PosReal} позволяет нам строить последовательно-параллельные электрические цепи с общей проводимостью ${C(z)}$ и проводимостями рёбер ${d_1z}$, ${1/d_2z}$, ${d_3z}$, ${1/d_4z}$, \dots, ${(d_mz)^{(-1)^m}}$; см.~рис.~\ref{fig-ladder}. Утверждение~(\ref{roots})$\implies$(\ref{continued}) доказано в~\cite[Lemma 4]{FLR}; доказательство использует результаты из~\cite{Wall}. Наше доказательство в~\S\ref{basicposreal}, следующее статье~\cite{Cauer}, проще.

\begin{figure}[htb]
%\documentclass{article}
%\usepackage{circuitikz}
%\begin{document}
\ctikzset{bipoles/length=1.0cm}
\begin{circuitikz}[scale=1.0]
\draw[color=brown]
(0,1) to[L=$d_1$, *-*] (2,1)
      to[L=$d_3$, *-*] (4,1)
      to[L=$d_5$, *-*] (6,1)
(8,1) to[L=$d_{m-1}$, *-] (10,1)
(10,0) to[C=$\frac{1}{d_{m}}$] (10,1)
(10,0)to[short,-*] (8,0) 
(6,0)	to[short,*-*] (4,0)
(4,0)	to[short,*-*] (2,0)
(2,0)	to[short,*-*] (0,0)		        
(0,0)to[sV=$\omega$](0,1)
(2,0)to[C=$\frac{1}{d_{2}}$, *-] (2,1)
(4,0)to[C=$\frac{1}{d_{4}}$, *-] (4,1)
(6,0)to[C=$\frac{1}{d_{6}}$, *-] (6,1)
%(8,0)to[C=$\frac{1}{d_{8}}$, *-] (8,1)
;
\draw[dashed,color=brown]
(6,0) to[short, *-*] (8,0)      
(6,1) to[short, *-*] (8,1)  
;\end{circuitikz}
%\end{document}
\caption{Последовательно-параллельная цепь, стоящая за леммой~\ref{PosReal}}
\label{fig-ladder}
\end{figure}

%A rational function satisfying condition~(\ref{preserve}) of Lemma~\ref{PosReal} is called \emph{positive real}.
%Implication (\ref{preserve})$\implies$(\ref{factorization}) is known as \emph{the Foster reactance theorem}~\cite{Foster}.
%The rest of the lemma is proved in \cite{Cauer}.

%Short proofs of Lemma~\ref{PosReal} и the following corollary are given in~\S\ref{basicposreal}.
%%A physical meaning of Lemma~\ref{PosReal} is synthesis of a circuit consisting of capacitors и inductors such that the conductance of the circuit is a given odd positive real function ${C(z)}$ in the frequency ${z=i\omega}$.

%Consider electrical circuits, in which all the edges have conductances ${z}$ и ${1/z}$, ${\Real z>0}$. (They have a natural physical meaning: circuits consisting of capacitors и inductors with incoming voltage of complex frequency ${z/i}$.) Let us describe possible conductances ${C(z)}$ of such electrical circuits. By Lemma~\ref{electrical}(\ref{el-1}), (\ref{el-2}) и~(\ref{el-5}) the functions ${C(z)}$ are \emph{positive real}, i.e., satisfy condition~(\ref{preserve}) of the following lemma.
%Юnote by ${\Real \infty=0}$, ${C(\infty)=\lim_{z\to 0}C(1/z)}$ и ${C'(\infty)=\lim_{z\to 0}(C(1/z))'}$.
%for certain complex number ${z}$.
%The conductance of such a network is a function ${C(z)}$ in one variable ${z}$. By Lemma~\ref{electrical}(1), (2) и (5)
%this function satisfies property~(2) of the following lemma.

%\subsection{A non-series-parallel network}

%The following result shows that Theorem~\ref{th2} cannot be %proved by our method for ${n\ge 3}$.

\begin{corollary} \label{series-parallel} \textup{\cite{Cauer}}
Если функция ${C(c_1,c_2)}$ удовлетворяет условиям~\textup{(\ref{th2-1})--(\ref{th2-3})} теоремы~\ref{th2}, то
${C(c_1,c_2)}$ является проводимостью некоторой последовательно-параллельной электрической цепи с проводимостями рёбер ${c_1}$ и ${c_2}$.
\end{corollary}

\begin{example} \label{heuristic} Для трёх переменных ${c_1}$, ${c_2}$, ${c_3}$ следствие\7\ref{series-parallel} перестаёт быть верным.
Например, возьмём цепь с четырьмя вершинами и проводимостями рёбер ${c_{13}=c_1}$, ${c_{23}=c_2}$, ${c_{24}=c_1}$, ${c_{14}=c_2}$, ${c_{34}=c_3}$. По лемме~\ref{electrical}(\ref{el-3}) и соображениям симметрии получаем, что ${\partial C(c_1,c_2,c_3)/\partial c_3=0}$, если ${c_1=c_2}$. Таким образом, ${C(c_1,c_2,c_3)}$ не может быть проводимостью последовательно-параллельной цепи с проводимостями рёбер ${c_1,c_2,c_3}$, потому что все рёбра такой цепи существенны.
\end{example}

\section{Доказательства основных результатов} \label{proofs}

\subsection{Доказательство теоремы \ref{th2}}

Далее в \textit{электрических цепях} в качестве проводимостей рёбер допускаются произвольные комплексные числа с положительной действительной частью. Это обобщение мотивировано \S\ref{alternating} (и позволяет описывать как цепи постоянного, так и переменного тока).
Теорема~\ref{th2} является прямым следствием результатов~\S\ref{prelim}:

\begin{proof}[Доказательство теоремы \ref{th2}]
Предположим, что прямоугольник с отношением сторон $c$ можно разрезать на прямоугольники с отношениями сторон ${c_1}$, ${\dots}$, ${c_n}$. По лемме \ref{correspondence-lemma} существует электрическая цепь с проводимостью $c$, состоящая из рёбер с проводимостями ${c_1}$, ${\dots}$, ${c_n}$. Для каждого ${k=1,\dots,n}$ заменим каждое ребро с проводимостью ${c_k}$ в цепи на ребро с комплексной проводимостью ${z_k}$, где ${\Real z_k>0}$. Пусть ${C(z_1,\dots,z_n)}$ --- проводимость полученной электрической цепи.
По лемме \ref{electrical}(\ref{el-1}), (\ref{el-2}) и~(\ref{el-5}) функция ${C(z_1,\dots,z_n)}$
обладает свойствами~(\ref{th2-1})--(\ref{th2-3}) из теоремы~\ref{th2}.
\end{proof}

\subsection{Доказательство теоремы~\ref{lfrs}} \label{prooflfrs}

\begin{proof}[Доказательство теоремы~\ref{lfrs}]
(\ref{lfrs-3})$\implies$(\ref{lfrs-1}) \cite{FLR} Предположим, что условие~(\ref{lfrs-3}) теоремы~\ref{lfrs} выполнено и, скажем, ${m}$ нечётно.
%Let us construct a tiling of a square by rectangles of ratios $c$ и ${1/c}$.
Рассмотрим квадрат с единичной площадью. Вертикальным разрезом отрежем от этого квадрата прямоугольник с отношением сторон ${d_1c}$. Оставшаяся часть является прямоугольником с отношением сторон $$
1-d_1c=\cfrac1{d_2c +\dots+\cfrac1{d_mc}}.
$$
Теперь отрежем от оставшейся части горизонтальным разрезом прямоугольник с отношением сторон ${1/d_2c}$. Останется прямоугольник с отношением сторон
$$
d_3 c+\cfrac1{d_4 c+\dots+\cfrac1{d_mc}}.
$$
Продолжим этот процесс, чередуя вертикальные и горизонтальные разрезы. Условие~(\ref{lfrs-3}) гарантирует нам, что после ${{(m-1)\mbox{-го}}}$ шага мы получим прямоугольник с отношением сторон ${d_mc}$. Мы получили разрезание квадрата на прямоугольники с отношениями сторон ${d_1c}$, ${1/d_2c}$, ${d_3c}$, ${1/d_4c}$, ${\dots}$, ${d_mc}$.
Поскольку все ${d_k\in\mathbb{Q}}$, то все прямоугольники в полученном разрезании можно разрезать на прямоугольники с отношениями сторон $c$ и ${1/c}$.

(\ref{lfrs-1})$\implies$(\ref{lfrs-2}). Предположим, что квадрат разрезан на прямоугольники с отношениями сторон $c$ и ${1/c}$.
%Compressing the figure we get a square dissected into squares и rectangles of ratio ${R^2}$.
По лемме \ref{correspondence-lemma} существует электрическая цепь с общей проводимостью ${1}$ и проводимостями рёбер $c$ и ${1/c}$.
% Consider the corresponding electrical circuit.
Заменим каждое ребро с проводимостью $c$ (соответственно ${1/c}$) в этой цепи на ребро с проводимостью ${z\in\mathbb{C}}$ (соответственно ${1/z}$).
Пусть ${C(z)}$ --- проводимость полученной цепи. Тогда по лемме \ref{electrical}(\ref{el-1}) получаем ${C(z)\in\mathbb{Q}(z)}$, а по лемме \ref{electrical}(\ref{el-2}) функция ${C(z)}$ нечётна.
%Let ${q(z)}$ be the entire resistance of the obtained network (defined by the Kirchhoff laws).

Так как ${C(c)=1}$ и функция ${C(z)\in\mathbb{Q}(z)}$ непостоянна, то $c$ является корнем ненулевого многочлена с целыми коэффициентами. Покажем, что такой многочлен минимальной степени --- искомый.
Пусть $z$ --- любой другой корень минимального многочлена числа $c$. Тогда и ${C(z)=1}$.

Докажем, что ${\Real z>0}$.
Сначала предположим, что ${\Real z<0}$. Тогда ${\Real (-z)>0}$ и ${\Real (-1/z)>0}$.
По лемме \ref{electrical}(\ref{el-5}) получаем ${0<\Real C(-z)=-\Real C(z)=-1}$ --- противоречие. Теперь пусть ${\Real z=0}$. Пусть последовательность $z_k$ стремится к $z$, где ${\Real z_k<0}$ для всех $k$. Вновь ${0<\Real C(-z_k)=-\Real C(z_k)\rightarrow-1}$ --- противоречие. Поэтому ${\Real z>0}$.

(\ref{lfrs-2})$\implies$(\ref{lfrs-3})
%\begin{lemma}\label{l1} Consider a rational function $$q(z)=\frac{\prod\limits_{k=1}^n (z^2+b_k^2)}{z\prod\limits_{k=1}^m(z^2+a_k^2)},$$ where ${m}$ equals ${n{-}1}$ or ${n}$, ${a_k}$ и ${b_k}$ are real numbers и ${0<b_1<a_1<b_2<a_2<\dots}$. тогда ${\Real q(z)>0}$ если ${\Real z>0}$ и ${\Real q(z)<0}$ если ${\Real z<0.$
%\end{lemma}
\cite{FLR}
Пусть ${p(z)}$ --- минимальный многочлен числа $c$. Пусть ${\deg p(z)}$ нечётна. Положим ${C(z):=\frac{p(-z)-p(z)}{p(-z)+p(z)}}$.
Тогда ${C(c)=1}$, ${C(z)\in \mathbb{Q}(z)}$, ${\lim_{z\to\infty}C(z)\ne 0}$, ${C(z)}$ нечётна и все корни уравнения ${C(z)=1}$ имеют положительную действительную часть.
По лемме~\ref{PosReal}(\ref{roots})$\implies$(\ref{continued}) функция ${C(z)}$ удовлетворяет условию~(\ref{continued})
леммы~\ref{PosReal}. Так как ${C(z)\in \mathbb{Q}(z)}$, то по алгоритму Евклида все ${d_k\in \mathbb{Q}}$. Подставляя ${z=c}$, получаем требуемое условие. В случае, если ${\deg p(z)}$ --- чётное число, нужно заменить ${C(z)}$ на ${1/C(z)}$ и повторить вышесказанное.
\end{proof}

\subsection{Доказательство теоремы~\ref{th1}}

Здесь мы используем идеи из \S\ref{prooflfrs} и \S\ref{basicposreal}.

\begin{proof}[Доказательство теоремы~\ref{th1}]
(\ref{th1-3})$\implies$(\ref{th1-1}) Аналогично доказательству теоремы~\ref{lfrs}, (\ref{lfrs-3})$\implies$(\ref{lfrs-1}).

(\ref{th1-1})$\implies$(\ref{th1-2}). Предположим, что прямоугольник с отношением сторон $c$ разрезан на прямоугольники с отношениями сторон $c$ и ${1/c}$.
Поворачивая разрезание на ${\pi/2}$ и растягивая его по горизонтали в $c$ раз, мы получим квадрат, 
разрезанный на квадраты и прямоугольники с отношением сторон ${c^2}$. По лемме~\ref{correspondence-lemma} существует электрическая цепь с общей проводимостью ${1}$ и проводимостями рёбер ${1}$ и ${c^2}$, в которой все рёбра существенны. Так как в исходном разрезании существует по крайней мере один прямоугольник с отношением сторон ${1/c}$, то рассматриваемая цепь содержит по крайней мере одно ребро с проводимостью ${c^2}$.
% Consider the corresponding electrical circuit.
Заменим каждое ребро с проводимостью ${c^2}$ (соответственно ${1}$) в цепи на ребро с проводимостью ${z\in\mathbb{C}}$ (соответственно ${w\in\mathbb{C}}$).
Пусть ${C(z,w)}$ --- проводимость полученной цепи, а ${C(z):=C(z,1)}$.
%Let ${q(z)}$ be the entire resistance of the obtained network (defined by the Kirchhoff laws).

Докажем, что ${c^2}$ является корнем ненулевого многочлена с целыми коэффициентами. Действительно, по лемме~\ref{electrical}(\ref{el-4}) имеем ${C'(c^2)>0}$,
потому что в цепи существует по крайней мере одно существенное ребро с проводимостью ${c^2}$. Поэтому функция ${C(z)\not\equiv\mbox{const}}$.
По лемме~\ref{electrical}(\ref{el-1}) получаем, что ${C(z)\in\mathbb{Q}(z)}$. Так как ${C(c^2)=1}$, то ${c^2}$ является корнем ненулевого многочлена с целыми коэффициентами. Возьмём такой многочлен минимальной степени. Пусть $z$ --- любой его корень, отличный от $c^2$. Тогда ${C(z,1)=C(c^2)=1}$.

Докажем, что ${z}$ является отрицательным действительным числом. Предположим сначала, что ${\Imaginary  {z}<0}$. Тогда ${\Real  iz>0}$. По лемме~\ref{electrical}(\ref{el-2}) получаем, что ${\Real C(iz,i)=\Real  (iC(z,1))=\Real  i=0}$. Так как ${C(iz,i)}$ --- рациональная функция, то любая окрестность точки ${iz}$ содержит такую точку ${z'}$, что ${\Real C(z',i)<0}$. Для достаточно малой окрестности справедливо ${\Real  z'>0}$, потому что ${\Real  iz>0}$. По непрерывности, окрестность точки ${i}$ содержит такую точку ${w'}$, что ${\Real  w'>0}$ и ${\Real C(z',w')<0}$. Эти неравенства противоречат лемме~\ref{electrical}(\ref{el-5}). Случай ${\Imaginary  {z}>0}$ рассматривается похожим образом. Предположим теперь, что ${z>c^2}$. Тогда по лемме~\ref{electrical}(\ref{el-4}) получаем ${1=C(z)>C(c^2)=1}$ --- противоречие. Случай ${0\le z<c^2}$ рассматривается аналогично. Поэтому ${z<0}$.

(\ref{th1-2})$\implies$(\ref{th1-3})
%\begin{lemma}\label{l1} Consider a rational function $$q(z)=\frac{\prod\limits_{k=1}^n (z^2+b_k^2)}{z\prod\limits_{k=1}^m(z^2+a_k^2)},$$ where ${m}$ equals ${n{-}1}$ or ${n}$, ${a_k}$ и ${b_k}$ are real numbers и ${0<b_1<a_1<b_2<a_2<\dots}$. тогда ${\Real q(z)>0}$ если ${\Real z>0}$ и ${\Real q(z)<0}$ если ${\Real z<0.$
%\end{lemma}
Пусть ${p(z)}$ --- минимальный многочлен числа ${c^2}$. Так как минимальный многочлен имеет лишь простые корни, то ${p(z^2)=(z^2-c^2)\prod_{k=1}^n (z^2+b_k^2)}$ для некоторых ${b_1>\dots>b_n>0.}$  Возьмём многочлен ${q(z)}$ с рациональными коэффициентами, такой, что ${q(z)=z\prod_{k=1}^n(z^2+a_k^2)}$, где ${a_1>b_1>a_2>\dots>b_n.}$ Рассмотрим рациональную функцию ${C(z)=q(z)/(zq(z)-p(z^2))}$. Получаем ${C(c)=1/c}$.

Проверим, что функция ${C(z)}$ удовлетворяет условию~(\ref{derivatives}) леммы~\ref{PosReal}.
Ясно, что ${C(z)}$ нечётна и ${\lim_{z\to\infty}C(z)\ne0}$.
Корнями ${C(z)}$ являются числа ${0, \pm ia_1,\dots,\pm ia_n}$.
Непосредственный подсчёт показывает, что для каждого ${l=1,\dots,n}$
$$
C'(\pm ia_l)= -\frac{q'(\pm ia_l)}{p(-a_l^2)} =  \frac{2a_l^2}{(c^2+a_l^2)(a_l^2-b_l^2)}\prod_{k\ne l} \frac{a_k^2-a_l^2}{b_k^2-a_l^2}>0
$$
по предположению, что ${a_1>b_1>a_2>\dots>b_n>0}$. Аналогично ${C'(0)=-q'(0)/p(0)>0}$.
%Thus ${C(z)}$ satisfies condition~(\ref{derivatives}) of Lemma~\ref{PosReal}.

Тогда по лемме~\ref{PosReal}(\ref{derivatives})$\implies$(\ref{continued}) функция ${C(z)}$ удовлетворяет условию~(\ref{continued}) леммы~\ref{PosReal}. Так как ${C(z)\in \mathbb{Q}(z)}$, то по алгоритму Евклида все ${d_k\in \mathbb{Q}}$. Подставляя ${z=c}$, получаем требуемое условие.
\end{proof}
%\begin{claim} The function ${C(z)}$ satisfies condition (3) of Lemma~\ref{PosReal}.
%}$\Real z>0}$ если и only если ${\Real C(z)>0.$
%\end{claim}
%Analyzing the behavior of ${C(z)}$ near the poles we obtain that for each ${y\in\mathbb{R}\setminus \{0\}}$ the equation ${C(ix)=iy}$ has at least ${2n{+}1}$ real solutions и the equation ${C(ix)=0}$ has at least ${2n}$ real solutions. For each ${w\in\mathbb{C}}$ the equation ${C(z)=w}$ has at most ${2n{+}1}$ complex solutions и the equation ${C(z)=0}$ has at most ${2n}$ complex solutions. Therefore for each ${x\in\mathbb{R}}$ all solutions of the equation ${C(iz)=ix}$ are real. Note that ${C(c)>0}$ и ${C(-c)<0}$. Thus by continuity of the function ${C(z)}$ the claim follows.

\section{Вариации}\label{var}

%In this section we study tilings of arbitrary (not necessarily convex) polygons by rectangles. This is equivalent to investigation of planar electrical networks with several terminals.
%This section consists mostly of open problems.

\subsection{Разрезания многоугольников на прямоугольники} \label{polygons}

В этом пункте мы рассмотрим следующую задачу.
%generalization of Problem~\ref{main}.
%The results of \S\ref{var} are proved in \S\ref{finalproofs}.

\begin{problem}\label{general}
Какие многоугольники
%When a given (not necessarily convex) polygon ${P$
можно разрезать на прямоугольники с отношениями сторон ${c_1, \dots, c_n}$?
\end{problem}

Близкая задача о разрезаниях \emph{с учётом знаков} решена в~\cite{KK}.

Частный случай ${n=1}$, ${c_1=1}$ этой задачи, состоящий в описании всех многоугольников, которые могут быть разрезаны на квадраты, поставлен в~\cite{FHTW}.
Р.\7Кеньоном был получен ответ для этого частного случая для шестиугольников (см. рис.~\ref{ell}).
%; in the following theorem a misprint in the original statement is corrected). %, cf.~ \cite{FHTW,KK}.

\begin{figure}[htbp]
\definecolor{zzttqq}{rgb}{0.6,0.2,0}
\begin{tikzpicture}[line cap=round,line join=round,>=triangle 45,x=0.5cm,y=0.5cm]
\clip(-4.88,-1.69) rectangle (5.72,6.85);
\fill[color=zzttqq,fill=zzttqq,fill opacity=0.1] (1,6) -- (-4,6) -- (-4,-1) -- (5,-1) -- (5,4) -- (1,4) -- cycle;
\draw [color=zzttqq] (1,6)-- (-4,6);
\draw [color=zzttqq] (-4,6)-- (-4,-1);
\draw [color=zzttqq] (-4,-1)-- (5,-1);
\draw [color=zzttqq] (5,-1)-- (5,4);
\draw [color=zzttqq] (5,4)-- (1,4);
\draw [color=zzttqq] (1,4)-- (1,6);
\draw [color=zzttqq] (0,6)-- (0,3.04);
\draw [color=zzttqq] (0,-1)-- (0,3.48);
\draw [color=zzttqq] (0,4)-- (1,4);
\draw[color=zzttqq] (-2.35,2.94) node[anchor=north west] {$x$};
\draw[color=zzttqq] (0.05,5.4) node[anchor=north west] {$y$};
\draw[color=zzttqq] (2.2,2) node[anchor=north west] {$z$};
\fill [color=zzttqq] (1,6) circle (0.5pt);
\draw[color=zzttqq] (1.28,6.45) node {$A_2$};
\fill [color=zzttqq] (-4,6) circle (0.5pt);
\draw[color=zzttqq] (-3.96,6.45) node {$A_3$};
\fill [color=zzttqq] (-4,-1) circle (0.5pt);
\draw[color=zzttqq] (-4.34,-1.14) node {$A_4$};
\fill [color=zzttqq] (5,-1) circle (0.5pt);
\draw[color=zzttqq] (5.32,-1.06) node {$A_5$};
\fill [color=zzttqq] (5,4) circle (0.5pt);
\draw[color=zzttqq] (5.32,4.4) node {$A_6$};
\fill [color=zzttqq] (1,4) circle (0.5pt);
\draw[color=zzttqq] (1.46,4.4) node {$A_1$};
\end{tikzpicture} 
\caption{$L$-образный шестиугольник.}
\label{ell}
\end{figure}

\begin{theorem} \label{th-kenyon}
\textup {(ср.~\cite[Theorem~9]{K})} Пусть $A_1A_2A_3A_4A_5A_6$ --- $L$-образный шестиугольник, вершины которого занумерованы, начиная с вершины невыпуклого угла \textup {(}см. рис.\7\ref{ell}\textup {)}. Шестиугольник $A_1A_2A_3A_4A_5A_6$ можно разрезать на квадраты, если и только если система линейных уравнений
\begin{equation}
\label{syst}
\begin{cases}
A_3A_4\cdot x+ A_1A_2\cdot y&=A_2A_3,\\
A_5A_6\cdot z-A_1A_2\cdot y&=A_6A_1;
\end{cases}
\end{equation}
имеет решение в неотрицательных рациональных числах $x,y,z$.
\end{theorem}

Доказательство достаточности довольно просто: если рациональные числа ${x,y,z\ge 0}$ удовлетворяют системе\7(\ref{syst}), то шестиугольник $A_1A_2A_3A_4A_5A_6$ можно разрезать на ${3}$ прямоугольника с отношениями сторон ${x,y,z}$, а затем и на квадраты (см. рис.~\ref{ell}).

Наша ближайшая цель --- дать подобный теореме\7\ref{th-kenyon} критерий для гораздо более широкого класса многоугольников.

%Let us give a necessary condition for an arbitrary ``generic'' polygon ${P}$ to be tilable by squares.

Пусть $P$ --- \emph{ортогональный} многоугольник, то есть многоугольник со сторонами, параллельными координатным осям. Будем считать, что многоугольник $P$ --- \emph{простой}, то есть его граница $\partial P$ состоит из одной связной компоненты. Занумеруем стороны многоугольника, параллельные оси $x$, против часовой стрелки в порядке их следования на границе $\partial P$.
Обозначим через $b$ число горизонтальных сторон у многоугольника~ $P$.
Обозначим через $I_u$ \emph{ориентированную длину} стороны $u$, то есть длину стороны $u$ со знаком \2$+$\3 (соответственно, \2$-$\3), если многоугольник $P$ локально лежит ниже (соответственно, выше) стороны $u$. Обозначим через $U_u$ координату стороны $u$ по оси $Oy$. Будем считать, что многоугольник $P$ --- \emph{общего положения}, то есть числа ${U_1, \dots, U_b}$ попарно различны. (Понятие \2многоугольника общего положения\3 используется только в указанном смысле, а не в смысле какой-либо меры в пространстве всех многоугольников.)

Нам потребуются следующие понятия \cite{CM}. Назовём последовательность
${(p_1, \dots, p_{2k})}$
целых чисел
%сторон многоугольника $P$
\emph{циклически упорядоченной}, если %a циклическая перестановка
некоторая её циклическая перестановка строго упорядочена по возрастанию.
%Например, для многоугольника на рисунке~\ref{kenyon} последовательность сторон $(c+d,a,d,a+b,c,b)$ --- циклическая.

% (то есть, увеличивающаяся последовательность).
Обозначим через $\Omega_b$ множество матриц $C_{uv}$ размера ${b\times b}$ с действительными элементами, удовлетворяющих следующим условиям:
\begin{enumerate}[(1)]
\item \label{omega-1} матрица $C_{uv}$ симметрична, то есть ${C_{uv}=C_{vu}}$;
\item \label{omega-2} сумма элементов матрицы $C_{uv} $ в каждой строке равна нулю;
\item \label{omega-3} для любой циклически упорядоченной последовательности ${(p_1, \dots, p_{2k})}$ выполнено неравенство
${(-1)^k \det\{ C_{p_ip_{2k-j+1}} \}_{i,j=1}^k\ge0}$.
\end{enumerate}
Последнее условие есть ограничение на знаки некоторых миноров матрицы $C_{uv}$. Например, множество $\Omega_3$ состоит из матриц вида
${\small\left(\begin{matrix}
x+y & -x & -y\\
-x  & x+z& -z\\
-y  & -z &y+z
\end{matrix}\right) }$,
где ${x,y,z>0}$.

%We need the following notation \cite{CM}.
%Let ${\Omega_{b}}$ be the set of symmetric ${b\times b}$ matrices ${C_{uv}}$ with zero sum of each column satisfying the following condition: for each ${1\le k\le l\le b}$ all the minors of the submatrix ${\{-C_{uv}\}_{k\le u\le l, \ v<k \text{ or }v>l}}$ are nonnegative. Юnote ${-\bar\Omega_{b}=\{-C:C\in\bar\Omega_{b}\}}$.

%Here ${2b}$ denotes the number of vertices in ${P}$.

\begin{theorem} \label{squaring-polygons}
Пусть $P$ --- ортогональный многоугольник общего положения с $b$ горизонтальными сторонами, имеющими ориентированные длины ${I_1, \dots, I_b}$ и ${y\mbox{-ко}}$орди\-на\-ты ${U_1, \dots, U_b}$.
%Suppose, что многоугольник $P$ может крыться черепицей квадратами; тогда есть матрица $C _ {uv} \in \Omega _ {b} $ с рациональными записями таким образом что $I_v =\sum _ {u=1} ^ {b} C _ {uv} U_u$ для каждого $v=1, \dots, b$.
Тогда следующие два условия эквивалентны:
\begin{enumerate}[(1)]
\item \label{s-p-1} многоугольник $P$ можно разрезать на квадраты;
\item \label{s-p-2} существует матрица ${C_{uv}\in \Omega_{b}}$ с рациональными элементами, такая, что ${I_v=\sum_{u=1}^{b} C_{uv}U_u}$ для каждого ${v=1,\dots, b}$.
\end{enumerate}
\end{theorem}

%\begin{problem} Is the converse theorem true for ${b\ge 4}$?
%\end{problem}

При ${b=2}$ эта теорема эквивалентна теореме~\ref{dehn}, а при
${b=3}$ --- теореме~\ref{th-kenyon}.
%\begin{example}\label{nongeneric}

Обсудим алгоритмические аспекты данного результата.

Теорема ~\ref{squaring-polygons} сводит вопрос о существовании разрезания данного многоугольника на квадраты к вопросу о существовании рационального решения у некоторой системы нелинейных неравенств. Последний вопрос очень сложен.  Известно, что в общем случае не существует алгоритма для ответа на него.
Поэтому Теорема~\ref{squaring-polygons} не даёт алгоритма, который бы по \emph{произвольному} многоугольнику позволил бы определить, можно ли его разрезать на квадраты.

Однако можно выделить \emph{обширный класс} многоугольников, для которого теорема~\ref{squaring-polygons} даёт такой алгоритм. А именно, рассмотрим такие многоугольники $P$, для которых числа ${U_1, \dots, U_b}$ линейно независимы над $\mathbb {Q}$. Тогда условие ${I_v =\sum _ {u=1} ^ {b} C _ {uv} U_u}$ из теоремы~\ref{squaring-polygons} однозначно определяет матрицу $C_{uv}$. Чтобы выяснить, можно ли многоугольник $P$ разрезать на квадраты, достаточно проверить, удовлетворяет ли матрица $C_{uv}$ условиям (\ref{omega-1})--(\ref{omega-3}) из определения множества $\Omega_b$.

Теорема \ref{squaring-polygons} % (1) $\implies$ (2)
не обязательно выполнена для многоугольников \emph{не общего положения}, например, для ортогонального многоугольника с параметрами
$$U_1=U_3=0, \quad U_2=2, \quad U_4 =-4, \quad I_1 =\sqrt {2}, \quad I_2=2, \quad I_3=2-\sqrt {2}, \quad I_4 =-4.$$

Мы докажем теорему~\ref{squaring-polygons} в \S\ref{finalproofs}. Мы также дадим короткое доказательство следующего результата.

\begin{theorem} \textup{\cite{SuDi}} \label{polygon-square} Ортогональный многоугольник общего положения с рациональными координатами вершин можно разрезать на прямоугольники с отношениями сторон $c$ и ${1/c}$ тогда и только тогда, когда квадрат можно разрезать на прямоугольники с отношениями сторон $c$ и~${1/c}$.
\end{theorem}

%We prove Theorems~\ref{squaring-polygons}  %\ref{rectangling-polygons}
%and~\ref{polygon-square} in the end of this section.

%\begin{theorem} \label{th3}
%Suppose that a polygon of area ${S}$ can be dissected into rectangles of ratios ${R_1,\dots,R_n}$. Suppose that the $x$-coordinates of all vertices of the polygon are rational. тогда there exists a function ${E(z_1,\dots,z_n)\in\mathbb{Q}(z_1,\dots,z_n)}$ such that the following conditions hold:

%(1) ${E(R_1, R_2, \dots, R_n)=S}$;

%(2) ${E(cz_1,cz_2, \dots, cz_n)\equiv cE(z_1,z_2, \dots, z_n)}$;

%(3) если ${\Real z_k>0}$ for all ${k=1,\dots,n}$ тогда ${\Real E(z_1,z_2, \dots, z_n)>0}$.
%\end{theorem}

%\begin{conjecture} The converse theorem is also true.
%\end{conjecture}

%\subsection{Electrical circuits with several terminals} \label{responsesection}

\subsection{Электроимпедансная томография}\label{eit}

Наш подход к задаче \ref{general} следует идеям из \cite{K, CGV} и использует электрические цепи с несколькими входами.

Далее будем рассматривать замкнутые электрические цепи с несколькими граничными вершинами ${1}$, \dots, ${b}$ и с входными потенциалами ${U_1}$, \dots, ${U_b}$.
Замкнутую электрическую цепь будем называть \emph{плоской}, если её можно изобразить без самопересечений в единичном круге так, что граничные вершины лежали бы на границе круга. В этом случае будем нумеровать граничные вершины против часовой стрелки, если следовать по границе круга. Мы \emph{не} предполагаем, что замкнутая электрическая цепь связна, но требуем, чтобы каждая связная компонента цепи содержала хотя бы одну граничную вершину. Потенциалы и токи в такой цепи определяются законом Ома~(C) и правилом Кирхгофа~(I) из~\S\ref{prelim}.
%Let ${1,\dots,b}$ be the boundary vertices of the electrical network., и ${b+1,\dots,n}$ are the remaining ({\it interior}) ones. Let ${U_1,U_2,\dots,U_b}$ be the given voltages at the boundary vertices, и ${U_{b+1},\dots,U_n}$ be the voltages at the interior vertices determined by the Ohm и Kirchhoff laws.
%Let ${r_{ij}}$ be the conductance of the edge between nodes ${i}$ и ${j}$.

Рассмотрим линейное отображение ${\mathbb{C}^b\to \mathbb{C}^b}$, переводящее вектор потенциалов ${(U_1,\dots,U_b)}$ в вектор \emph{входных токов} ${(I_1, \dots, I_b)=(\sum_{k=1}^{n} I_{1k},\dots,\sum_{k=1}^{n} I_{bk})}$, текущих внутрь цепи через вершины ${1, \dots, b}$ соответственно. Матрица ${C_{uv}}$ этого линейного отображения называется {\it реакцией} замкнутой электрической цепи.
Эта матрица симметрична~\cite{CM}. При ${b=2}$ реакция цепи равна ${{\tiny\left(\begin{matrix} C & -C\\ -C & C\end{matrix}\right)}}$,
где $C$ --- проводимость цепи.
%(This definition differs by a sign from the one in \cite{C}.)
%We have the following generalization of the results of \S\ref{prelim}, cf.~\cite{Duffin}.

%The value ${E=\sum_{1\le u,v\le b} C_{uv} U_u\bar U_v}$ is called {\it energy dissipation}.
%We consider ${E(r_{ij})}$ as the function of the conductances ${r_{ij}}$ for fixed voltages ${v_i}$.
%A symmetric ${b\times b}$ matrix $c$ is \textit{almost positively definite}, если it is non-negatively definite и ${x^{T}Cx=0}$ only если ${Cx=0}$.
%The following (probably new) results are proved completely analogously to the proof of the results of \S\ref{prelim} above.

%\subsection{Relation with tilings}

%\subsection{Electrical impedance tomography}

%Analogously to the continuous case the following problems for discrete harmonic functions are natural to study.
%The rest of \S\ref{var} consists of open problems и a few examples.

Мы сводим результаты \S\ref{polygons} к
\emph{обратной задаче} теории электрических цепей, которая состоит
%following problems even more interesting in themselves:
%\begin{itemize}
%\item \emph{Direct problem.} Describe possible responses of electrical networks.
%\item \emph{Inverse problem.}
в определении неизвестных сопротивлений элементов электрической цепи по измерениям силы тока и потенциалов в заданных узлах.
%to describe all the networks having a given response.
%\end{itemize}
Она представляет собой дискретный аналог
\emph{электроимпедансной томографии} \cite{Cal,SU}. Эта задача поставлена в \cite{LSy} и решена для плоских цепей постоянного тока в \cite{CIM, CM, C, CGV}. Сформулируем некоторые глубокие результаты И.~Колина де~Вердье, Э.~Куртиса и Дж.~А.~Морроу.
%For instance, the following deep results are known.
%Let us state the following deep results.

%We need the following definitions \cite{C,CM}. Let ${G}$ be a network и ${u_1,\dots,u_k,v_1,\dots,v_k}$ be a sequence of boundary vertices. The sequence is \emph{connected through ${G$} если ${G}$ contains ${k}$ disjoint paths ${p_1,\dots,p_k}$ such that ${p_l}$ starts at ${u_l}$, ends at ${v_l}$ и passes through no other boundary vertices. Denote by ${\pi(G)}$ be the set of all the sequences connected through ${G}$. The network ${G}$ is \emph{minimal} если for each network ${G^-}$ obtained from ${G}$ by either deleting or contracting an edge we have ${\pi(G^-)\ne \pi(G)}$. Юnote by ${\mathbb{R}_+}$ the set of positive real numbers.

\begin{theorem} \label{CGV5} %\textup{\cite[Proposition 10]{C}}
\textup{\cite[Theorem~5]{CM, CIM, CGV}}
%The response of a planar electrical network with ${b}$ boundary vertices и positive edge conductances belongs to the set ${\Omega_b}$. Each element of the set ${\Omega_b}$ is the response of a minimal planar electrical network with ${b}$ boundary vertices и positive edge conductances.
Множество всех возможных реакций плоских электрических сетей с $b$ граничными вершинами и положительными проводимостями рёбер --- это в точности множество $\Omega_b$ из~\S\ref{polygons}.
%The following ${3}$ conditions are equivalent:
%\begin{enumerate}[(1)]
%\item ${C_{uv}}$ is the response of a planar electrical network with ${b}$ boundary vertices и positive edge conductances;
%\item ${C_{uv}}$ is the response of a \textup{critical} planar electrical network with ${b}$ boundary vertices и positive edge conductances;
%\item ${C_{uv}\in\Omega_b}$.
%\end{enumerate}
\end{theorem}

Назовём электрическую цепь \emph{минимальной} (или \emph{критической}), если она имеет минимальное число рёбер среди всех плоских электрических цепей, имеющих такую же реакцию, как и исходная цепь. Минимальность цепи зависит только от её графа\7\cite{CGV}. В\7\cite[\S9]{CM, CIM} приводится алгоритм нахождения проводимостей рёбер минимальной цепи по её реакции. Этот алгоритм влечёт следующий результат (Д.~Ингерман, Дж.~А.~Морроу, частное сообщение).

\begin{theorem} \label{inverse}
%\textup{\cite[\S6.4]{CM}}
Проводимости рёбер минимальной электрической цепи единственным образом определяются по её реакции.  Проводимость каждого ребра есть рациональная функция с рациональными коэффициентами
от элементов матрицы реакции.
\end{theorem}

Для цепей \emph{переменного} тока обратная задача, вероятно, является нерешённой \cite{Foster}.
Сформулируем простой фольклорный результат. Оставшаяся часть~\S\ref{var} не используется в доказательстве сформулированных выше теорем. Однако авторы считают, что для каких-либо дальнейших продвижений в задачах~\ref{mainproblem} и~\ref{general} необходимы
%the solution of the inverse problem.
обобщения результатов, сформулированных ниже.

%Let us state some results for alternating-current circuits (the rest of \S\ref{var} is не used in the proofs of the above results). %of \S\ref{polygons}).

% The following notion is inappropriate because
% we need ${C U=0}$ instead of ${Re C=0}$ !!!
%A real ${b\times b}$ matrix ${R}$ is \emph{almost positively definite}, если for any ${U\in\mathbb{R}^b}$ we have either ${U^T R U>0}$ or else ${R U=0}$.

\begin{theorem}\label{detection} При ${b=2}$ или ${b=3}$ следующие два условия эквивалентны\textup{:}
\begin{enumerate}[(1)]
\item \label{det-1} ${C_{uv}}$ является реакцией связной электрической цепи с ${b}$ граничными вершинами, проводимости рёбер которой имеют положительную действительную часть\textup{;}
\item \label{det-2} ${C_{uv}}$ --- комплексная ${b\times b}$ матрица, обладающая следующими свойствами\textup{:}
\begin{itemize}
\item \emph{симметричность}\textup{:} ${C_{uv}=C_{vu}}$\textup{;}
\item \emph{сумма элементов в каждой строке равна нулю}\textup{:} ${\sum_{1\le u\le b}C_{uv}=0}$\textup{;}
\item \emph{положительность вещественной части}\textup{:} ${\sum_{1\le u,v\le b}(\Real C_{uv})U_u U_v\ge 0}$ для любых ${U_1,\dots,U_b\in\mathbb{R}}$\textup{;}
\item \emph{почти положительная определённость}\textup{:}    последнее неравенство является строгим, кроме случая ${U_1=\dots=U_b}$.
\end{itemize}
\end{enumerate}
\end{theorem}

\begin{problem} Верен ли этот результат для ${b\ge 4}$?
\end{problem}

%Theorem~\ref{detection} is proved at the end of the section.

%The problem of \textit{determining the structure of a network by its response} might be of practical interest.
%One can see that for each direct-current network there is a (not necessarily planar) direct-current network without nonboundary vertices и with the same response. Whereas

%Unlike nonplanar direct-current networks  \emph{nonboundary vertices in alternating-current networks can be detected by the response}.
%For instance, by Theorem~\ref{detection} there are electrical networks with response
%{\tiny ${\left(\begin{matrix}
%2 & 1 & -3\\
%1 & 2 & -3\\
%-3& -3& 6 \\
%\end{matrix}\right)$}; any such network
%necessarily has nonboundary vertices.

%(Indeed, the response of a network without nonboundary vertices cannot have positive entries outside the diagonal.)
%\end{example}

%A result concerning the inverse problem is presented in \S\ref{inverse}.

%\begin{problem}
%Which is the least possible number of vertices in an electrical network with a given response?
%\end{problem}

%\begin{problem} Suppose that ${C_{uv}(c_{12}, c_{13}, \dots, c_{n-1,n} )}$ is the response of an electrical network with ${b}$ boundary vertices и with edge conductances ${c_{12}, c_{13}, \dots, c_{n-1,n}}$. Let ${C_{uv}}$ be a complex ${b\times b}$ matrix. Is it true that the equation ${C_{uv}(c_{12}, c_{13}, \dots, c_{n-1,n} )=C_{uv}}$ has either no or unique or infinitely many (continuum of) solutions ${(c_{12}, c_{13}, \dots, c_{n-1,n} )\in\mathbb{C}^{b(b-1)/2}}$?
%\end{problem}

\subsection{Случайные блуждания}

{\it Случайным блужданием} по электрической цепи (или по \emph{взвешенному графу}) называется марковская цепь с матрицей перехода ${P_{kl}=c_{kl}/\sum_{j=1}^{n}c_{jk}}$. %Here we require the same assumptions on the graph ${G}$ и the numbers ${c_{kl}}$ as in \S\ref{prelim}. Under these assumptions the
Такая марковская цепь является эргодической и обратимой.
Обозначим через ${k_1l_1,\dots,k_ml_m}$ все рёбра марковской цепи.
Следующая теорема позволяет перевести результаты \S1--\S2 на \2язык\3 случайных блужданий.

\begin{theorem} \textup{\cite[с. 42]{DS}}
Пусть ${P(c_{k_1l_1}, \dots,c_{k_ml_m})}$ --- вероятность того, что случайное блуждание, стартующее из вершины ${1}$, достигнет вершины ${2}$ до первого возвращения в вершину ${1}$. Пусть ${C(c_{k_1l_1}, \dots,c_{k_ml_m})}$ --- проводимость электрической цепи \textup{(}с граничными вершинами ${1}$ и ${2}$\textup{)}.
Тогда $${P(c_{k_1l_1}, \dots,c_{k_ml_m})={C(c_{k_1l_1}, \dots,c_{k_ml_m})}/{(c_{12}+\dots+c_{1n})}}.$$
\end{theorem}

%\begin{figure}[htbp]
%\includegraphics{walksnetworks.pdf}
%\caption{Walks и networks}
%\label{walksnetworks}
%\end{figure}

%по лемме~\ref{electrical}(1) ${P(c_{kl})}$ is a rational function in ${c_{kl}}$. A ``translation'' of Lemma~\ref{electrical}(5) is:

%\begin{theorem}
%If ${\Real c_{kl}>0}$ for all ${1\le k<l\le n}$ тогда ${\Real \left((c_{12}+\dots+c_{1n})P(c_{kl})\right)>0}$.
%\end{theorem}

%\begin{problem} Does this result remain true for nonreversible Markov chains?
%\end{problem}

Например, переводом лемм~\ref{electrical}(\ref{el-1}) и~(\ref{el-5}) будет следующий результат.

\begin{corollary} \label{random}
Вероятность ${P(c_{k_1l_1}, \dots,c_{k_ml_m})}$ есть рациональная функция от нge weights
${c_{k_1l_1}, \dots,c_{k_ml_m}}$, которая обладает следующим свойством\textup{:} если ${\Real c_{k_1l_1}, \dots,\Real c_{k_ml_m}>0}$, то ${\Real \left((c_{12}+\dots+c_{1n})P(c_{k_1l_1}, \dots,c_{k_ml_m})\right)>0}$.
\end{corollary}

Последний результат не обязательно выполнен для \emph{необратимых} марковских цепей: например, для марковской цепи с вершинами ${1}$, ${2}$, ${3}$, ${4}$ и ориентированными рёбрами ${14}$, ${42}$, ${43}$.
%we have ${P(c_{14},c_{42},c_{43})=c_{42}/(c_{42}+c_{43})}$.
%A natural question is: \emph{does this result remain true for nonreversible Markov chains}?

В \cite{Fomin} получен близкий результат для необратимых марковских цепей. Доказательство использует результаты электроимпедансной томографии, сформулированные в~\S\ref{eit}.

Необратимые плоские марковские цепи могут быть интерпретированы, как разрезания трапеций на трапеции \cite{K} (здесь \emph{трапецией} называется четырёхугольник, две стороны которого параллельны оси\7$x$). %The \emph{ratio} of the trapezoid is the length of the horizontal middle edge divided by the hight. 
Естественными задачами являются следующие: \emph{обобщение результатов статьи на
%\begin{itemize}
разрезания на трапеции; бесконечные разрезания; разрезания с учётом знаков}.
%\end{itemize}
%Probably, electrical networks with edges having different conductances in different directions (or conductances of arbitrary sign --- in case of signed tilings) might be useful for answering the question.

\section{Обобщение основных идей}\label{final}

\subsection{Электрические цепи} \label{elcirc}

%\subsection{Proof of the results of \S\ref{prelim}}

Наш подход основан на обобщении результатов \S\ref{prelim} на электрические цепи с ${b}$ входами.
Короткие доказательства результатов \S\ref{prelim} получаются в этом параграфе как частный случай ${b=2}$. Наше доказательство леммы~\ref{responseclaim}(\ref{rc-3}), которая является обобщением леммы~\ref{electrical}(\ref{el-3}), вероятно, является новым.
Для простоты в этом пункте предполагаем, что цепи не имеют кратных рёбер.
%We prove the results of \S\ref{prelim} in this general situation.
%Our proofs of Lemma~\ref{correspondence-general} и Lemma~\ref{responseclaim}(3), generalizing
%Lemma~\ref{correspondence-lemma} и Lemma~\ref{electrical}(3), respectively, are probably new.
Все доказательства основаны на следующем фундаментальном \textit{законе сохранения энергии}.

\begin{lemma} \label{energy} \label{fullpower} Пусть ${E(U,I)}$ --- билинейная или полубилинейная функция.
Рассмотрим электрическую цепь с вершинами ${1, \dots, n}$, такую, что вершины ${1, \dots, b}$ --- граничные.
Предположим, что числа ${U_k}$, где ${1\le k\le n}$, и ${I_{kl}}$, где ${1\le k,l\le n}$, удовлетворяют закону Ома~\textup{(C)} и правилу Кирхгофа~\textup{(I)} из \S\ref{prelim}. Положим ${I_u=\sum_{k=1}^n I_{uk}}$. Тогда
$$\sum\limits_{1\le k<l\le n}E(U_k-U_l,I_{kl})=\sum\limits_{1\le u\le b}E(U_u,I_u).$$
\end{lemma}

%\begin{claim} \label{fullpower} Let ${E(U,I)}$ be a bilinear function.
%Suppose that the numbers ${U_k}$, ${2< k\le n}$, и ${I_{kl}}$, ${1\le k,l\le n}$ satisfy the laws \textup{(C),(I)}. Let ${I_1=\sum_{k=2}^n I_{1k}}$. тогда
%$$\sum\limits_{1\le k<l\le n}E(U_k-U_l,I_{kl})=E(U_1-U_2,I_1).$$
%\end{claim}

Обычно мы будем применять эту лемму к функции \emph{диссипации энергии} ${E(U,I)=\Real (U\bar I)}$.

\begin{proof}[Доказательство леммы \ref{energy}] По закону\7(C) справедливо ${I_{lk}=-I_{kl}}$. Следовательно, по правилу\7(I) справедливо
$$
\sum\limits_{1\le k<l\le n}E(U_k-U_l,I_{kl})=\sum\limits_{k=1}^n E(U_k,\sum\limits_{l=1}^n I_{kl})\sum\limits_{1\le u\le b}E(U_u,I_u).
$$
\end{proof}

%\begin{proof}[Proof of Claim \ref{fullpower}] Clearly, by law (I) we have ${\sum_{k=1}^n I_{2k}=-I_1}$. Thus
%$$\sum\limits_{1\le k<l\le n}E(U_k-U_l,I_{kl})=\sum\limits_{k=1}^n E(U_k,\sum\limits_{l=1}^n I_{kl})=E(U_1-U_2,I_1).$$
%\end{proof}

%\begin{claim} (Energy preserving law) \label{fullpower} Let the numbers ${U_k}$, ${3\le k\le n}$, и ${I_{kl}}$, ${1\le k,l\le n}$ satisfy the Ohm и the Kirchhoff laws (1)--(2). Let ${I=-\sum_{k=2}^n I_{1k}}$. тогда
%}$\sum\limits_{1\le k<m\le n}(U_m-U_k)\bar{I}_{km}=(U_2-U_1)\bar{I}}$.
%\end{claim}

%\begin{proof} Clearly, by Kirchhoff law (2) we have ${\sum_{k=1}^n I_{2k}=I}$. Thus
%\sum\limits_{1\le k<m\le n}(U_m-U_k)\bar{I}_{km}=\sum\limits_{k=1}^n U_k\sum\limits_{m=1}^n (-\bar{I}_{km})=(U_2-U_1)\bar{I}.
%\end{proof}

Докажем теорему~\ref{weyl} в более общей ситуации, когда у электрической цепи ${b}$ граничных вершин, а проводимости рёбер --- комплексные числа с положительной действительной частью.

\begin{proof}[Доказательство теоремы \ref{weyl}] \textit{Единственность}. %Suppose that for each ${k,l}$ either ${\Real c_{kl}>0}$ or ${c_{kl}=0}$.
Предположим, что два набора токов ${I^{I,II}_{kl}}$ и потенциалов ${U^{I,II}_k}$ подчиняются (C) и (I). Тогда их разность ${I_{kl}=I^{I}_{kl}-I^{II}_{kl}}$ и ${U_{k}=U^{I}_{k}-U^{II}_{k}}$ подчиняется (C) и (I) при нулевых входных потенциалах ${U_1=\dots=U_b=0}$. Тогда по %Claim~\ref{fullpower}
лемме~\ref{energy} имеем
$$
\sum\limits_{1\le k<l\le n}\Real \bar c_{kl}|U_k-U_l|^2=\sum\limits_{1\le k<l\le n}\Real ((U_k-U_l)\bar{I}_{kl})\sum\limits_{1\le u\le b}\Real (U_u\bar I_u)=0.
$$
Для всех ${k,l}$ либо ${\Real c_{kl}>0}$, либо ${c_{kl}=0}$. Таким образом, каждое слагаемое ${\Real \bar c_{kl}|U_k-U_l|^2=0}$. Поскольку любая компонента связности цепи содержит граничную вершину, то все потенциалы ${U_{k}}$ равны между собой. Следовательно, ${U_k=0}$, ${I_{kl}=0}$, и поэтому ${I^{I}_{kl}=I^{II}_{kl}}$, ${U^{I}_{k}=U^{II}_{k}}$ для всех ${k,l}$.

\textit{Существование}. Количество уравнений в системе (C), (I) равно количеству неизвестных.
Мы доказали, что
%By the above previous paragraph
эта система имеет единственное решение при ${U_1=\dots=U_b=0}$. По конечномерной альтернативе Фредгольма она имеет решение для любых ${U_1}$, \dots, ${U_b}$.
\end{proof}

%\begin{remark} The assumption ${\Real c_{kl}>0}$ in Theorem~\ref{weyl} is essential by example in figure~\ref{examples} for ${a=i,b=-i}$.
%\end{remark}

%Let us summarize the properties of the formulas for response.
Следующий результат обобщает лемму~\ref{electrical}.

\begin{lemma} \label{responseclaim} Предположим, что электрическая цепь имеет ${b}$ граничных вершин и ${m}$ рёбер с проводимостями ${c_{1},\dots, c_{m}}$. Тогда реакция цепи ${C_{uv}(c_{1},\dots, c_{m})}$ обладает следующими свойствами\textup{:}
\begin{enumerate}[(1)]
\item \label{rc-1}
    %\emph{rationality}\textup{:}
    ${C_{uv}(c_{1},\dots, c_{m})\in \mathbb{Q}(c_{1},\dots, c_{m})^{b \times b}}$\textup{;}
    %\item ${C_{uv}(c_{1},\dots, c_{m})}$ is linear-fractional with respect to each variable ${c_{j}}$;
\item \label{rc-2}
    %\emph{homogeneity}\textup{:}
    ${C_{uv}(tc_{1},\dots, tc_{m})=tC_{uv}(c_{1},\dots, c_{m})}$\textup{;}
\item \label{rc-3}
    ${\frac{\partial}{\partial c_{j}}C_{uv}(c_{1},\dots, c_{m}) =(V_{ku}-V_{lu})(V_{kv}-V_{lv})}$, где ${k}$ и ${l}$ являются концами ребра ${j}$, а ${\left(V_{pq}\right)^{n,b}_{p,q=1}}$ --- матрица линейного отображения ${(U_1,\dots,U_b)\mapsto (U_1,\dots,U_n)}$\textup{;}
\item \label{rc-4}
    %\emph{monotonicity}\textup{:}
    если ${c_{1},\dots, c_{m}>0}$, то ${\frac{\partial}{\partial c_{j}}C_{uv}(c_{1},\dots, c_{m})}$ неотрицательно определена\textup{;}
%\item ${C_{uv}(c_{12}, c_{13}, \dots, c_{n-1,n} )}$ is symmetric;
\item \label{rc-5}
    %\emph{}\textup{:}
    если ${\Real c_{1},\dots,\Real c_{m}>0}$, то ${\Real C_{uv}(c_{1},\dots, c_{m})}$ неотрицательно определена.
%\item \label{rc-6} ${C_{uv}(c_{1},\dots, c_{m})=C_{vu}(c_{1},\dots, c_{m})}$.
\end{enumerate}
\end{lemma}

\begin{proof}[Доказательство леммы \ref{responseclaim}] (\ref{rc-1})
По теореме~\ref{weyl} и методу Крамера решение ${\{I_{kl}(U_1,\dots,U_b)\}}$ системы линейных уравнений (C), (I)
состоит из линейных функций от ${U_1, \dots, U_b}$ с коэффициентами, являющимися рациональными функциями от ${c_{1}, \dots, c_m}$.
Поэтому элементы матрицы линейного отображения ${(U_1,\dots,U_b)\mapsto{\sum_{k=1}^n I_{uk}(U_1,\dots,U_b)}}$ являются рациональными функциями от ${c_{1}, \dots, c_m}$.

%(1) Fix an edge ${j}$. Consider the matrix ${M(c_j)}$ of the system of linear equations (C), (I) as a function in the parameter ${c_{j}}$. The matrix ${M(c_{j})}$ contains a ${2\times 2}$ minor
%{\tiny ${\left(
%\begin{matrix}
%-c_{j} & c_{j} \\
%c_j    & -c_{j}
%\end{matrix}
%\right)$}.
%The entries of ${M(c_{j})}$ outside this ${2\times 2}$ minor do не depend on ${c_{j}}$. Hence ${\det M(c_{j})}$ being the denominator of ${C_{uv}(c_{1},\dots, c_{m})}$ is linear in ${c_{j}}$. Analogously, the nominator is linear in ${c_{j}}$.

(\ref{rc-2}) Рассмотрим систему линейных уравнений, полученную из условий (C), (I) подстановкой ${t c_{1},\dots, t c_{m}}$ вместо ${c_{1},\dots, c_{m}}$. Она определяет те же потенциалы, что и в исходной цепи, а токи изменяются в ${t}$ раз. Поэтому ${C(tc_{1},\dots, tc_{m})=tC(c_{1},\dots, c_{m})}$. %\sum\limits_{k=2}^n tc_{1k}\frac {U_1-U_k}=tC(c_{kl})}$.

(\ref{rc-3})
Зафиксируем потенциалы ${U_1,\dots,U_b}$ и все проводимости ${c_{pq}}$, кроме ${c_{kl}}$.
Рассмотрим потенциалы и токи в цепи как функции от ${c_{kl}}$. Положим ${E(U,I):=U\frac{\partial I}{\partial c_{kl}}-\frac{\partial{U}}{\partial c_{kl}}I}$. %(\textit{Lagrange function}).
Тогда ${E(U_k-U_l,I_{kl})=(U_k-U_l)^2}$ и ${E(U_p-U_q,I_{pq})=0}$ для ${pq\ne kl}$.
Поэтому по лемме~\ref{fullpower} имеем
\begin{multline*}\sum_{1\le u,v\le b}\frac{\partial C_{uv}}{\partial c_{kl}}U_uU_v=\sum_{1\le u\le b}E(U_u,I_u)=\sum\limits_{1\le p<q\le n}E(U_p-U_q,I_{pq})=\\
=(U_k-U_l)^2=\sum_{1\le u,v\le b}(V_{ku}-V_{lu})(V_{kv}-V_{lv})U_uU_v.
\end{multline*}

%Let ${E(U,I)=\frac{\partial{U}}{\partial c_{uv}}I-U\frac{\partial I}{\partial c_{uv}}}$. %(\textit{Lagrange function}).
%тогда ${E(U_u-U_v,I_{uv})=(U_u-U_v)^2}$ и ${E(U_k-U_l,I_{kl})=0}$ for ${kl\ne uv}$.
%Thus by Claim~\ref{fullpower} we have
%$$\frac{\partial C}{\partial c_{uv}}(U_1-U_2)^2=E(U_1-U_2,I)=\sum\limits_{1\le k<l\le n}E(U_k-U_l,I_{kl})=(U_u-U_v)^2.$$

(\ref{rc-4}) Доказательство напрямую вытекает из предыдущей формулы.

(\ref{rc-5}) Для всех ${k,l}$ имеем либо ${\Real c_{kl}>0}$, либо ${c_{kl}=0}$. Возьмём ${U_1,\dots,U_b\in\mathbb{R}}$ (хотя доказательство ниже работает и для комплексных потенциалов).
По лемме~\ref{energy} имеем
\begin{multline*}
\sum_{1\le u,v\le b} (\Real C_{uv})U_u \bar U_v\sum_{1\le u,v\le b} \Real(U_u \bar C_{uv}\bar U_v)\sum_{1\le u\le b} \Real (U_u \bar I_u)=\\
=\sum\limits_{1\le k<l\le n}\Real ((U_k-U_l)\bar I_{kl})\sum\limits_{1\le k<l\le n}\Real c_{kl}|U_k-U_l|^2\ge 0.
\end{multline*}
\end{proof}

\begin{remark} \label{strict} Для связной цепи последнее неравенство является строгим, кроме единственного случая, когда ${U_1=\dots=U_b}$. Ещё одно свойство матрицы реакции --- ${C_{uv}=C_{vu}}$ (что также доказывается с помощью подходящего рассуждения о сохранении энергии).
\end{remark}

%\begin{remark}
%The classical proof of Lemma \ref{main}(4) is based on the following observation: \emph{if all ${c_{kl}>0}$ тогда the solution  ${(U_{b+1},\dots,U_n)}$ of the system of linear equations (C), (I) minimizes the energy dissipation function ${\sum_{1\le k<l\le n}c_{kl}(U_k-U_l)^2}$ for fixed ${U_1,\dots,U_b$}. \cite{Rayleigh}
%\end{remark}

%\subsection{Inverse problem} \label{inverse}

%The following result is probably new.

%\begin{lemma} \label{inverse-lemma} не CORRECT!!! Suppose that ${C_{uv}(c_1, \dots, c_{m})}$ is the response of an electrical network with ${b}$ boundary vertices и with ${m}$ edges of conductances ${c_{1}, \dots, c_{m}}$. тогда there is an integer number ${0\le r\le m}$ и ${m-r}$ functions
%$$
%c_{k}(c_1,\dots,c_r,C_{11},C_{12},\dots,C_{b-1,b})\in\mathbb{Q}(c_1,\dots,c_r,C_{11},C_{12},\dots,C_{b-1,b}),
%\qquad r+1\le k\le m
%$$
%such that the system of equations
%$$
%C_{uv}(c_1, \dots, c_{m})=C_{uv}, \qquad 1\le u\le v\le b,
%$$
%in variables ${c_1, \dots, c_{m}}$ is equivalent to
%$$
%c_k=c_{k}(c_1,\dots,c_r,C_{11},C_{12},\dots,C_{b-1,b}), \qquad r+1\le k\le m.
%$$
%\end{lemma}

\subsection{Цепи и разрезания} \label{correspond}

Мы обобщим идеи из \cite{BeSch,K}, чтобы получить следующее обобщение леммы~\ref{correspondence-lemma}. Основная сложность доказательства --- гарантировать отсутстве \2седловых\3 вершин в цепи. 
%Part (\ref{c-g-2})$\implies$(\ref{c-g-1}) of the following result is probably new, cf.~\cite{BeSch,K}.

\begin{lemma} \label{correspondence-general} Пусть $P$ --- ортогональный многоугольник общего положения с $b$ горизонтальными сторонами, имеющими ориентированные длины ${I_1, \dots, I_b}$ и ${y\mbox{-ко}}$орди\-на\-ты ${U_1, \dots, U_b}$. Тогда следующие $2$ условия эквивалентны\textup{:}
\begin{enumerate}[(1)]
\item \label{c-g-1} многоугольник $P$ можно разрезать на $m$ прямоугольников с отношениями сторон $c_1, \dots, c_m$\textup{;}
%a polygon with sides parallel to the coordinate axes, vertical sides having lengths ${I_1,\dots, I_b}$ и $x$-coordinates ${U_1, \dots, U_b}$, can be tiled by rectangles of ratios ${c_1,\dots,c_n}$;
\item \label{c-g-2} существует плоская замкнутая электрическая цепь с $b$ граничными вершинами, $m$ существенными рёбрами проводимостей ${c_1,\dots,c_m>0}$, входными потенциалами ${U_1, \dots, U_b}$ и входными токами ${I_1, \dots, I_b}$.
\end{enumerate}
\end{lemma}

\begin{remark}
Само по себе условие~(\ref{c-g-2}) не гарантирует существование ортогонального многоугольника с $b$ горизонтальными сторонами, имеющими ориентированные длины ${I_1,\dots,I_b}$ и ${y\mbox{-ко}}$орди\-на\-ты ${U_1,\dots,U_b}$. %Example~\ref{nongeneric}
%\end{remark}
%\begin{remark}
%Example after Theorem~\ref{squaring-polygons} shows that
Лемма \ref{correspondence-general}(\ref{c-g-1})$\implies$(\ref{c-g-2}) не обязательно выполняется для многоугольников
\emph{не общего положения}.
\end{remark}

\begin{figure}[htbp]
\definecolor{qqwwff}{rgb}{0,0.4,1}
\definecolor{ffqqqq}{rgb}{1,0,0}
\definecolor{zzttqq}{rgb}{0.6,0.2,0}
\begin{tikzpicture}[line cap=round,line join=round,>=triangle 45,x=0.7cm,y=0.7cm]
\clip(-1.66,0.8) rectangle (11.1,5.19);
\fill[color=zzttqq,fill=zzttqq,fill opacity=0.1] (-1.5,3) -- (1.5,3) -- (1.5,5) -- (4.5,5) -- (4.5,1) -- (-1.5,1) -- cycle;
\draw [color=zzttqq] (-1.5,3)-- (1.5,3);
\draw [color=zzttqq] (1.5,3)-- (1.5,5);
\draw [color=zzttqq] (1.5,5)-- (4.5,5);
\draw [color=zzttqq] (4.5,5)-- (4.5,1);
\draw [color=zzttqq] (4.5,1)-- (-1.5,1);
\draw [color=zzttqq] (-1.5,1)-- (-1.5,3);
\draw [color=zzttqq] (1.5,3)-- (3.5,3);
\draw [color=zzttqq] (3.5,5)-- (3.5,1);
\draw [line width=1pt,color=ffqqqq] (2,5)-- (4,5);
\draw [line width=1pt,color=ffqqqq] (-1,3)-- (3,3);
\draw [line width=1pt,color=ffqqqq] (-1,1)-- (4,1);
\draw [line width=1pt,dash pattern=on 4pt off 4pt,color=qqwwff] (4,5)-- (4,1);
\draw [line width=1pt,dash pattern=on 4pt off 4pt,color=qqwwff] (2.5,5)-- (2.5,3);
\draw [line width=1pt,dash pattern=on 4pt off 4pt,color=qqwwff] (1,3)-- (1,1);
\draw (5.48,3.7) node[anchor=north west] {$\mathbf{ \rightarrow}$};
\draw [color=zzttqq,fill=zzttqq,fill opacity=0.1] (9,3) circle (1.98);
\draw [line width=1pt,dash pattern=on 4pt off 4pt,color=qqwwff] (9,4.98)-- (9,1.02);
\draw [line width=1pt,dash pattern=on 4pt off 4pt,color=qqwwff] (9,4.98)-- (7.02,3);
\draw [line width=1pt,dash pattern=on 4pt off 4pt,color=qqwwff] (7.02,3)-- (9,1.02);
\fill [color=ffqqqq] (9,4.98) circle (2.0pt);
\fill [color=ffqqqq] (7.02,3) circle (2.0pt);
\fill [color=ffqqqq] (9,1.02) circle (2.0pt);
\end{tikzpicture} 
\caption{построение электрической цепи по данному разрезанию многоугольника на прямоугольники.}
\label{construction}
\end{figure}

\begin{proof}[Доказательство леммы~\ref{correspondence-general}] (\ref{c-g-1})$\implies$(\ref{c-g-2}). %Юnote by ${1,2,\dots,b}$ all the vertical sides of the polygon ${P}$. Denote by ${b+1,\dots,n}$ all the maximal vertical cuts in the tiling.
%Let us construct a graph with vertices ${1,\dots,n}$ in the plane as follows.
Рассмотрим ортогональный многоугольник $P$ общего положения, разрезанный на прямоугольники.
%It is convenient to consider it as a subset of the following auxiliary polygon ${\hat P}$. For each vertical side ${I}$ of ${P}$ take an isosceles triangle ${\hat I}$ with base ${I}$ и a small height such that ${P\cap \hat I=I}$. Set ${\hat P}$ to be the union of ${P}$ и all the triangles ${\hat I}$.

Построим сначала граф требуемой замкнутой электрической цепи (см. рис.~\ref{construction}).
%Take the union of the drawn segments и the horizontal sides of all rectangles ${P_e}$.
%Remove the ''bottom'' horizontal sides of the polygon ${P}$.
Рассмотрим объединение горизонтальных сторон всех прямоугольников разрезания. Это объединение распадается на несколько непересекающихся отрезков, которые мы назовём \emph{горизонтальными разрезами}. Покрасим в красный цвет (жирные линии на рис.~\ref{construction}) все горизонтальные разрезы, кроме небольших окрестностей их концов. Покрасим в синий цвет (пунктирные линии) вертикальную среднюю линию каждого прямоугольника разрезания.

Стянем каждый красный отрезок в точку. При этой операции синее множество \2перейдёт\3 в некоторый граф $G$, а многоугольник $P$ \2перейдёт\3 в некоторый топологический диск $D$
(так как ${y\mbox{-ко}}$ординаты горизонтальных сторон многоугольника $P$ различны, то пересечение каждого красного отрезка с границей $\partial P$ состоит не более, чем из одной связной части).
Обозначим через ${1,\dots,b}$ вершины графа $G$, полученные при стягивании красных отрезков в тех горизонтальных разрезах, которые целиком содержат некоторые стороны многоугольника $P$. Обозначим через ${b+1,\dots,n}$ остальные вершины графа. Ясно, что ${G\subset D}$, ${G\cap \partial D=\{1,\dots,b\}}$ и
что каждая связная компонента графа $G$ содержит граничную вершину. Таким образом, $G$ --- граф некоторой плоской цепи.

Зададим потенциалы, токи и проводимости в данной цепи. Для каждой вершины ${k=1,\dots,n}$ графа $G$ положим потенциал $U_k$ равным ${y\mbox{-ко}}$ординате горизонтального красного отрезка, из которого получилась данная вершина. Для каждого ребра $kl$ графа $G$, полученного из вертикальной средней линии прямоугольника разрезания, положим ток $I_{kl}$ и проводимость $c_{kl}$ равными горизонтальной стороне (с соответствующим знаком) и отношению этого прямоугольника, соответственно. Теперь аксиомы~(C), (I) проверяются непосредственно. Построенная замкнутая электрическая цепь --- искомая.

\begin{figure}[htbp]
\definecolor{qqwwff}{rgb}{0,0.4,1}
\definecolor{zzttqq}{rgb}{0.6,0.2,0}
\begin{tikzpicture}[line cap=round,line join=round,>=triangle 45,x=0.7cm,y=0.7cm]
\clip(0,-1.5) rectangle (15.85,6.21);
\fill[color=zzttqq,fill=zzttqq,fill opacity=0.1] (8.5,3) -- (11.5,3) -- (11.5,5) -- (14.5,5) -- (14.5,1) -- (8.5,1) -- cycle;
\draw [color=zzttqq] (8.5,3)-- (11.5,3);
\draw [color=zzttqq] (11.5,3)-- (11.5,5);
\draw [color=zzttqq] (11.5,5)-- (14.5,5);
\draw [color=zzttqq] (14.5,5)-- (14.5,1);
\draw [color=zzttqq] (14.5,1)-- (8.5,1);
\draw [color=zzttqq] (8.5,1)-- (8.5,3);
\draw [color=zzttqq] (11.5,3)-- (13.5,3);
\draw [color=zzttqq] (13.5,5)-- (13.5,1);
\draw (5.82,3.51) node[anchor=north]
{$\mathbf{ \rightarrow}$};
\draw [color=zzttqq,fill=zzttqq,fill opacity=0.1] (3,3) circle (1.4cm);
\draw [line width=1pt,color=qqwwff] (2.97,5)-- (1,3);
\draw [line width=1pt,color=qqwwff] (2.97,5)-- (3,1);
\draw [line width=1pt,color=qqwwff] (3,1)-- (1,3);
\draw [->] (7.5,0) -- (7.5,6);
\draw [->] (7.5,0) -- (15.5,0);
\draw (6.89,5.44) node[anchor=north] {$U_1$};
\draw (6.89,3.44) node[anchor=north] {$U_2$};
\draw (6.89,1.39) node[anchor=north] {$U_3$};
\draw (8.5,-0.13) node[anchor=north] {$I_{f(3)}$};
\draw (11.51,-0.13) node[anchor=north] {$I_{f(2)}$};
\draw (13.5,-0.13) node[anchor=north] {$I_f$};
\draw (14.51,-0.13) node[anchor=north] {$I_{f(1)}$};
\draw (3.90,3.62) node[anchor=north] {$f(1)$};
\draw (0.55,5.49) node[anchor=north] {$f(2)$};
\draw (0.55,1.77) node[anchor=north] {$f(3)$};
\draw [color=qqwwff](12.60,4.51) node[anchor=north] {$P_{12}$};
\draw [color=qqwwff](11.05,2.58) node[anchor=north] {$P_{23}$};
\draw [color=qqwwff](13.98,3.62) node[anchor=north] {$P_{13}$};
\draw [->] (1.06,4.67) -- (1.65,4.08);
\draw [->] (1.06,1.42) -- (1.67,1.97);
\draw (2.19,3.62) node[anchor=north] {$f$};
\draw[color=zzttqq] (10.60,4.97) node {$P$};
\fill [color=qqwwff] (2.97,5) circle (2.0pt);
\draw[color=qqwwff] (3.17,5.55) node {1};
\fill [color=qqwwff] (3,1) circle (2.0pt);
\draw[color=qqwwff] (3.14,0.49) node {3};
\fill [color=qqwwff] (1,3) circle (2.0pt);
\draw[color=qqwwff] (0.52,3.31) node {2};
\draw [color=black] (7.5,5)-- ++(-1.0pt,0 pt) -- ++(2.0pt,0 pt) ++(-1.0pt,-1.0pt) -- ++(0 pt,2.0pt);
\draw [color=black] (7.5,3)-- ++(-1.0pt,0 pt) -- ++(2.0pt,0 pt) ++(-1.0pt,-1.0pt) -- ++(0 pt,2.0pt);
\draw [color=black] (7.5,1)-- ++(-1.0pt,0 pt) -- ++(2.0pt,0 pt) ++(-1.0pt,-1.0pt) -- ++(0 pt,2.0pt);
\draw [color=black] (8.5,0)-- ++(-1.0pt,0 pt) -- ++(2.0pt,0 pt) ++(-1.0pt,-1.0pt) -- ++(0 pt,2.0pt);
\draw [color=black] (11.5,0)-- ++(-1.0pt,0 pt) -- ++(2.0pt,0 pt) ++(-1.0pt,-1.0pt) -- ++(0 pt,2.0pt);
\draw [color=black] (13.5,0)-- ++(-1.0pt,0 pt) -- ++(2.0pt,0 pt) ++(-1.0pt,-1.0pt) -- ++(0 pt,2.0pt);
\draw [color=black] (14.5,0)-- ++(-1.0pt,0 pt) -- ++(2.0pt,0 pt) ++(-1.0pt,-1.0pt) -- ++(0 pt,2.0pt);
\end{tikzpicture} 
\caption{построение разрезания многоугольника на прямоугольники по данной замкнутой электрической цепи.}
\label{construction2}
\end{figure}

(\ref{c-g-2})$\implies$(\ref{c-g-1}). Рассмотрим замкнутую электрическую цепь, такую, как в условии~(\ref{c-g-2}) леммы.
%Suppose that ${U_1}$, ${\dots}$, ${U_b}$ и ${I_1, \dots, I_b}$ are the incoming voltages и currents, respectively.
%Assume without loss of generality that the minimal $x$-coordinate of the points of the polygon ${P}$ is ${0}$.
%Let ${P}$ be the rectangle with the vertices
%$(0,U_1), (0,U_2), (I_1,U_2), (I_1,U_1)}$.
Построим разрезание многоугольника $P$ на прямоугольники следующим образом.

Пусть $e$ --- ребро цепи. Обозначим через ${e\uparrow}$ (соответственно, ${e\downarrow}$) конец ребра $e$ с более высоким (соответственно, более низким) напряжением. Это определение корректно, так как в условии леммы предполагается, что все ребра являются существенными. \emph{Гранью} назовём связную компоненту дополнения к цепи в единичном диске $D$.
Обозначим через ${e\leftarrow}$ (соответственно, ${e\rightarrow}$) грань, которая остаётся с левой (соответственно, с правой) стороны от ребра $e$
при движении вдоль ребра $e$ от вершины ${e\uparrow}$ к ${e\downarrow}$.

Согласно аксиоме (I) каждой грани $f$ можно приписать действительное число $I_f$ так, чтобы для всех рёбер $kl$ выполнялось равенство ${I_{kl\rightarrow}-I_{kl\leftarrow}=|I_{kl}|}$. Не теряя общности, будем считать, что ${\min_{f}I_f=\min_{(x,y)\in P}x}$, где минимум в левой части берётся по всем граням $f$, примыкающим к границе диска~$D$.

Пусть $P_e$ --- прямоугольник с вершинами ${(I_{e\leftarrow},U_{e\uparrow})}$, ${(I_{e\leftarrow},U_{e\downarrow})}$, ${(I_{e\rightarrow},U_{e\uparrow})}$, ${(I_{e\rightarrow},U_{e\downarrow})}$. Прямоугольники $P_e$, где $e$ пробегает все рёбра цепи, образуют искомое замощение многоугольника $P$ в силу следующих двух утверждений:\\1)~многоугольники $P_e$ покрывают $P$ --- это лемма~\ref{cla};\\2)~многоугольники $P_e$ не перекрываются между собой --- это лемма~\ref{clb}.
\end{proof}
%It suffices to prove the following:
%\begin{enumerate}[(a)]
%\item ${\partial \bigcup P_e\subset\partial P}$;
%\item ${\sum_e\mathop{Area}(P_e)=\mathop{Area}(P)}$.
%\end{enumerate}
%Indeed,

\begin{lemma}\label{cla} ${\bigcup_e P_e=P}$.
\end{lemma}

\begin{proof} Достаточно доказать, что ${\partial \bigcup_e P_e\subset\partial P}$. Так как ${\partial P}$ --- простая замкнутая кривая на плоскости и множество ${\bigcup_e P_e}$ ограничено, то тогда мы получим нужное утверждение.
%­d a vertex ${0}$ lying outside the disc ${D}$ to the network. Join the vertex ${0}$ with the boundary vertices by ${b}$ \emph{new} edges. %Set ${U_0=0}$ и define ${I_{0k}=I_k}$ to be the current flowing inside the network through a boundary vertex ${k}$.

Нам потребуется следующее описание границы ${\partial P}$ (см. рис.~\ref{construction2}).
Граничные вершины цепи разбивают единичную окружность $\partial D$ на $b$ дуг. Стартуя с вершины $b$, будем обходить окружность $\partial D$
против часовой стрелки. Занумеруем дуги в порядке, в котором они встречаются при движении.
Обозначим через $f(v) $ грань, содержащую дугу $v$. Обозначим через $H_v$ отрезок, соединяющий  точки ${(I_{f(v)},U_v)}$ и ${(I_{f({v+1})},U_v)}$.
Обозначим через $V_v$ отрезок, соединяющий  точки ${(I_{f(v)},U_{v-1})}$ и ${(I_{f({v})},U_{v})}$, где мы считаем ${U_0:=U_b}$.
Ясно, что ${\partial P=\bigcup_{v=1}^{b}(H_v\cup V_v)}$.

Возьмём точку ${p\in\partial\bigcup_e P_e}$ \2общего положения\3, скажем, на %the интерьере
горизонтальной стороне \2многоугольника\3 ${\bigcup_e P_e}$.
Точка $p$ обязана принадлежать горизонтальной стороне некоторого прямоугольника разрезания. Пусть, скажем, она принадлежит верхней стороне прямоугольника $P_e$. Положим ${v:=e\uparrow}$.
%Юnote by ${\mathop{St}v}$ the set of edges containing the vertex ${v}$.
%	note by ${S}$ the set of all the rectangles ${P_d}$ such that ${v}$ is an endpoint of the edge ${d}$.

Проведём горизонтальную прямую $H$ через верхнюю сторону прямоугольника $P_e$. %Consider все прямоугольники $P_d$ таким образом, что $d\in\mathop {С-} v$.
Мы говорим, что прямоугольник $P_d$ --- \emph{примыкающий}, если вершина $v$ является концом ребра $d$.
Примыкающие прямоугольники примыкают к прямой $H$ сверху или снизу.

Предположим сначала, что вершина $v$ не граничная. Простая индукция показывает, что к каждой точке прямой $H$ (за исключением некоторого конечного множества) примыкает одинаковое число примыкающих прямоугольников $P_d$ сверху и снизу. Так как прямоугольник $P_e$ примыкает к точке $p$ снизу и точка $p$ выбрана общего положения, то некоторый примыкающий прямоугольник $P_d$ примыкает к ней сверху. Таким образом, $p$ принадлежит внутренности ${\mathop {Int} P_e\cup P_f\subset \mathop {Int} \bigcup_e P_e}$ --- противоречие.

Значит, $v$ --- граничная вершина. Аналогично рассуждению выше, к каждой точке фигуры ${H-H_v}$ (за исключением некоторого конечного множества) примыкает одинаковое число примыкающих прямоугольников $P_d$ сверху и снизу. Следовательно, ${p\in H_v}$. Таким образом, ${p\in\partial P}$, что и требовалось.
\end{proof}
%Assume that say ${I_{e\leftarrow}<x<I_{e\rightarrow}}$, ${y=U_{e\uparrow}}$. Let ${u=e\uparrow}$. Let us show that ${u}$ is a boundary vertex и ${I_{u \leftarrow}<x<I_{u\rightarrow}}$.

%Let ${e=e_1,\dots,e_p}$ be all the edges for which ${u=e_k\uparrow}$ и ${e_{p+1},\dots,e_q}$ be all the edges for which ${u=e_k\downarrow}$. Consider \emph{positive} segments ${[I_{e_k\leftarrow},I_{e_k\rightarrow}]}$, ${1\le k\le p}$, и \emph{negative} segments ${[I_{e_k\rightarrow},I_{e_k\leftarrow}]}$, ${p+1\le k\le q}$.

%Assume that ${u}$ is nonboundary. тогда clearly any point of the line ${\mathbb{R}}$ is covered by the same number of positive и negative segments. In particular, there is a negative segment ${[I_{e_k\rightarrow},I_{e_k\leftarrow}]}$ containing ${x}$. тогда the point ${(x,y)}$ belongs to the horizontal side of the rectangle ${P_{e_k}}$ и lies in the interior of ${P_e\cup P_{e_k}}$, a contradiction.
%So ${u}$ is a boundary vertex.

%Assume ${x\notin [I_{u \leftarrow},I_{u\rightarrow}]}$. Still any point of the set ${\mathbb{R}-[I_{u \leftarrow},I_{u\rightarrow}]}$ is covered by the same number of positive и negative segments. Thus we again get a contradiction. Hence ${I_{u \leftarrow}<x<I_{u\rightarrow}}$. So ${(x,y)}$ belongs to a horizontal side of ${P}$.

%Analogously one can prove that a point ${(x,y)}$ in the interior of a vertical side of the polygon ${\bigcup_e P_e}$ necessarilylies in a vertical side of ${P}$. This proves the claim.

\begin{lemma}\label{clb}
${\sum_e\Area(P_e)=\Area(P)}$.
\end{lemma}

\begin{proof} Утверждение немедленно следует из леммы~\ref{energy}, так как площади прямоугольников выражаются формулами ${\Area(P_{kl})=(U_k-U_l)I_{kl}}$ и ${\Area(P)=\sum_{1\le u\le b}U_uI_u}$.
\end{proof}

\subsection{Обратные задачи}\label{basicposreal}

Сначала приведём доказательство леммы~\ref{PosReal}, следующее идеям работы~\cite{Cauer}. Обобщения на случай ${b>2}$ см. в \cite{Foster}.% и ссылки в ней.

\begin{proof}[Доказательство леммы~\ref{PosReal}]
(\ref{preserve})$\implies$(\ref{roots}). Действительно, если ${\Real z\le 0}$, то ${\Real C(z)=-\Real C(-z)\le 0}$ и поэтому ${C(z)\ne 1}$.
%So (\ref{roots}) holds.
%\cite[Lemma 4]{FLR} Let ${r ( z ) = C(z) + C ( - z )}$. We wish to show
%that ${r ( z )}$ is identically zero. To begin, note that ${r(z)}$ takes the imaginary axis into itself.
%Thus ${r ( i y ) = C(iy) + C ( - i y ) = C(iy) +  C(\overline{iy}) = C ( i y ) + \overline{C(iy)} = 0}$. Since the
%numerator of ${r ( z )}$ has an infinite number of roots, the numerator is identically zero,
%and hence so is ${r(z)}$. Finally, suppose ${C(z) = x}$ has a root with real part less than or
%equal to zero. если this root has real part less than zero, тогда~(\ref{preserve}) is treated. если this root is
%on the imaginary axis, тогда by continuity we can pick ${z}$ sufficiently close to the root but
%with negative real part so that ${C(z)}$ is close enough to ${x}$ so that ${\Real (C(z)) > 0}$, also
%violating~(\ref{preserve}).

(\ref{roots})$\implies$(\ref{preserve}). %Since ${C(z)}$ is odd it follows that all roots of ${C(z)=-1}$ have negative real parts.
Рассмотрим уравнение ${C(z)=w}$. Будем двигать ${w}$ непрерывно в полуплоскости ${\Real w>0}$. Корни нашего уравнения не могут пересекать прямую ${\Real z=0}$ (так как если ${\Real z=0}$, то ${\Real C(z)=0}$ для нечётной функции ${C(z)\in\mathbb{R}(z)}$). Поэтому для всех ${w}$ в полуплоскости ${\Real w>0}$ все корни уравнения ${C(z)=w}$ находятся в полуплоскости ${\Real z>0}$. Так как ${C(z)}$ нечётна, то то же самое верно и для полуплоскостей ${\Real w<0}$, ${\Real z<0}$. Поэтому утверждение\7(\ref{preserve}) справедливо.

(\ref{preserve})$\implies$(\ref{derivatives}). %Since ${C(z)}$ is odd it follows that ${\Real z<0}$ implies ${\Real C(z)<0}$.
Предположим, что ${C(z)=0}$. %, where ${z\in \mathbb{C}}$.
Тогда ${\Real z=0}$ (так как ${\Real z>0\implies\Real C(z)>0}$ и ${\Real z<0\implies\Real C(z)=-\Real C(-z)<0}$). Так как утверждение~(\ref{preserve}) и обратное к нему справедливы в окрестности точки ${z}$, то ${C'(z)>0}$.
%A simple limiting argument proves the same for ${z=\infty}$.

(\ref{derivatives})$\implies$(\ref{factorization})
%Assume for simplicity that ${C(\infty)\ne0}$.
Пусть ${z_1,\dots,z_m}$ --- нули функции ${C(z)}$. Так как ${C'(z_k)>0}$, то эти нули являются простыми. Поэтому ${C(z)}$ имеет не более, чем ${m}$ полюсов. Нули ${z_1,\dots,z_m}$ разделяют проективную прямую ${\Real z=0}$ на ${m}$ \2отрезков\3. Так как ${C'(z_k)>0}$, то для достаточно малых ${\epsilon>0}$ справедливо ${i\cdot C(z_k-i\epsilon)>0}$ и ${i\cdot C(z_k+i\epsilon)<0}$. По теореме о промежуточном значении каждый \2отрезок\3 содержит полюс функции ${C(z)}$. Поэтому все ${m}$ полюсов функции ${C(z)}$ принадлежат проективной прямой ${\Real z=0}$ и чередуются с нулями ${z_1,\dots,z_m}$. Значит, утверждение\7(\ref{factorization}) справедливо.
%Thus ${C(z)}$ satisfies condition~(\ref{factorization}).

(\ref{factorization})$\implies$(\ref{continued}). %The \emph{height} of a rational function is the sum of the degrees of the nominator и the denominator.
Обозначим через ${\mathop{ht}C(z)}$ сумму степеней числителя и знаменателя несократимой дроби ${C(z)}$.
Доказательство проведём индукцией по ${\mathop{ht}C(z)}$. Если ${\mathop{ht}C(z)=1}$, то доказывать нечего. Теперь предположим, что ${n\ge 1}$ и, скажем, ${b_n\ne 0}$ в условии~(\ref{factorization}).
%Let
%$$
%C(z)=\frac{\prod\limits_{k=1}^n
%(z^2+b_k^2)}{d_1z\prod\limits_{k=1}^n(z^2+a_k^2)}.
%$$
%(The other alternatives from condition (\ref{factorization}) are considered similarly.)

Положим ${r(z):=1/(C(z)-d_1z)}$ и ${q(z):=1/C(z)}$. Докажем, что ${r(z)}$ удовлетворяет условию~(\ref{derivatives}).
Действительно, ясно, что ${\lim_{z\to\infty}r(z)\ne 0}$ и ${r(z)}$ нечётна.
Нулями функции ${r(z)}$ являются числа ${\pm ib_1,\dots,\pm ib_n}$.
Для всех ${l=1,\dots,n}$
$$r'(\pm ib_l)=q'(\pm ib_l)=\frac{2}{d_1(a_l^2-b_l^2)}\prod\limits_{k\ne l}\frac{b_k^2-b_l^2}{a_k^2-b_l^2}>0
$$
в силу условия ${a_1>b_1>a_2>\dots>b_n\ge 0}$.

Следовательно, по лемме~\ref{PosReal}(\ref{derivatives})$\implies$(\ref{factorization}) функция
${r(z)}$ также удовлетворяет условию~(\ref{factorization}). С другой стороны, ${\mathop{ht}r(z)<\mathop{ht}C(z)}$. По предположению индукции, ${r(z)}$ удовлетворяет условию~(\ref{continued}).
Поэтому
${C(z)=d_1z+1/r(z)}$
также удовлетворяет условию~(\ref{continued}).

%(\ref{continued})$\implies$(\ref{series-parallel-network}). The construction of the required network is straightforward, cf.~\S\ref{prooflfrs}.

(\ref{continued})$\implies$(\ref{preserve}).
Это доказывается индукцией по ${m}$.
%This follows from Lemma~\ref{electrical}(\ref{el-5}).
\end{proof}

\begin{proof}[Доказательство теоремы~\ref{th-foster}]
(\ref{th-foster-1})$\implies$(\ref{th-foster-2}). Это немедленно следует из определений (см.~\S\ref{alternating}).

(\ref{th-foster-2})$\implies$(\ref{th-foster-3}). Это следует из леммы~\ref{electrical}(\ref{el-2}) и~(\ref{el-5}).

(\ref{th-foster-3})$\implies$(\ref{th-foster-1}).
Сначала предположим, что ${\lim_{z\to\infty}C(z)\ne0}$. Тогда по
лемме~\ref{PosReal}(\ref{preserve})$\implies$(\ref{continued})
условие~(\ref{continued}) выполнено.
Тогда требуемая последовательно-параллельная цепь с проводимостями рёбер ${d_1z}$, ${1/d_2z}$, ${d_3z}$, ${1/d_4z}$, \dots, ${(d_mz)^{(-1)^m}}$ сразу легко строится (ср.~\S\ref{prooflfrs}).

Теперь предположим, что ${\lim_{z\to\infty}C(z)=0}$. Применим
лемму~\ref{PosReal}(\ref{preserve})$\implies$(\ref{continued})
к функции ${1/C(z)}$. Тогда функция ${C(z)}$ удовлетворяет условию~(\ref{continued}) леммы~\ref{PosReal} с той лишь разницей, что ${d_1=0}$.
Снова применяя конструкцию предыдущего абзаца, приходим к требуемой цепи.
\end{proof}

\begin{remark}\label{rationality} Если ${C(z)\in\mathbb{Q}(z)}$, то цепь в условии~(\ref{th-foster-2}) теоремы~\ref{th-foster} можно считать
последовательно-параллельной, а каждое число ${d_j}$ рациональным. %Indeed, by The Euclidean algorithm it follows that all ${d_j\in\mathbb{Q}}$ in expression~(\ref{continued}) of Lemma~\ref{PosReal}.
\end{remark}

\begin{proof}[Доказательство следствия~\ref{series-parallel}] %Let ${C(c_1,c_2)}$ be the conductance of the given network. тогда по лемме~\ref{electrical}(2) we get ${C(c_1,c_2)=\sqrt{c_1c_2}\, C(z,1/z)}$, where ${z=\sqrt{c_1/c_2}}$.
По условию~(\ref{th2-2}) теоремы~\ref{th2} имеем ${C(c_1,c_2)=\sqrt{c_1c_2}\, C(z)}$, где ${C(z)=C(z,1/z)}$ и ${z=\sqrt{c_1/c_2}}$.
По теореме~\ref{th-foster}(\ref{th-foster-3})$\implies$(\ref{th-foster-1})
и замечанию~\ref{rationality}
существует последовательно-параллельная цепь с проводимостью ${C(z)}$, проводимость каждого ребра $j$ которой равна
либо ${d_jz}$, либо ${1/d_jz}$ для некоторых ${d_1,\dots,d_m\in\mathbb{Q}}$. Заменим каждое ребро на соответствующую последовательно-параллельную подцепь, чтобы получить новую цепь с проводимостями рёбер ${z}$ и ${1/z}$. Умножая проводимость каждого ребра на ${\sqrt{c_1c_2}}$, получаем требуемую цепь.
%по лемме~\ref{PosReal}(\ref{preserve})$\implies$(\ref{continued}) it satisfies condition (\ref{continued}).
%Therefore, say, for ${m}$ even и ${C(0)=0}$,
%$$
%C(c_1,c_2)=d_1c_1+\cfrac1{\cfrac{d_2}{c_2} +\cfrac1{d_3 c_1+\dots+\cfrac{d_m}{ c_2}}}.
%$$
%All the numbers ${d_k\in \mathbb{Q}}$ by the Euclidean algorithm. Now the required series-parallel network is constructed analogously to the proof теоремы~\ref{lfrs}(\ref{lfrs-3})$\implies$(\ref{lfrs-1}).
\end{proof}

\section{Доказательство вариаций}\label{finalproofs}

\subsection{Доказательство теоремы \ref{squaring-polygons}} \label{proof-general}

%The proof теоремы~\ref{squaring-polygons} follows the idea of \cite[\S8.1]{K} и uses the following results of \cite{CGV}.

%\begin{theorem} \label{CGV5} \textup{\cite[Theorem~5]{CGV}} The set ${-\bar\Omega_b}$ is the set of all possible responses of planar networks with ${b}$ boundary vertices и positive edge conductances.
%\end{theorem}

%\begin{theorem} \label{CGV3} \textup{\cite[Theorem~3]{CGV}} Let ${G}$ be a network with a given response ${C_{uv}}$ having minimal number of edges. Consider the response of ${G}$ as a matrix-valued function in edge conductances ${c_{kl}}$. тогда this function is injective.
%\end{theorem}

\begin{proof}[Доказательство теоремы~\ref{squaring-polygons}] (\ref{s-p-1})$\implies$(\ref{s-p-2}).
Предположим, что многоугольник ${P}$ разрезан на квадраты.
%rectangles of ratios ${c_1, \dots,c_n}$.
По лемме~\ref{correspondence-general} существует плоская электрическая цепь с единичными проводимостями рёбер,
%$c_1, \dots,c_n}$,
входными потенциалами ${U_1}$, \dots, ${U_b}$ и входными токами ${I_1}$, \dots, ${I_b}$. Пусть
$C_{uv}$ --- реакция этой цепи. Тогда ${I_v=\sum C_{uv} U_u}$. По лемме~\ref{responseclaim}(\ref{rc-1}) все элементы матрицы ${C_{uv}}$ рациональны. По теореме~\ref{CGV5} получаем ${C_{uv}\in \Omega_b}$.
%\end{proof}

(\ref{s-p-2})$\implies$(\ref{s-p-1}). Пусть ${C_{uv}\in \Omega_b}$ --- матрица с рациональными элементами, такая, что ${I_v=\sum C_{uv} U_u}$. По теореме~\ref{CGV5} существует плоская электрическая цепь с реакцией ${C_{uv}}$. Рассмотрим минимальную цепь с таким свойством. По теореме~\ref{inverse}
проводимость каждого ребра минимальной цепи рациональна.
%Докажем that all the conductances ${c_{kl}\in\mathbb{Q}}$. Indeed, otherwise perform a nontrivial automorphism of the field ${\mathbb Q(c_{12},\dots,c_{n,n-1})}$. The obtained network has the same graph и the same response по лемме~\ref{responseclaim}(1). This contradicts to minimality of the number of edges in the circuit по теореме~\ref{CGV3}. Thus all ${c_{kl}\in\mathbb Q}$.
Пусть входные потенциалы равны ${U_1}$, \dots, ${U_b}$. Тогда входные токи равны ${I_1}$, \dots, ${I_b}$. Удалим все несущественные рёбра цепи. По лемме~\ref{correspondence-general} многоугольник ${P}$ можно разрезать на прямоугольники с рациональными отношениями сторон, а следовательно, и на квадраты.
\end{proof}

%\begin{remark} The converse theorem would follow from Lemma~\ref{correspondence-general}, \cite[Theorem~5]{CGV} и the following conjecture: \emph{if the response of a network with real edge conductances is rational тогда существует a network with rational edge conductances и the same response}.
%\end{remark}

\begin{corollary} \label{mostgeneral} \textup{(лемм~\ref{responseclaim}, \ref{correspondence-general} и теоремы~\ref{CGV5})} Пусть ${P}$ --- ортогональный многоугольник общего положения с ${b}$ горизонтальными сторонами, имеющими ориентированные длины ${I_1,\dots,I_b}$ и ${y\mbox{-ко}}$ординаты ${U_1,\dots,U_b}$. Предположим, что многоугольник ${P}$ можно разрезать на прямоугольники с отношениями сторон ${c_1}$, ${\dots}$, ${c_n}$. Тогда существует функция ${C_{uv}(z_1,\dots,z_n)}$, удовлетворяющая условиям~\textup{(\ref{rc-1}), (\ref{rc-2}),} и \textup{(\ref{rc-5})} леммы~\ref{responseclaim}, такая, что ${C_{uv}(c_1,\dots,c_n)\in \Omega_b}$ и ${I_v=\sum_{1\le u\le b} C_{uv}(c_1,\dots,c_n)U_u}$ для всех ${v=1,\dots,b}$.
\end{corollary}

\subsection{Доказательство теоремы~\ref{polygon-square}}

\begin{proof}[Доказательство теоремы~\ref{polygon-square}] ${\impliedby}$. Это верно, так как любой многоугольник с рациональными координатами вершин можно разрезать на квадраты.

$\implies$. Пусть ${P}$ разрезан на прямоугольники с отношениями сторон $c$ и ${1/c}$. Докажем аналогично доказательству теоремы~\ref{lfrs}(\ref{lfrs-1})$\implies$(\ref{lfrs-2}), что все корни минимального многочлена числа $c$ имеют положительную действительную часть. Тогда теорема~\ref{polygon-square} будет следовать из теоремы~\ref{lfrs}(\ref{lfrs-2})$\implies$(\ref{lfrs-1}).

Рассмотрим замкнутую электрическую цепь, описанную в лемме \ref{correspondence-general}.
Заменим каждое ребро с проводимостью $c$ (соответственно, ${1/c}$) в цепи на ребро с проводимостью ${z\in\mathbb{C}}$ (соответственно, ${1/z}$). Пусть ${C_{uv}(z)}$ --- реакция полученной цепи. Рассмотрим функцию \emph{диссипации энергии} ${E(z)=\sum_{1\le u,v\le b} C_{uv}(z)U_u U_v}$. Так как каждое ${U_u\in\mathbb{Q}}$, то по лемме~\ref{responseclaim}(\ref{rc-1}) получаем, что ${E(z)\in\mathbb{Q}(z)}$. По лемме~ \ref{responseclaim}(\ref{rc-2}) получаем, что
${E(z)}$ нечётна.
Ясно, что ${E(c)=\sum_{1\le u\le b} I_u U_u=\Area(P)}$. Поэтому ${E(c)\in\mathbb{Q}}$ и ${E(c)>0}$.
%We have ${E(c)>0}$ по лемме \ref{responseclaim}(4) (in fact ${E(c)}$ is the area of ${P}$).

Так как ${E(z)\in\mathbb{Q}(z)}$, ${E(z)\not\equiv\mbox{const}}$ и ${E(c)\in\mathbb{Q}}$, то $c$ --- корень многочлена с целыми коэффициентами. Возьмём такой многочлен минимальной степени. Пусть ${z}$ --- любой другой его комплексный корень. Тогда ${E(z)=E(c)>0}$.

Докажем, что ${\Real z>0}$.
Сначала предположим, что ${\Real z<0}$. Тогда по лемме \ref{responseclaim}(\ref{rc-5}) имеем ${0\le\Real E(-z)=-\Real E(z)<0}$ --- противоречие. Простой предельный переход показывает, что предположение ${\Real z=0}$ также приводит к противоречию.
Поэтому ${\Real z>0}$.
\end{proof}

%\begin{proof}[Доказательство теоремы~\ref{rectangling-polygons}]
%Suppose that a polygon ${P}$ is tiled by rectangles of ratios ${c_1,\dots,c_n}$.
%Consider the electrical circuit given по лемме~\ref{correspondence-general}.
%For each ${k=1,\dots,n}$ replace each edge of conductance ${c_k}$ in the network by an edge of conductance ${z_k}$. Пусть ${C_{uv}(z_1,\dots,z_n)}$ be the response of the obtained network. тогда Theorem~\ref{rectangling-polygons}(1)--(4)
%follows from Lemma~\ref{responseclaim} и \cite[Theorem~5]{CGV}.
%\end{proof}

\subsection{Доказательство теоремы~\ref{detection}}

\begin{proof}[Доказательство теоремы~\ref{detection}] (\ref{det-1})$\implies$(\ref{det-2}). Это следует из леммы~\ref{responseclaim}(\ref{rc-5}) и замечания~\ref{strict}.

(\ref{det-2})$\implies$(\ref{det-1}). При ${b=2}$ доказывать нечего. Предположим, что ${b=3}$. Пусть ${\delta>0}$ --- достаточно малое число, ${r_{uv}:=-\Real C_{uv}-\delta}$, ${m_{uv}:=-\Imaginary C_{uv}}$. По почти положительной определённости получаем, что
${\Real }${\tiny ${\left(\begin{matrix}
C_{11} & C_{12} \\
C_{21} & C_{22}
\end{matrix}\right)}$}
положительно определена. Поэтому
{\tiny ${\left(\begin{matrix}
r_{31}+r_{12} & -r_{12} \\
-r_{12}       & r_{12}+r_{23}
\end{matrix}\right)}$}
положительно определена для достаточно малого ${\delta}$. Следовательно, ${r_{12}+r_{23}, r_{31}+r_{12}, r_{12}r_{23}+r_{23}r_{31}+r_{31}r_{12}> 0}$. Аналогично ${r_{23}+r_{31}>0}$. Поэтому по крайней мере два числа из $r_{12}$, $r_{23}$, $r_{31}$ являются положительными.

Если ${r_{12}, r_{23}, r_{31}>0}$, то требуемая цепь --- треугольник с вершинами ${1}$, ${2}$, ${3}$ и проводимостями рёбер ${c_{kl}=r_{kl}+im_{kl}+\delta}$.

Теперь предположим, что ровно одно из чисел $r_{12}$, $r_{23}$, $r_{31}$, скажем, ${r_{31}}$, неположительно. Возьмём достаточно большое число ${M}$ и положим ${\Delta_M:=r_{12}r_{23}+r_{23}r_{31}+r_{31}r_{12}+iM(r_{23}+r_{12})}$.
Требуемая цепь является полным графом с вершинами ${1}$, ${2}$, ${3}$, ${4}$ и проводимостями рёбер
$$
\begin{tabular}{llll}
$c_{12}=$& ${im_{12} +\delta}$,
\quad & ${c_{14}=}$& ${\Delta_M/r_{23}}$,\\
$c_{23}=$& ${im_{23} +\delta}$,
\quad & ${c_{34}=}$& ${\Delta_M/r_{12}}$,\\
$c_{31}=$& ${im_{31} +\delta-iM}$,
\quad & ${c_{24}=}$& ${\Delta_M/(r_{31}+iM)}$.
\end{tabular}
$$
Ясно, что при ${M^2>(r_{12}r_{23}+r_{23}r_{31}+r_{31}r_{12})|r_{31}|/(r_{23}+r_{12})}$ все ${\Real  c_{kl}> 0}$.

Покажем при помощи преобразования данной электрической цепи, что она имеет реакцию ${C_{uv}}$. Заменим \2звезду\3, состоящую из рёбер ${14}$, ${24}$ и ${34}$, на треугольник, состоящий из трёх новых рёбер с проводимостями ${c'_{12}=r_{12}}$, ${c'_{23}=r_{23}}$ и ${c'_{31}=r_{31}+iM}$. Это \emph{преобразование звезда-треугольник} не изменяет реакцию цепи \cite[p. 12]{K}. Полученная цепь имеет три пары кратных рёбер. Поэтому она имеет ту же реакцию, что и треугольник с проводимостями рёбер ${r_{12}+im_{12}+\delta}$, ${r_{23}+im_{23}+\delta}$, ${r_{31}+im_{31}+\delta}$. Значит, реакция цепи равна ${C_{uv}}$.конкурсом Мёбиуса
\end{proof}

\section*{Благодарности} Авторы благодарны А.~Акопяну, М.~Хьюстону и Б.~Джорджу за полезные обсуждения и Ф. Шарову за перевод на русский язык.
М.~Скопенков частично поддержан грантом ИНТАС 06-1000014-6277,
грантом российского фонда фундаментальных исследований
06-01-72551-НЦНИЛ-а, 
%РФФИ-01-00648-a,
%President of the Russian Federation grant NSh-4578.2006.1,
%Agency for Education и Science grant RNP-2.1.1.7988,
фондом поддержки молодых учёных \2Конкурс Мёбиуса\3 и фондом Эйлера.

%\section*{References}

\bibliographystyle{elsarticle-num}